\documentclass[11pt,feqno]{article}
\usepackage{amssymb,bbm,amsmath,amsthm,amscd}
\usepackage{vmargin}
\usepackage[dvipsnames]{xcolor}
\usepackage{array}
\usepackage{soul}
 
 \usepackage{graphicx}
 \usepackage{subfigure}
\DeclareGraphicsExtensions{.pdf,.eps}

\catcode`\@=11
\newbox\bz@
\newdimen\bdimz@
\def\linethrough#1{\setbox\bz@=\hbox{#1}
\bdimz@=\ht\bz@ \divide\bdimz@ by 5 \advance\bdimz@ by -\dp\bz@ \ht\bz@=\bdimz@
\leavevmode\hbox{$\overline{\overline{\box\bz@}}$\relax}}
\catcode`\@=12

\def\nbZ{{\mathchoice {\hbox{$\sf\textstyle Z\kern-0.4em Z$}}
{\hbox{$\sf\textstyle Z\kern-0.4em Z$}}
{\hbox{$\sf\scriptstyle Z\kern-0.3em Z$}}
{\hbox{$\sf\scriptscriptstyle Z\kern-0.2em Z$}}}}

\newcounter{compteur}[subsection]

\usepackage[colorlinks=true]{hyperref}

\newtheorem{theorem}{Theorem}[section]
\newtheorem{lemma}[theorem]{Lemma}

\newtheorem{definition}[theorem]{Definition}
\newtheorem{hypothesis}[theorem]{Hypothesis}
\newtheorem{remark}[theorem]{Remark}

\usepackage{macros}

\makeatletter 

\@addtoreset{equation}{section}

\begin{document}

\author{Samuele Rubino \thanks{Departamento EDAN \& IMUS, Universidad de Sevilla, Spain. {\tt samuele@us.es}}}
\title{Numerical analysis of a projection-based stabilized POD-ROM for incompressible flows}
\maketitle

\begin{abstract}
In this paper, we propose a new stabilized projection-based POD-ROM for the numerical
simulation of incompressible flows. The new method draws inspiration from successful
numerical stabilization techniques used in the context of Finite Element (FE) methods, 
such as Local Projection Stabilization (LPS). In particular, the new LPS-ROM is
a velocity-pressure ROM that uses pressure modes as well to compute the reduced order
pressure, needed for instance in the computation of relevant quantities, such as drag and
lift forces on bodies in the flow. The new LPS-ROM circumvents the standard discrete inf-sup condition 
for the POD velocity-pressure spaces, whose fulfillment can be rather expensive
in realistic applications in Computational Fluid Dynamics (CFD). Also, the velocity modes does not have to be
neither strongly nor weakly divergence-free, which allows to use snapshots generated
for instance with penalty or projection-based stabilized methods. The numerical analysis
of the fully Navier--Stokes discretization for the new LPS-ROM is presented, by mainly
deriving the corresponding error estimates. Numerical studies are performed to discuss
the accuracy and performance of the new LPS-ROM on a two-dimensional laminar unsteady flow past
a circular obstacle.
\end{abstract}


\medskip


\section{Introduction}\label{sec:Intro}

Reduced Order Models (ROM) have been applied to numerical design in modern engineering as a tool that is wide-spreading in the scientific community in the recent years in order to solve complex realistic 
multi-parameters, multi-physics and multi-scale problems. Among the most popular ROM approaches, Proper Orthogonal Decomposition (POD) strategy provides optimal (from the energetic point of view) bases or modes to represent the dynamics from a given database 
(snapshots) obtained by a full order system. Onto these reduced bases, a Galerkin projection of the governing equations can be employed to obtain a low-order dynamical system for the bases coefficients. This has led researchers to apply POD-ROM to a variety of physical and engineering problems, including Computational Fluid Dynamics (CFD) 
problems in order to model the Navier--Stokes Equations (NSE), see e.g. \cite{Baiges13, Gunzburger06, Iollo04, Barone12, Beran03, Wang12}.

\medskip

In this context, POD velocity modes are usually assumed to be at least weakly divergence-free. To be this assumption true, the POD velocity modes should be generated, for instance, by a Full Order Model (FOM) which consists in a NSE space discretization using inf-sup stable Finite Element (FE) for the velocity-pressure pair. In this way, the contribution of the pressure formally drops out from the ROM, which thus only approximates the velocity field through the POD velocity modes.
Despite the appealing computational efficiency of only velocity ROM, there exist however important settings in which the pressure should be somehow considered. Indeed, the pressure is needed in many CFD applications, e.g. in the computation of relevant physical quantities, such as drag and lift forces on bodies in the flow, and for incompressible shear flows, as the mixing layer or the wake flow \cite{NPM05}, where neglecting it may lead to large amplitude errors. On the other hand, note also that the weakly divergence-free property does not hold for many popular discretizations of the NSE. This is the case, for instance, of using equal order FE for the velocity-pressure pair, for which a suitable numerical stabilization becomes essential to circumvent the violation of the standard discrete inf-sup condition, as considered in this work. Furthermore, the pressure approximation allows the computation of the residual associated to the strong form of the NSE, often needed in stabilized discretizations (cf. \cite{BergmannIollo09}). Altogether, all these reasons pushed us to propose and fully analyze a new robust and stable ROM that directly incorporates an approximation of the pressure, driven by numerical stabilization motivations, and recovers it correctly, avoiding spurious pressure oscillations.

\medskip

The new method draws inspiration from successful
numerical stabilization techniques used in the context of FE methods, 
such as Local Projection Stabilization (LPS) methods (cf. \cite{IMAJNA,ARCME}). In particular, the new LPS-ROM is
a coupled velocity-pressure ROM that uses pressure modes as well to compute the reduced order
pressure. In order to avoid pressure instabilities sources, the new LPS-ROM circumvents the standard discrete inf-sup condition 
for the POD velocity-pressure spaces, whose fulfillment can be rather expensive
in realistic applications in CFD, see for instance \cite{Rozza15,RozzaStabile18}, where an offline strategy based on the supremizer enrichment of the reduced velocity space has been proposed and applied in the POD context, adapted from the Reduced Basis (RB) method framework. Also, with respect to other proposals existing in the current ROM literature that provides velocity-pressure approximations, the velocity modes does not have to be
neither strongly nor weakly divergence-free for the new LPS-ROM, which allows to use snapshots generated
for instance with penalty or projection-based stabilized methods. This is not the case, for instance, of ROM based on a pressure Poisson equation approach (see, for instance, the first two methods investigated in \cite{IliescuJohn14} and also \cite{RozzaStabile18}), for which the velocity snapshots, and hence the POD velocity modes must be at least weakly divergence-free. This requirement also holds for the last method proposed and investigated in \cite{IliescuJohn14}, which uses a residual-based stabilization mechanism in order to overcome a possible violation of the discrete inf-sup condition in the ROM framework by considering a decoupled approach for the reduced velocity-pressure pair.

\medskip

The main contribution of the present paper has been to perform a stability and convergence analysis of the arising fully discrete LPS-ROM applied to the unsteady incompressible NSE, by mainly deriving the proof of a rigorous error estimate that considers all contributions: the spatial discretization error (due to the FE discretization), the temporal discretization error (due to the backward Euler method), and the POD truncation error. To the best of our knowledge, the LPS-ROM introduced is novel, and the numerical analysis for a stabilization-motivated ROM to take into account the violation of the discrete inf-sup condition cannot be found in the literature so far. Indeed, on the one hand a thorough numerical analysis has been recently performed for a stabilization-motivated ROM accounting for instabilities due to convection-dominated phenomena \cite{ACR17}, which only uses velocity POD modes. On the other hand, only few numerical investigations of stabilization-motivated ROM can be found in the literature (cf. \cite{BergmannIollo09} and last method in \cite{IliescuJohn14}). In particular, in \cite{IliescuJohn14} the authors showed that adding a Pressure-Stabilizing Petrov--Galerkin (PSPG) term to account for stabilizing the violated discrete inf-sup condition and recover the reduced pressure provided more efficient and accurate results with respect to ROM based on a pressure Poisson equation approach, which would also require an {\em ad-hoc} treatment of the pressure boundary conditions. Parallel and independently to the current paper, a velocity-pressure ROM has been very recently proposed and analyzed in \cite{IliescuSchneier19}, which however relies on an artificial compression method to compute the reduced velocity-pressure approximations for which the pressure must be initialized.

\medskip

Note that the detailed numerical analysis corroborated with numerical studies makes apparent an interesting link between the number of POD velocity-pressure modes used in the LPS-ROM and the angle $\theta$ between the space spanned by the divergence of the POD velocity modes and the POD pressure space. Indeed, for the numerical example proposed, where the same numerical stabilization technique used for the ROM is initially applied also to the FOM (LPS-FOM) to generate the snapshots, so that these latter are not weakly divergence-free, we have found that for small values of $r$, which is common in practice, the saturation constant $\alpha=cos(\theta)$ \cite{ChaconDominguez00} is rather small, and this allows to ease the convergence order reduction due to the violation of the discrete inf-sup stability condition.

\medskip

Numerical studies performed on a two-dimensional laminar unsteady flow past a circular obstacle have also been used to assess the accuracy and efficiency of the new LPS-ROM. Despite the fact that the discrete inf-sup condition is not fulfilled by the new LPS-ROM, using a small equal number of POD velocity-pressure modes already provides accurate approximations, close to the LPS-FOM results, and theoretical considerations suggested by the numerical analysis are recovered in practice.

\medskip

The outline of the paper is as follows: In Section \ref{sec:VF}, we introduce the model problem and its continuous
variational formulation for time-dependent NSE. In Section \ref{sec:FOM}, we consider the FE-LPS full order discretization used to generate the snapshots for the online phase. Section \ref{sec:POD-ROM} briefly describes the POD methodology and introduce the formulation of the new LPS-ROM for the NSE. The stability and error analysis for the full discretization (FE in space and backward Euler in time) of the new model is presented in Section \ref{sec:NA}. The numerical investigation of the new method is proposed in Section \ref{sec:NumStud} for the simulation of a two-dimensional flow past a circular obstacle, in order to test on the one hand some theoretical predictions of the performed numerical analysis and to show on the other hand the accuracy and efficiency of the proposed method in preventing spurious pressure instabilities due to the violation of the discrete inf-sup condition for the reduced system. Finally, Section \ref{sec:Concl} presents the main conclusions of this work and ongoing research directions.

\section{Time-dependent NSE: model problem and variational formulation}\label{sec:VF}

We introduce an Initial-Boundary Value Problem (IBVP) for the incompressible evolution NSE. 
For the sake of simplicity, we just impose the homogeneous Dirichlet boundary condition on the whole boundary.

\medskip

Let $[0,T]$ be the time interval and $\Om$ a bounded polyhedral domain in $\mathbb{R}^{d}$, $d=2$ or $3$, with a Lipschitz-continuous boundary $\Ga=\partial\Om$. The transient NSE for an incompressible fluid are given by:

\medskip
\hspace{1cm} {\em Find $\uv:\Om\times (0,T)\longrightarrow\mathbb{R}^{d}$ and $p:\Om\times (0,T)\longrightarrow\mathbb{R}$ such that:}
\BEQ\label{eq:uNS}
\left \{
\begin{array}{rcll}
\partial_{t}\uv + (\uv \cdot\nabla)\uv - \nu \Delta\uv + \nabla p&=&\fv & \qmbx{in} \Om\times (0,T),\\
\div \uv &=&0 & \qmbx{in} \Om\times (0,T),\\
\uv &=& \bf{0} & \qmbx{on} \Ga\times (0,T),\\
\uv(\xv,0)&=&\uv_{0}(\xv) & \qmbx{in} \Om.
\end{array}
\right .
\EEQ
The unknowns are the velocity $\uv(\xv,t)$ and the pressure $p(\xv,t)$ of the incompressible fluid. The data are the source term $\fv(\xv,t)$, which represents a body force per mass unit (typically  the gravity), the kinematic viscosity $\nu$ of the fluid, which is a positive constant, and the initial velocity $\uv_{0}(\xv)$. 

\medskip

To define the weak formulation of problem (\ref{eq:uNS}), we need to introduce some useful notations for functional spaces \cite{Brezis11}. We consider the Sobolev spaces $H^{s}(\Om)$, $s\in \mathbb{R}$, $L^{p}(\Om)$ and $W^{m,p}(\Om)$, $m\in \mathbb{N}$, $1\leq p\leq\infty$. We shall use the following notation for vector-valued Sobolev spaces: ${\bf H}^{s}$, ${\bf L}^{p}$ and ${\bf W}^{m,p}$ respectively shall denote $[H^{s}(\Om)]^{d}$, $[L^{p}(\Om)]^{d}$ and $[W^{m,p}(\Om)]^{d}$ (similarly for tensor spaces of dimension $d\times d$). Also, the parabolic Bochner function spaces $L^{p}(0,T;X)$ and $L^{p}(0,T;{\bf X})$, where $X$ (${\bf X}$) stands for a scalar (vector-valued) Sobolev space, shall be denoted by $L^{p}(X)$ and $L^{p}({\bf X})$, respectively. In order to give a variational formulation of problem (\ref{eq:uNS}), let us consider the velocity space:
$$
\Xv={\bf H}_{0}^{1}=[H_{0}^{1}(\Om)]^{d}=\left\{\vv\in [H^{1}(\Om)]^{d} : \vv={\bf 0} \text{ on } \Ga\right\}. 
$$
This is a closed linear subspace of ${\bf H}^{1}$ and thus a Hilbert space endowed with the ${\bf H}^{1}$-norm. 
Thanks to Poincar\'e inequality, the ${\bf H}^{1}$-norm is equivalent on ${\bf H}_{0}^{1}$ to the norm $\nor{\vv}{{\bf H}_{0}^{1}}=\nor{\nabla\vv}{{\bf L}^{2}}$. Also, let us consider the pressure space:
$$
Q=L_{0}^{2}(\Om)=\left\{q\in L^{2}(\Om) : \int_{\Om}q\,d\xv=0\right\}. 
$$
Note that the null mean condition for the pressure is usually introduced in order to fix the constant the pressure is determined up through the formulation. However, this condition could be relaxed by adding a suitable penalty term to the variational formulation of \eqref{eq:uNS} that allows to simply take $Q=L^{2}(\Om)$ as pressure space.

\medskip

We shall thus consider the following variational formulation of \eqref{eq:uNS}:

\medskip
\hspace{1cm} {\em Given $\fv\in L^{2}({\bf H}^{-1})$, find $\uv: (0,T)\longrightarrow\Xv$, $p: (0,T)\longrightarrow Q$ such that}
\BEQ\label{eq:fvNS} 
\left \{ 
\begin{array}{rcll}
\disp\frac{d}{dt}(\uv,\vv) + b(\uv,\uv,\vv) + \nu (\nabla\uv,\nabla\vv) - (p,\div\vv) &=&\langle \fv,\vv\rangle & \forall\vv\in\Xv,\quad \text{in } \mathcal{D}'(0,T),\\
(\div \uv,q) + \sigma(p,q)&=&0 & \forall q \in Q, \quad \text{a.e. in } (0,T),\\
\uv(0)&=&\uv_{0},
\end{array} 
\right .
\EEQ
where $( \cdot,\cdot )$ stands for the $L^2$-inner product in $\Om$, $\langle \cdot,\cdot \rangle$ stands for the duality pairing between $\Xv$ and its dual $\Xv^{'}={\bf H}^{-1}$, and $\mathcal{D}'(0,T)$ is the space of distributions in $(0,T)$. The trilinear form
$b$ is given by: {\it for $\uv,\, \vv,\,\wv \in \Xv$}
\BEQ\label{eq:formb}
b(\uv,\vv,\wv)={1 \over 2} \left[   ( \uv \cdot \nabla\, \vv, \wv) - ( \uv \cdot \nabla\,  \wv, \vv)\right].
\EEQ
The term with factor $\sigma$ denotes the penalty term that permits to fix the constant the pressure is determined up through the formulation, for a small positive value of $\sigma$ (e.g., $\sigma = \mathcal{O}(10^{-6})$), and thus $Q=L^{2}(\Om)$ in \eqref{eq:fvNS} and hereafter. Note that in the continuous model problem \eqref{eq:uNS} the incompressibility condition is no more satisfied exactly, and thus we are going to search for a velocity field approximation that does not have to be neither strongly nor weakly divergence-free.

\section{Finite element full order model}\label{sec:FOM}

In order to give a FE approximation of \eqref{eq:fvNS}, let $\{{\cal T}_{h}\}_{h>0}$ be a family of affine-equivalent, conforming (i.e., without hanging nodes) and regular triangulations of $\overline{\Om}$, 
formed by triangles or quadrilaterals ($d=2$), tetrahedra or hexahedra ($d=3$). For any mesh cell $K \in {\cal T}_{h}$,
its diameter will be denoted by $h_K$ and $h = \max_{K \in {\cal T}_{h}} h_K$. We consider $\Xv_{h}\subset\Xv$, $Q_{h}\subset Q$ being suitable FE spaces for velocity and pressure, respectively. 
The FE approximation of \eqref{eq:fvNS} can be written as follows:
\medskip

\hspace{1cm}{\em Find $(\uhv,p_h):(0,T)\longrightarrow \Xv_{h}\times Q_{h}$ such that}
\BEQ\label{eq:FEapprox} 
\left \{ 
\begin{array}{rcll}
\disp\frac{d}{dt}(\uhv,\vhv) + b(\uhv,\uhv,\vhv) + \nu (\nabla\uhv,\nabla\vhv) - (p_h,\div\vhv) &=&\langle \fv,\vhv\rangle & \text{in } \mathcal{D}'(0,T),\\
(\div \uhv,q_h) + \sigma(p_h,q_h)&=&0 & \text{a.e. in } (0,T),\\
\uhv(0)&=&\uv_{0h},
\end{array} 
\right .
\EEQ
for any $(\vhv,q_h)\in \Xv_{h}\times Q_{h}$, and the initial condition $\uv_{0h}$ is some stable approximation to $\uv_{0}$ in $L^2$-norm belonging to $\Xv_{h}$.

\medskip

In order to circumvent the standard discrete inf-sup condition and thus use equal order interpolation for velocity and pressure, and also provide an extra-control on the high frequencies components of the pressure gradient that could lead to unstable discretizations, we introduce a filtered pressure stabilizing term of high-order. In this way, the considered FE method falls into the class of Local Projection Stabilization (LPS) methods (cf. \cite{IMAJNA,ARCME}). The stabilization effect is achieved by adding a least-square term that give a weighted control on the fluctuations of the pressure gradient, based upon a specific locally stable projection or interpolation operator on a continuous buffer space. This provides an efficient discretization with a reduced computational cost that keeps the same high-order accuracy with respect to standard projection-stabilized methods.

\medskip

To describe this approach, we define hereafter the specific choice of FE spaces done both for the numerical analysis and practical computations in the present work. Given an integer $l\geq 2$ and a mesh cell $K \in {\cal T}_{h}$, denote by $\mathbb{R}^{l}(K)$ either $\mathbb{P}^{l}(K)$ (i.e., the space of Lagrange polynomials of degree $\leq l$, defined on $K$), if the grids are formed by triangles ($d=2$) or tetrahedra ($d=3$), or $\mathbb{Q}^{l}(K)$ (i.e., the space of Lagrange polynomials of degree $\leq l$ on each variable, defined on $K$), if the family of triangulations is formed by quadrilaterals ($d=2$) or hexahedra ($d=3$). We consider the following FE spaces for the velocity:
\BEQ\label{eq:defVelSp}
\left \{ 
\begin{array}{lll}
&Y_{h}^{l}=V_{h}^{l}(\Om) = \{v_{h} \in C^{0}(\overline{\Om}) : v_{h}|_{K} \in \mathbb{R}^{l}(K),\,\forall K\in {\cal T}_{h}\},&\\ \\
&\Yhv^{l}=[Y_{h}^{l}]^{d} = \{\vhv \in [C^{0}(\overline{\Om})]^{d} : \vv_{h}|_{K}\in [\mathbb{R}^{l}(K)]^{d},\,\forall K\in {\cal T}_{h}\},&\\ \\
&\Xhv = \Yhv^{l}\cap {\bf H}_{0}^{1}.&
\end{array} 
\right .
\EEQ
Hereafter, $\Yhv^{l}$ (resp., $Y_{h}^{l}$) will constitute the discrete foreground vector-valued (resp., scalar) spaces in which we will work on. Also, $Q_h=Y_{h}^{l}\subset L^{2}$, since we use equal order FE. We define the scalar product:
$$
(\cdot,\cdot)_{\tau}:{\bf L}^{2}(\Om)\times {\bf L}^{2}(\Om) \to \mathbb{R},\quad
(\vv,\wv)_{\tau} = \sum_{K\in{\cal T}_{h}}\tau_{K}(\vv,\wv)_{K},
$$
and its associated norm:
$$
\nor{\vv}{\tau}=(\vv,\vv)_{\tau}^{1/2},
$$
where for any $K\in{\cal T}_{h}$, $\tau_{K}$ is in general a positive local stabilization parameter.

\medskip

The considered FE-FOM is given by: 
\medskip

\hspace{1cm}{\em Find $(\uhv,p_h):(0,T)\longrightarrow \Xv_{h}\times Q_{h}$ such that}
\BEQ\label{eq:FOM} 
\left \{ 
\begin{array}{rcll}
\disp\frac{d}{dt}(\uhv,\vhv) + b(\uhv,\uhv,\vhv) + \nu (\nabla\uhv,\nabla\vhv) - (p_h,\div\vhv) &=&\langle \fv,\vhv\rangle & \text{in } \mathcal{D}'(0,T),\\
(\div \uhv,q_h) + (\Pi_{h}^{*}(\nabla p_h),\Pi_{h}^{*}(\nabla q_h))_{\tau} + \sigma(p_h,q_h)&=&0 & \text{a.e. in } (0,T),\\
\uhv(0)&=&\uv_{0h},
\end{array} 
\right .
\EEQ
for any $(\vhv,q_h)\in \Xv_{h}\times Q_{h}$, where $\Pi_{h}^{*}=Id-\Pi_{h}$ is the ``fluctuation operator'', being $Id$ the identity operator, and $\Pi_h$ some locally stable (in $L^2$-norm) projection or interpolation operator from ${\bf L}^2(\Om)$ on the foreground vector-valued space $\Yhv^{l-1}$ (also called ``buffer space'' in this context), satisfying optimal error estimates (cf. \cite{Ciarlet02}).
In practical implementations, we choose $\Pi_{h}$ as a Scott--Zhang-like \cite{ScottZhang90} linear interpolation operator in the space $\Yhv^{l-1}$ (see \cite{Chacon13}, Sect. 4 for its construction), implemented by the software FreeFEM \cite{Hecht12}.

\medskip

To state the full space-time discretization of the unsteady LPS-FOM \eqref{eq:FOM}, consider a positive integer number $N$ and define $\Delta t = T/N$, $t_{n}=n\Delta t$, $n=0,1,\ldots,N$. We compute the approximations $\uhv^{n}$, $\ph^{n}$ to $\uv^{n}=\uv(\cdot,t_{n})$ and $p^{n}=p(\cdot,t_{n})$ by using, for simplicity of the analysis, a backward Euler scheme:
\begin{itemize}
\item {\bf Initialization.} Set: $\uhv^{0}=\uohv.$
\item {\bf Iteration.} For $n=0,1,\ldots,N-1$: {\em Given $\uhv^{n}\in\Xhv$, find $(\uhv^{n+1},\ph^{n+1})\in\Xhv\times Q_h$ such that:}
\BEQ\label{eq:discTempFOM} 
\left \{ 
\begin{array}{rcl}
\left(\disp\frac{\uhv^{n+1}-\uhv^{n}}{\Delta t},\vhv\right)+b(\uhv^{n+1},\uhv^{n+1},\vhv)+\nu(\nabla\uhv^{n+1},\nabla\vhv)& &\\
-(\ph^{n+1},\div\vhv) &=& \langle \fv^{n+1}, \vhv \rangle,\\ \\
(\div\uhv^{n+1},\qh)+(\Pi_{h}^{*}(\nabla \ph^{n+1}),\Pi_{h}^{*}(\nabla q_h))_{\tau}+\sigma(\ph^{n+1},\qh) &=& 0,
\end{array} 
\right .
\EEQ
for any $(\vhv,\qh) \in \Xhv\times Q_h$.
\end{itemize}


\section{Proper orthogonal decomposition reduced order model}\label{sec:POD-ROM}

We briefly describe the POD method, following \cite{KunischVolkwein01}, and apply it to the projection-based stabilized FOM \eqref{eq:discTempFOM}. 

\medskip

Let us consider the ensembles of velocity snapshots $\chi^{v}=\text{span}\left\{\uhv^{1},\ldots,\uhv^{N}\right\}$ and pressure snapshots $\chi^{p}=\text{span}\left\{\ph^{1},\ldots,\ph^{N}\right\}$, given by the FE solutions to \eqref{eq:discTempFOM} at time $t_{n}$, $n=1,\ldots,N$. The POD method seeks low-dimensional bases $\left\{\boldsymbol{\varphi}_{1},\ldots,\boldsymbol{\varphi}_{r_v}\right\}$ and $\left\{\psi_{1},\ldots,\psi_{r_p}\right\}$ in real Hilbert spaces $\mathcal{H}_{v}$, $\mathcal{H}_{p}$ that optimally approximate the velocity and pressure snapshots in the following sense:
\BEQ\label{eq:PODmethVel}
\min\Delta t\sum_{n=1}^{N}\left\| \uhv^{n} - \sum_{i=1}^{r_v}\left(\uhv^{n},\boldsymbol{\varphi}_{i}\right)_{\mathcal{H}_{v}}\boldsymbol{\varphi}_{i} \right\|_{\mathcal{H}_{v}}^2,
\EEQ
and
\BEQ\label{eq:PODmethPres}
\min\Delta t\sum_{n=1}^{N}\left\| \ph^{n} - \sum_{i=1}^{r_p}\left(\ph^{n},\psi_{i}\right)_{\mathcal{H}_{p}}\psi_{i} \right\|_{\mathcal{H}_{p}}^2,
\EEQ
subject to the conditions $\left(\boldsymbol{\varphi}_{j},\boldsymbol{\varphi}_{i}\right)_{\mathcal{H}_{v}}=\delta_{ij}$, $1\leq i,j \leq r_v$ and $\left(\psi_{j},\psi_{i}\right)_{\mathcal{H}_{p}}=\delta_{ij}$, $1\leq i,j \leq r_p$, where $\delta_{ij}$ is the Kronecker delta. To solve the optimization problems \eqref{eq:PODmethVel}-\eqref{eq:PODmethPres}, one can respectively consider the eigenvalue problems:
\BEQ\label{eq:eigenVel}
K^{v}{\bf a}_{i}=\lambda_i{\bf a}_{i},\text{ for } 1,\ldots,r_{v},
\EEQ
and
\BEQ\label{eq:eigenPres}
K^{p}{\bf b}_{i}=\gamma_i{\bf b}_{i},\text{ for } 1,\ldots,r_{p},
\EEQ 
where $K^{v}, K^{p}\in \mathbb{R}^{N\times N}$ are the velocity, pressure snapshots correlation matrices, respectively with entries: 
$$
K_{mn}^{v}=\Delta t\left(\uhv^{n},\uhv^{m}\right)_{\mathcal{H}_{v}},\text{ for } m,n=1,\ldots,N,
$$ 
and
$$
K_{mn}^{p}=\Delta t\left(\ph^{n},\ph^{m}\right)_{\mathcal{H}_{p}},\text{ for } m,n=1,\ldots,N,
$$ 
${\bf a}_{i}, {\bf b}_{i}$ are the $i$-th eigenvector, and $\lambda_{i},\gamma_{i}$ are the associated eigenvalues. The eigenvalues are positive and sorted in descending order: $\lambda_{1}\geq\ldots\geq\lambda_{r_v}>0$ and $\gamma_{1}\geq\ldots\geq\gamma_{r_p}>0$. It can be shown that the solutions of \eqref{eq:PODmethVel}-\eqref{eq:PODmethPres}, i.e. the POD velocity-pressure bases functions, are respectively given by:
\BEQ\label{eq:PODbasisVel}
\boldsymbol{\varphi}_{i}(\cdot)=\frac{1}{\sqrt{\lambda_i}}\sqrt{\Delta t}\sum_{n=1}^{N}({\bf a}_{i})_{n}\uhv^{n},\quad 1\leq i\leq r_v,
\EEQ
and
\BEQ\label{eq:PODbasisPres}
\psi_{i}(\cdot)=\frac{1}{\sqrt{\gamma_i}}\sqrt{\Delta t}\sum_{n=1}^{N}({\bf b}_{i})_{n}\ph^{n},\quad 1\leq i\leq r_p,
\EEQ
where $({\bf a}_{i})_{n},({\bf b}_{i})_{n}$ are the $n$-th components of the eigenvectors ${\bf a}_{i}, {\bf b}_{i}$, respectively. It can also be shown that the following POD projection error formulas hold \cite{Holmes96, KunischVolkwein01}:
\BEQ\label{eq:PODerrVel}
\Delta t\sum_{n=1}^{N}\left\| \uhv^{n} - \sum_{i=1}^{r_v}\left(\uhv^{n},\boldsymbol{\varphi}_{i}\right)_{\mathcal{H}_{v}}\boldsymbol{\varphi}_{i} \right\|_{\mathcal{H}_{v}}^{2} = \sum_{i=r_{v}+1}^{M_{v}}\lambda_{i},
\EEQ
and
\BEQ\label{eq:PODerrPres}
\Delta t\sum_{n=1}^{N}\left\| \ph^{n} - \sum_{i=1}^{r_p}\left(\ph^{n},\psi_{i}\right)_{\mathcal{H}_{p}}\psi_{i} \right\|_{\mathcal{H}_{p}}^{2} = \sum_{i=r_{p}+1}^{M_{p}}\gamma_{i},
\EEQ
where $M_{v},M_{p}$ are the rank of $\chi^v$ and $\chi^p$, respectively. Although $\mathcal{H}_{v},\mathcal{H}_{p}$ can be any real Hilbert spaces, in what follows we consider $\mathcal{H}_{v}={\bf L}^{2}$ and $\mathcal{H}_{p}=L^{2}$. Also, we are going to take the same number of velocity and pressure POD bases functions, i.e. $r_{v}=r_{p}=r$ in what follows. Thus, we expect that also the POD velocity-pressure spaces will not satisfy the standard discrete inf-sup condition, and the POD-Reduced Order Model (POD-ROM) we are going to consider must circumvent it. To overcome this restriction, we draw inspiration from the FOM \eqref{eq:discTempFOM} in order to construct the new projection-based stabilized POD-ROM.

\medskip

We respectively consider the following velocity and pressure spaces for the POD setting:
$$
\Xv_{r}=\text{span}\left\{\boldsymbol{\varphi}_{1},\ldots,\boldsymbol{\varphi}_{r}\right\}\subset \Xhv,
$$
and
$$
Q_{r}=\text{span}\left\{\psi_{1},\ldots,\psi_{r}\right\}\subset Q_h.
$$

\begin{remark}
Since, as shown in \eqref{eq:PODbasisVel}, the POD velocity modes are linear combinations of the velocity snapshots, the POD velocity modes satisfy the boundary conditions in \eqref{eq:uNS}. This is because of the particular choice we have made at the beginning to work with homogeneous Dirichlet boundary conditions. In general, one has to manipulate the velocity snapshots set. This is the case, for instance, of steady-state non-homogeneous Dirichlet boundary conditions, for which is preferable to consider a proper lift in order to generate POD velocity modes for the lifted velocity snapshots, satisfying homogeneous Dirichlet boundary conditions. This would lead to work with centered-trajectory method in the POD-ROM setting \cite{IliescuJohn15}.
\end{remark}

The standard Galerkin projection-based POD-ROM uses both Galerkin truncation and Galerkin projection.
The former yields an approximation of the velocity and pressure fields by a linear combination of the corresponding truncated POD basis:
\BEQ\label{eq:PODsolVel}
\uv(\xv,t)\approx \uv_{r}(\xv,t)=\sum_{i=1}^{r}a_{i}(t)\boldsymbol{\varphi}_{i}(\xv),
\EEQ
and
\BEQ\label{eq:PODsolVel}
p(\xv,t)\approx p_{r}(\xv,t)=\sum_{i=1}^{r}b_{i}(t)\psi_{i}(\xv),
\EEQ
where $\left\{a_{i}(t)\right\}_{i=1}^{r}$ and $\left\{b_{i}(t)\right\}_{i=1}^{r}$
are the sought time-varying coefficients representing the POD-Galerkin velocity and pressure trajectories. Note that $r<<\mathcal{N}$, where $\mathcal{N}$ denotes the number of degrees of freedom (d.o.f.) of the equal order FE velocities-pressure in FOM \eqref{eq:discTempFOM}. Replacing the velocity-pressure FE pair $(\uhv,\ph)$ with $(\uv_{r},p_r)$ in the FE approximation \eqref{eq:discTempFOM} and projecting the resulted equations onto the POD product space $(\Xv_{r},Q_r)$ using the POD basis $\left(\left\{\boldsymbol{\varphi}_{i}\right\}_{i=1}^{r},\left\{\psi_{i}\right\}_{i=1}^{r}\right)$, the full space-time discretization of the new projection-based stabilized POD-ROM reads as: 
\begin{itemize}
\item {\bf Initialization.} Set: $\uv_{r}^{0}=\disp\sum_{i=1}^{r}(\uohv,\boldsymbol{\varphi}_{i})\boldsymbol{\varphi}_{i}.$
\item {\bf Iteration.} For $n=0,1,\ldots,N-1$: {\em Given $\uv_{r}^{n}\in\Xv_{r}$, find $(\uv_{r}^{n+1},p_{r}^{n+1})\in\Xv_{r}\times Q_r$ such that:}
\BEQ\label{eq:discTempROM} 
\left \{ 
\begin{array}{rcl}
\left(\disp\frac{\uv_{r}^{n+1}-\uv_{r}^{n}}{\Delta t},\boldsymbol{\varphi}\right)+b(\uv_{r}^{n+1},\uv_{r}^{n+1},\boldsymbol{\varphi})+\nu(\nabla\uv_{r}^{n+1},\nabla\boldsymbol{\varphi})& &\\
-(p_{r}^{n+1},\div\boldsymbol{\varphi}) &=& \langle \fv^{n+1}, \boldsymbol{\varphi} \rangle,\\ \\
(\div\uv_{r}^{n+1},\psi)+(\Pi_{h}^{*}(\nabla p_{r}^{n+1}),\Pi_{h}^{*}(\nabla \psi))_{\tau}+\sigma(p_{r}^{n+1},\psi) &=& 0,
\end{array} 
\right .
\EEQ
for any $(\boldsymbol{\varphi},\psi) \in \Xv_{r}\times Q_r$.
\end{itemize}

\medskip

An alternative time discretization could be given by the semi-implicit Euler method, where the trilinear form in \eqref{eq:discTempROM} is discretized by $b(\uv_{r}^{n},\uv_{r}^{n+1},\boldsymbol{\varphi})$. Note that considering a semi-implicit time discretization of the new projection-based stabilized POD-ROM is less costly from the computational point of view with respect to a fully implicit one, which yields a nonlinear algebraic system of equations to be solved. However, the numerical analysis will be performed in detail for the more technical case of the fully implicit time discretization given by \eqref{eq:discTempROM}.

\section{Analysis of the projection-based stabilized POD-ROM}\label{sec:NA}

In this section, we perform the numerical analysis of the proposed unsteady POD-ROM \eqref{eq:discTempROM}, which we will
call in the sequel LPS-ROM.

\subsection{Technical background}\label{subsec:TB}

This section provides some technical results that are required for the numerical analysis. Throughout the paper, we shall denote by $C$ a positive constant that may vary from a line to another, but which is always independent of the FE mesh size $h$, the FE velocity-pressure equal interpolation order $l$, the time step $\Delta t$, and the velocity, pressure eigenvalues $\lambda_i$, $\gamma_i$.

\medskip

\begin{lemma}[See Lemma 13 in \cite{Layton08}]\label{lm:formb}
For any function $\uv,\vv,\wv\in\Xv$, the skew-symmetric trilinear form $b(\cdot,\cdot,\cdot)$ satisfies:
\BEQ\label{eq:formb}
b(\uv,\vv,\vv)=0,
\EEQ 
\BEQ\label{eq:formb1}
b(\uv,\vv,\wv)\leq C\nor{\nabla\uv}{{\bf L}^{2}}\nor{\nabla\vv}{{\bf L}^{2}}\nor{\nabla\wv}{{\bf L}^{2}}.
\EEQ
\end{lemma}

\begin{definition}
Let $X$ be a Hilbert space and $Y$, $Z$ two finite-dimensional subspaces of $X$ with intersection reduced to the zero function. The pair of finite-dimensional spaces $(Y,Z)$ is called to satisfy the saturation property if there exists a positive constant $C$ such that:
\BEQ\label{eq:SatProp}
\nor{y}{X}+\nor{z}{X}\leq C\nor{y+z}{X}\quad \forall y\in Y, z\in Z.
\EEQ
\end{definition}

Thus, the saturation property can be viewed as an inverse triangular inequality.

\begin{lemma}\label{lm:AngleCond}
The saturation property is equivalent to the existence of a constant $\alpha<1$ such that:
\BEQ\label{eq:AngleCond}
|(y,z)_X|\leq \alpha\nor{y}{X}\nor{z}{X}\quad \forall y\in Y, z\in Z.
\EEQ
\end{lemma}

Actually, we may take $\alpha=1-2/C^2$ (see Remark 2 in \cite{ChaconDominguez00}), and in the sequel we will call $\alpha$ the saturation constant. Then, we can interpret the saturation property in the sense that the angle between spaces $Y$ and $Z$, defined by:
\BEQ\label{eq:AngleDef}
\theta=\arccos\left(\sup_{y\in Y\backslash \{0\},\, z\in Z\backslash \{0\}}\frac{(y,z)_{X}}{\nor{y}{X}\nor{z}{X}}\right),
\EEQ
is uniformly bounded from below by a positive angle, and $\alpha=\cos(\theta)$. 

\begin{remark}
Lemma \ref{lm:AngleCond} will be essential in Theorem \ref{th:PODEE} to bound the error term co\-ming from the continuity equation, which cannot be removed by using the standard Stokes projection since the reduced velocity-pressure spaces in the LPS-ROM \eqref{eq:discTempROM} proposed in this work violate the standard discrete inf-sup condition. 

\medskip

In the numerical studies performed in Section \ref{sec:NumStud} we will observe that the saturation constant $\alpha$, for the chosen numerical setup, starts with small values for small $r$ and seems to experience a flattening effect around $10^{-2}$ when adding more POD modes. 

\medskip

Note that the argument of saturation property has been used in \cite{ChaconDominguez00} to develop a stabilized post-processing of the Galerkin FE solution of convection-dominated flows and very recently extended \cite{ACR19} to POD-ROM approximations to propose a cure for instabilities due to advection-dominance in POD solution to advection-diffusion-reaction equations \cite{ESAIM}.
\end{remark}

To ensure error estimates in Theorem \ref{th:PODEE} (main result of the present paper), we make the following regularity assumption on the continuous solution: 
\begin{hypothesis}\label{hp:ContRegularity}
We assume that the continuous solution $(\uv,p)$ of the unsteady NSE \eqref{eq:fvNS} has augmented regularity, i.e. $(\uv,p)\in L^{\infty}({\bf H}^{s+1})\times L^{\infty}(H^{s})$, $s\geq 1$, such that $\partial_t^2 \uv\in L^2({\bf L}^2)$.
\end{hypothesis}

For the subsequent numerical analysis, we need the following technical hypothesis on the stabilization parameters $\tau_K$:
\begin{hypothesis}\label{hp:stabCoef}
The stabilization parameters $\tau_K$ satisfy the following condition:
\BEQ
\tau_K\leq C h_{K}^{2} \quad \forall K\in{\cal T}_h.
\EEQ
\end{hypothesis}
\begin{remark}
The question whether the stabilization parameters should depend on the
spatial resolution of the underlying FE space, or on the number of POD basis functions
used has been addressed in \cite{IliescuJohn15}, by means of numerical analysis arguments. In that work,
numerical investigations using both definitions suggested that the one based on estimates
from the underlying FE discretization provides a better suppression of numerical oscillations, 
and thus guarantees a more effective numerical stabilization. For this reason, we
make here assumption \ref{hp:stabCoef} on the stabilization parameters, which is also essential for the
subsequent numerical analysis.
\end{remark}
%

\begin{definition}
Let $\vv_r$ be the $L^2$-orthogonal projection of $\uv$ on the reduced velocity space $\Xv_r$:
\BEQ\label{eq:L2projVel}
(\uv-\vv_r,\boldsymbol{\varphi})=0 \quad \forall \boldsymbol{\varphi}\in \Xv_r,
\EEQ
and $z_r$ be the $L^2$-orthogonal projection of $p$ on the pressure space $Q_r$:
\BEQ\label{eq:L2projPres}
(p-z_r,\psi)=0 \quad \forall \psi\in Q_r.
\EEQ
\end{definition}

We have the following error estimates for $\vv_r$ and $z_r$ (see \cite{IliescuWang14}, Lemma 3.3):
\begin{lemma}[$L^2$-projection error estimates]
\BEQ\label{eq:PODlmL2vel}
\Delta t \sum_{n=1}^{N}\nor{\uv^{n}-\vv_r^n}{{\bf L}^2}^2\leq C\left(h^{2s}+\Delta t ^{2} + \sum_{i=r+1}^{M_v}\lambda_{i}\right),
\EEQ
\BEQ\label{eq:PODlmL2pres}
\Delta t \sum_{n=1}^{N}\nor{p^{n}-z_r^n}{L^2}^2\leq C\left(h^{2s}+\Delta t ^{2} + \sum_{i=r+1}^{M_p}\gamma_{i}\right).
\EEQ
\end{lemma}

\begin{lemma}[$H^1$-projection error estimates]\label{lm:PODlmH1}
\BEQ\label{eq:PODlmH01vel}
\Delta t \sum_{n=1}^{N}\nor{\nabla(\uv^{n}-\vv_r^n)}{{\bf L}^2}^2\leq C\left[\left(1+\nor{S^v_r}{2}\right)\left(h^{2s}+\Delta t ^{2}\right) + \sum_{i=r+1}^{M_v}\nor{\nabla\boldsymbol{\varphi}_{i}}{{\bf L}^2}^{2}\lambda_{i}\right],
\EEQ
with $\nor{S^v_r}{2}$ denoting the $2$-norm of the stiffness velocity matrix with entries $[S^v_r]_{ij}=(\nabla\boldsymbol{\varphi}_{j},\nabla\boldsymbol{\varphi}_{i})$, $i,j=1,\ldots,r$.
\BEQ\label{eq:PODlmH1pres}
\Delta t \sum_{n=1}^{N}\nor{\nabla(p^{n}-z_r^n)}{{\bf L}^2}^2\leq C\left[\left(1+\nor{S^p_r}{2}\right)\left(h^{2s}+\Delta t ^{2}\right) + \sum_{i=r+1}^{M_p}\nor{\nabla\psi_{i}}{{\bf L}^2}^{2}\gamma_{i}\right],
\EEQ
with $\nor{S^p_r}{2}$ denoting the $2$-norm of the stiffness pressure matrix with entries $[S^p_r]_{ij}=(\nabla\psi_{j},\nabla\psi_{i})$, $i,j=1,\ldots,r$.
\end{lemma}
The appearance of $\nor{S^v_r}{2}$, $\nor{S^p_r}{2}$ in Lemma \ref{lm:PODlmH1} comes from the use of the POD inverse estimates (see \cite{KunischVolkwein01}, Lemma 2):
\begin{eqnarray}
\nor{\nabla\boldsymbol{\varphi}}{{\bf L}^2}&\leq& \nor{S^v_r}{2}^{1/2}\nor{\boldsymbol{\varphi}}{{\bf L}^2} \quad \forall \boldsymbol{\varphi}\in\Xv_r,\\
\nor{\nabla\psi}{{\bf L}^2}&\leq& \nor{S^p_r}{2}^{1/2}\nor{\psi}{L^2} \quad \forall \psi\in Q_r.
\end{eqnarray}

\subsection{Existence and stability results for LPS-ROM}\label{subsec:STAB}

We have the following existence and unconditional stability result for the LPS-ROM \eqref{eq:discTempROM}:
\begin{theorem}\label{th:stabEst}
Problem \eqref{eq:discTempROM} admits a solution that satisfies the following bound:
\begin{eqnarray}\label{eq:stabEst}
&&\nor{\uv_{r}^{k}}{{\bf L}^{2}}^{2}+\sum_{n=0}^{N-1}\nor{\uv_{r}^{n+1}-\uv_{r}^{n}}{{\bf L}^{2}}^{2} +\Delta t\sum_{n=0}^{N-1}\left(\nu\nor{\nabla\uv_{r}^{n+1}}{{\bf L}^{2}}^{2} 
+ \nor{\Pi_{h}^{*}(\nabla p_{r}^{n+1})}{\tau}^{2}+\sigma\nor{p_{r}^{n+1}}{L^{2}}^{2}\right) \nonumber \\
&\leq&
\nor{\uv_{r}^{0}}{{\bf L}^{2}}^{2}+\frac{4\Delta t}{\nu}\sum_{n=0}^{N-1}\nor{\fv^{n+1}}{{\bf H}^{-1}}^{2},
\end{eqnarray}
for $k=0,\ldots,N$.
\end{theorem}
{\bf Proof.}
Problem \eqref{eq:discTempROM} can be written as:
\BEQ\label{eq:discTempROMsteady} 
\left \{ 
\begin{array}{rcl}
b(\uv_{r}^{n+1},\uv_{r}^{n+1},\boldsymbol{\varphi})+\widetilde{a}(\uv_{r}^{n+1},\boldsymbol{\varphi})
-(p_{r}^{n+1},\div\boldsymbol{\varphi}) &=& \langle \widetilde{\fv}^{n+1}, \boldsymbol{\varphi} \rangle,\\ \\
(\div\uv_{r}^{n+1},\psi)+(\Pi_{h}^{*}(\nabla p_{r}^{n+1}),\Pi_{h}^{*}(\nabla \psi))_{\tau}+\sigma(p_{r}^{n+1},\psi) &=& 0,
\end{array} 
\right .
\EEQ
for any $(\boldsymbol{\varphi},\psi)\in \Xv_{r}\times Q_r$, where $\widetilde{a}(\uv_{r}^{n+1},\boldsymbol{\varphi})=\Delta t^{-1}(\uv_{r}^{n+1},\boldsymbol{\varphi})+\nu (\nabla\uv_{r}^{n+1},\nabla\boldsymbol{\varphi})$ and $ \widetilde{\fv}^{n+1}=\fv^{n+1}+\Delta t^{-1}(\uv_{r}^{n},\boldsymbol{\varphi})$. 
This problem fits into the same functional framework as for implicit discretizations of the steady NSE (see \cite{CMAME15}, for instance), since $\widetilde{a}$ is an inner product on space $\Xv$ that generates a norm equivalent to the ${\bf H}^{1}$-norm.
Then, the existence of a solution follows from Brouwer's fixed point theorem \cite{Brezis11} as for the steady case.

\begin{itemize}
\item {\bf Velocity-pressure estimate.}
\end{itemize}

To prove estimate \eqref{eq:stabEst}, we set $\boldsymbol{\varphi}=2\Delta t\uv_{r}^{n+1}$, $\psi=2\Delta t p_{r}^{n+1}$ in \eqref{eq:discTempROM} and add both equations. Using the polarization identity:
$$
\left(\disp\frac{\uv_{r}^{n+1}-\uv_{r}^{n}}{\Delta t},2\Delta t\uv_{r}^{n+1}\right)=\nor{\uv_{r}^{n+1}}{{\bf L}^{2}}^{2}-\nor{\uv_{r}^{n}}{{\bf L}^{2}}^{2}+\nor{\uv_{r}^{n+1}-\uv_{r}^{n}}{{\bf L}^{2}}^{2},
$$
and noting that $b(\uv_{r}^{n+1},\uv_{r}^{n+1},\uv_{r}^{n+1})=0$ by \eqref{eq:formb}, we obtain:
\begin{eqnarray}\label{eq:stabEst1}
&&\nor{\uv_{r}^{n+1}}{{\bf L}^{2}}^{2}+\nor{\uv_{r}^{n+1}-\uv_{r}^{n}}{{\bf L}^{2}}^{2}+2\Delta t\nu\nor{\nabla\uv_{r}^{n+1}}{{\bf L}^{2}}^{2}+2\Delta t\nor{\Pi_{h}^{*}(\nabla p_{r}^{n+1})}{\tau}^{2}\nonumber \\ \nonumber \\
&+& 2\Delta t\sigma\nor{p_{r}^{n+1}}{L^{2}}^{2}=\nor{\uv_{r}^{n}}{{\bf L}^{2}}^{2}+2\Delta t\langle \fv^{n+1},\uv_{r}^{n+1}\rangle.
\end{eqnarray}
By definition of the dual norm and Young's inequality, from \eqref{eq:stabEst1} we get:
\begin{eqnarray}\label{eq:stabEst2}
&&\nor{\uv_{r}^{n+1}}{{\bf L}^{2}}^{2}+\nor{\uv_{r}^{n+1}-\uv_{r}^{n}}{{\bf L}^{2}}^{2}+\Delta t\nu\nor{\nabla\uv_{r}^{n+1}}{{\bf L}^{2}}^{2}+\Delta t\nor{\Pi_{h}^{*}(\nabla p_{r}^{n+1})}{\tau}^{2}+\Delta t\sigma\nor{p_{r}^{n+1}}{L^{2}}^{2}\nonumber \\
&\leq& \nor{\uv_{r}^{n}}{{\bf L}^{2}}^{2}+\frac{4\Delta t}{\nu}\sum_{n=0}^{N-1}\nor{\fv^{n+1}}{{\bf H}^{-1}}^{2}.
\end{eqnarray}
Then, the stability estimate \eqref{eq:stabEst} follows by summing \eqref{eq:stabEst2} from $n=0$ to $k\leq N-1$.
\qed

\begin{remark}\label{rm:POD-G-ROM}
The stability estimate \eqref{eq:stabEst}, which makes apparent the estimate of the
pressure stabilization term, guarantees an extra-control on the high frequencies of the pressure gradient that could lead to unstable discretization.
\end{remark}

\begin{remark}\label{rm:PresEst}
An alternative unconditional pressure stability estimate for LPS-ROM \eqref{eq:discTempROM} could be obtained if we do not add the penalty term with factor $\sigma$ to the variational formulation. In this case, when considering for instance Dirichlet boundary conditions on the whole boundary, the constant the pressure is determined up through the formulation could be fixed prescribing its value at a point of the boundary. When considering do nothing boundary conditions on a part of the boundary (outflow), the constant the pressure is determined up through the formulation is already fixed by the prescribed boundary conditions.

\medskip

To do so, one sets $\psi=0$ in \eqref{eq:discTempROM}. For any $\boldsymbol{\varphi}\in\Xv_r$, this yields:
$$
(p_{r}^{n+1},\div\boldsymbol{\varphi})=\left(\disp\frac{\uv_{r}^{n+1}-\uv_{r}^{n}}{\Delta t},\boldsymbol{\varphi}\right)+b(\uv_{r}^{n+1},\uv_{r}^{n+1},\boldsymbol{\varphi})+\nu(\nabla\uv_{r}^{n+1},\nabla\boldsymbol{\varphi})-\langle \fv^{n+1}, \boldsymbol{\varphi} \rangle.
$$
Let $P_{r}^{n+1}=\disp\sum_{k=0}^{n}\Delta t p_{r}^{k+1}$, then summation over the discrete times gives:
$$
(P_{r}^{n+1},\div\boldsymbol{\varphi})=\left(\uv_{r}^{n+1}-\uv_{r}^{0},\boldsymbol{\varphi}\right)+\sum_{k=0}^{n}\Delta t\left[b(\uv_{r}^{k+1},\uv_{r}^{k+1},\boldsymbol{\varphi})+\nu(\nabla\uv_{r}^{k+1},\nabla\boldsymbol{\varphi})-\langle \fv^{k+1}, \boldsymbol{\varphi} \rangle\right].
$$ 
Thus, we get:
\begin{eqnarray}\label{eq:ModInfSup}
&&\sup_{\boldsymbol{\varphi}\in\Xv_r}\frac{(P_{r}^{n+1},\div\boldsymbol{\varphi})}{\nor{\nabla\boldsymbol{\varphi}}{{\bf L}^2}}+\nor{\Pi_{h}^{*}(\nabla P_{r}^{n+1})}{\tau}\nonumber \\
&\leq& C\left[\nor{\uv_{r}^{n+1}}{{\bf L}^{2}}+\nor{\uv_{r}^{0}}{{\bf L}^{2}}+\sum_{k=0}^{n}\Delta t\left(\nor{\nabla\uv_{r}^{k+1}}{{\bf L}^{2}}^{2}+\nu\nor{\nabla\uv_{r}^{k+1}}{{\bf L}^{2}}+\nor{\fv^{k+1}}{{\bf H}^{-1}}\right)\right]\nonumber \\
&&+\nor{\Pi_{h}^{*}(\nabla P_{r}^{n+1})}{\tau}.
\end{eqnarray}
where we have applied triangle inequality, Cauchy--Schwarz inequality, the definition of the dual norm, and the standard estimate \eqref{eq:formb1} for the convective term. So, if we define the norm:
\BEQ\label{eq:NewNorm}
|||\cdot|||=\sup_{\boldsymbol{\varphi}\in\Xv_r}\frac{(\cdot,\div\boldsymbol{\varphi})}{\nor{\nabla\boldsymbol{\varphi}}{{\bf L}^2}}+\nor{\Pi_{h}^{*}(\nabla \cdot)}{\tau},
\EEQ
by Cauchy--Schwarz and triangle inequalities, from \eqref{eq:ModInfSup} we have:
\begin{eqnarray}\label{eq:ModInfSup1}
|||P_{r}^{n+1}|||&\leq&C\left[\max_{k=0,\ldots,N}\nor{\uv_{r}^{k}}{{\bf L}^{2}}+\sum_{k=0}^{N-1}\Delta t\nor{\nabla\uv_{r}^{k+1}}{{\bf L}^{2}}^{2}+\left(\sum_{k=0}^{N-1}\Delta t\nu\nor{\nabla\uv_{r}^{k+1}}{{\bf L}^{2}}^{2}\right)^{1/2}\right.\nonumber \\
&&+\left.\left(\sum_{k=0}^{N-1}\Delta t\nor{\fv^{k+1}}{{\bf H}^{-1}}^{2}\right)^{1/2}
+\left(\sum_{k=0}^{N-1}\Delta t\nor{\Pi_{h}^{*}(\nabla p_{r}^{k+1})}{\tau}^{2}\right)^{1/2}\right].
\end{eqnarray}
Using estimate \eqref{eq:stabEst} to bound the terms on the r.h.s. of \eqref{eq:ModInfSup1}, we obtain:
\BEQ\label{eq:ModInfSup2}
|||P_{r}^{n+1}|||\leq C\left[\nor{\uv_{r}^{0}}{{\bf L}^{2}}+\frac{1}{\sqrt{\nu}}\left(\sum_{k=0}^{N-1}\Delta t\nor{\fv^{k+1}}{{\bf H}^{-1}}^{2}\right)^{1/2}
+\frac{\nor{\uv_{r}^{0}}{{\bf L}^{2}}^{2}}{\nu}+\frac{1}{\nu^{2}}\sum_{k=0}^{N-1}\Delta t\nor{\fv^{k+1}}{{\bf H}^{-1}}^{2}\right],
\EEQ
which proves the unconditional stability of the time primitive of the pressure in the space induced by the norm \eqref{eq:NewNorm}, which we will call $\widetilde{Q}_{r}$ in the sequel. Indeed, due to the stability of $\uv_{0h}$ in $L^2$-norm and the regularity of $\fv$, we have:
\BEQ\label{eq:ModInfSup3}
\nor{P_{r}}{{L^{\infty}(\widetilde{Q}_{r})}}\leq \frac{C}{\nu^2},
\EEQ
and we have denoted $P_{r}(t)=\disp\int_{0}^{t}\widetilde{p}_{r}(s)\,ds$, being $\widetilde{p}_{r}$ the piecewise constant in time function that takes the value $p_{r}^{n+1}$ on $(t_{n},t_{n+1})$. Note that for the ROM pressure we can only obtain:
\BEQ\label{eq:stabPresEst}
\nor{\widetilde{p}_{r}}{{L^{1}(\widetilde{Q}_{r})}}\leq \frac{C}{\nu^2\sqrt{\Delta t}},
\EEQ
which is a not uniform bound (with respect to $\Delta t$) in $L^{1}(\widetilde{Q}_{r})$ space of space-time functions (cf. \cite{ChaconLewandowski14}, Remark 10.2), thus the technical trick of considering the time primitive of the pressure, as originally introduced in \cite{ChaconLewandowski14} for a FE-FOM, to prove its unconditional stability in this case. On the other side, when considering the penalty term with factor $\sigma> 0$ to fix the constant the pressure is determined up through the formulation (for Dirichlet boundary conditions on the whole boundary), as in Theorem \ref{th:stabEst}, thus we have:
\BEQ\label{eq:stabPresEst}
\nor{\widetilde{p}_{r}}{{L^{2}(L^{2})}}\leq \frac{C}{\sqrt{\nu\sigma}}.
\EEQ
\end{remark}

\subsection{Error estimates for LPS-ROM}\label{subsec:EE}

We are now in position to prove the following error estimate result for the LPS-ROM defined by  \eqref{eq:discTempROM}:
\begin{theorem}\label{th:PODEE}
Under the regularity assumption on the continuous solution (Hypothesis \ref{hp:ContRegularity}), the assumption on 
the stabilization parameters (Hypothesis \ref{hp:stabCoef}), and supposing that $\nor{\uv^{0}-\uv_{r}^{0}}{{\bf L}^{2}}=\mathcal{O}(h^{s})$, the solution of the LPS-ROM \eqref{eq:discTempROM}
satisfies the following error estimate: 
\begin{eqnarray}\label{eq:PODEEF}
&&\Delta t\sum_{n=0}^{N-1}\nor{\uv^{n+1}-\uv_{r}^{n+1}}{{\bf L}^{2}}^{2}+\nu\Delta t\sum_{n=0}^{N-1}\nor{\nabla(\uv^{n+1}-\uv_{r}^{n+1})}{{\bf L}^{2}}^{2} + \sigma\Delta t\sum_{n=0}^{N-1}\nor{p^{n+1}-p_{r}^{n+1}}{L^{2}}^{2}\nonumber \\
&\leq& C^{*}\left[\sigma+\frac{h^{2s}}{\nu^3}
+ \left(\frac{1}{\nu}+\frac{\alpha^2}{\sigma}\right)\left(h^{2s}+\Delta t^2\right)\left(1+\nor{S^v_r}{2}\right)+\frac{1}{\nu}\left(h^{2s}+\Delta t^2\right)\nor{S^p_r}{2}\right]\nonumber \\
&&+C^{*}\left[\left(\frac{1}{\nu}+\frac{\alpha^2}{\sigma}\right)\sum_{i=r+1}^{M_v}\left(1+\nor{\nabla\boldsymbol{\varphi}_{i}}{{\bf L}^2}^{2}\right)\lambda_i + \frac{1}{\nu}\sum_{i=r+1}^{M_p}\left(1+\nor{\nabla\psi_{i}}{{\bf L}^2}^{2}\right)\gamma_i\right],
\end{eqnarray}
where $C^*$ is a positive constant that will be determined throughout the proof.
\end{theorem}

{\bf Proof.} 
We start deriving the error bounds by splitting the error for the velocity and the pressure into two terms:
\BEQ\label{eq:splitErrVel}
\uv^{n+1}-\uv_{r}^{n+1}=(\uv^{n+1}-\vv_{r}^{n+1})-(\uv_{r}^{n+1}-\vv_{r}^{n+1})=\boldsymbol{\eta}^{n+1}-\boldsymbol{\phi}_{r}^{n+1},
\EEQ
\BEQ\label{eq:splitErrPres}
p^{n+1}-p_{r}^{n+1}=(p^{n+1}-z_{r}^{n+1})-(p_{r}^{n+1}-z_{r}^{n+1})=\rho^{n+1}-s_{r}^{n+1}.
\EEQ
In \eqref{eq:splitErrVel}, the first term, $\boldsymbol{\eta}^{n+1}=\uv^{n+1}-\vv_{r}^{n+1}$, represents the difference between $\uv^{n+1}$ and its $L^2$-orthogonal projection on $\Xv_r$ \eqref{eq:L2projVel}.
The second term, $\boldsymbol{\phi}_{r}^{n+1}=\uv_{r}^{n+1}-\vv_{r}^{n+1}$, is the remainder. Similarly, in \eqref{eq:splitErrPres}, the first term, $\rho^{n+1}=p^{n+1}-z_{r}^{n+1}$, represents the difference between $p^{n+1}$ and its $L^2$-orthogonal projection on $Q_r$ \eqref{eq:L2projPres}.
The second term, $s_{r}^{n+1}=p_{r}^{n+1}-z_{r}^{n+1}$, is the remainder. 

\medskip

Next, we construct the error equation. We first evaluate the weak formulation of the continuous NSE \eqref{eq:fvNS} at $t=t_{n+1}$ and let $\vv=\boldsymbol{\varphi}$, then subtract the 
LPS-ROM \eqref{eq:discTempROM} from it. For any $\boldsymbol{\varphi}\in\Xv_r$, we obtain:
\BEQ\label{eq:ErrEq}
\left \{ 
\begin{array}{lll}
&&(\partial_{t}\uv^{n+1},\boldsymbol{\varphi})-\disp\frac{1}{\Delta t}(\uv_{r}^{n+1}-\uv_{r}^{n},\boldsymbol{\varphi})+
\nu\left(\nabla(\uv^{n+1}-\uv_{r}^{n+1}),\nabla\boldsymbol{\varphi}\right)\\ \\
&+&b(\uv^{n+1},\uv^{n+1},\boldsymbol{\varphi})-b(\uv_{r}^{n+1},\uv_{r}^{n+1},\boldsymbol{\varphi})-(p^{n+1}-p_{r}^{n+1},\div\boldsymbol{\varphi})\\ \\
&+&\left(\div(\uv^{n+1}-\uv_{r}^{n+1}),\psi\right)-\left(\Pi_{h}^{*}(\nabla p_{r}^{n+1}),\Pi_{h}^{*}(\nabla\psi)\right)_{\tau}+\sigma(p^{n+1}-p_{r}^{n+1},\psi)=0.
\end{array}
\right.
\EEQ
By adding and subtracting the difference quotient term $\disp\frac{1}{\Delta t}(\uv^{n+1}-\uv^{n},\boldsymbol{\varphi})$ in \eqref{eq:ErrEq}, and applying the decompositions \eqref{eq:splitErrVel}-\eqref{eq:splitErrPres}, 
we get, for any $(\boldsymbol{\varphi},\psi)\in\Xv_{r}\times Q_r$:
\BEQ\label{eq:ErrEq0}
\left \{ 
\begin{array}{lll}
&&(\partial_{t}\uv^{n+1}-\disp\frac{\uv^{n+1}-\uv^{n}}{\Delta t},\boldsymbol{\varphi})+\disp\frac{1}{\Delta t}(\boldsymbol{\eta}^{n+1}-\boldsymbol{\phi}_{r}^{n+1},\boldsymbol{\varphi})
-\disp\frac{1}{\Delta t}(\boldsymbol{\eta}^{n}-\boldsymbol{\phi}_{r}^{n},\boldsymbol{\varphi})\\ \\
&+&\nu\left(\nabla(\boldsymbol{\eta}^{n+1}-\boldsymbol{\phi}_{r}^{n+1}),\nabla\boldsymbol{\varphi}\right)+b(\uv^{n+1},\uv^{n+1},\boldsymbol{\varphi})-b(\uv_{r}^{n+1},\uv_{r}^{n+1},\boldsymbol{\varphi})\\ \\
&-&(\rho^{n+1}-s_{r}^{n+1},\div\boldsymbol{\varphi})+\left(\div(\boldsymbol{\eta}^{n+1}-\boldsymbol{\phi}_{r}^{n+1}),\psi\right)-\left(\Pi_{h}^{*}(\nabla p_{r}^{n+1}),\Pi_{h}^{*}(\nabla\psi)\right)_{\tau}\\ \\
&+&\sigma(\rho^{n+1}-s_{r}^{n+1},\psi)=0.
\end{array}
\right.
\EEQ
Note that $\disp\frac{1}{\Delta t}(\boldsymbol{\eta}^{n+1}-\boldsymbol{\eta}^{n},\boldsymbol{\varphi})=0$, since $\vv_{r}^{n+1}$ is the $L^2$-orthogonal projection of $\uv^{n+1}$ on $\Xv_r$, and $\sigma(\rho^{n+1},\psi)=0$, since $z_{r}^{n+1}$ is the $L^2$-orthogonal projection of $p^{n+1}$ on $Q_r$. 
Choosing $\boldsymbol{\varphi}=2\Delta t\boldsymbol{\phi}_{r}^{n+1}$, $\psi=2\Delta t s_{r}^{n+1}$ in \eqref{eq:ErrEq0} and letting $\boldsymbol{c}^{n}=\partial_{t}\uv^{n+1}-\disp\frac{\uv^{n+1}-\uv^{n}}{\Delta t}$, we get:
\BEQ\label{eq:ErrEq1}
\left \{ 
\begin{array}{lll}
&&\left(\disp\frac{\boldsymbol{\phi}_{r}^{n+1}-\boldsymbol{\phi}_{r}^{n}}{\Delta t},2\Delta t\boldsymbol{\phi}_{r}^{n+1}\right)
+2\Delta t\nu\nor{\nabla\boldsymbol{\phi}_{r}^{n+1}}{{\bf L}^{2}}^{2}+2\Delta t\sigma\nor{s_{r}^{n+1}}{L^2}^{2}\\ \\
&=2\Delta t&\left[\langle\boldsymbol{c}^{n},\boldsymbol{\phi}_{r}^{n+1}\rangle+b(\uv^{n+1},\uv^{n+1},\boldsymbol{\phi}_{r}^{n+1})-b(\uv_{r}^{n+1},\uv_{r}^{n+1},\boldsymbol{\phi}_{r}^{n+1})\right.\\ \\
&&+\nu(\nabla\boldsymbol{\eta}^{n+1},\nabla\boldsymbol{\phi}_{r}^{n+1})-(\rho^{n+1},\div\boldsymbol{\phi}_{r}^{n+1})+(\div\boldsymbol{\eta}^{n+1},s_{r}^{n+1})\\ \\
&&\left.-\left(\Pi_{h}^{*}(\nabla p_{r}^{n+1}),\Pi_{h}^{*}(\nabla s_{r}^{n+1})\right)_{\tau}\right].
\end{array}
\right.
\EEQ
Using the polarization identity:
$$
\left(\disp\frac{\boldsymbol{\phi}_{r}^{n+1}-\boldsymbol{\phi}_{r}^{n}}{\Delta t},2\Delta t\boldsymbol{\phi}_{r}^{n+1}\right)=
\nor{\boldsymbol{\phi}_{r}^{n+1}}{{\bf L}^{2}}^{2}-\nor{\boldsymbol{\phi}_{r}^{n}}{{\bf L}^{2}}^{2}+\nor{\boldsymbol{\phi}_{r}^{n+1}-\boldsymbol{\phi}_{r}^{n}}{{\bf L}^{2}}^{2},
$$
from \eqref{eq:ErrEq1} we obtain:
\BEQ\label{eq:ErrEq2}
\left \{ 
\begin{array}{lll}
&&\nor{\boldsymbol{\phi}_{r}^{n+1}}{{\bf L}^{2}}^{2}-\nor{\boldsymbol{\phi}_{r}^{n}}{{\bf L}^{2}}^{2}+\nor{\boldsymbol{\phi}_{r}^{n+1}-\boldsymbol{\phi}_{r}^{n}}{{\bf L}^{2}}^{2}
+2\Delta t\nu\nor{\nabla\boldsymbol{\phi}_{r}^{n+1}}{{\bf L}^{2}}^{2}+2\Delta t\sigma\nor{s_{r}^{n+1}}{L^2}^{2}\\ \\
&=2\Delta t&\left[\langle\boldsymbol{c}^{n},\boldsymbol{\phi}_{r}^{n+1}\rangle+b(\uv^{n+1},\uv^{n+1},\boldsymbol{\phi}_{r}^{n+1})-b(\uv_{r}^{n+1},\uv_{r}^{n+1},\boldsymbol{\phi}_{r}^{n+1})\right.\\ \\
&&+\nu(\nabla\boldsymbol{\eta}^{n+1},\nabla\boldsymbol{\phi}_{r}^{n+1})-(\rho^{n+1},\div\boldsymbol{\phi}_{r}^{n+1})+(\div\boldsymbol{\eta}^{n+1},s_{r}^{n+1})\\ \\
&&\left.-\left(\Pi_{h}^{*}(\nabla p_{r}^{n+1}),\Pi_{h}^{*}(\nabla s_{r}^{n+1})\right)_{\tau}\right].
\end{array}
\right.
\EEQ
We estimate the terms on the r.h.s. of \eqref{eq:ErrEq2} one by one. By definition of the dual norm and Young's inequality, we get for the first term on the r.h.s. of \eqref{eq:ErrEq2}:
\BEQ\label{eq:timederErrEq}
\langle\boldsymbol{c}^{n},\boldsymbol{\phi}_{r}^{n+1}\rangle\leq \nor{\boldsymbol{c}^{n}}{{\bf H}^{-1}}\nor{\nabla\boldsymbol{\phi}_{r}^{n+1}}{{\bf L}^2}\leq \frac{\varepsilon_{1}^{-1}}{4}\nor{\boldsymbol{c}^{n}}{{\bf H}^{-1}}^{2}+\varepsilon_{1}\nor{\nabla\boldsymbol{\phi}_{r}^{n+1}}{{\bf L}^{2}}^{2},
\EEQ
for some small positive constant $\varepsilon_{1}$ (to be determined later). 

\medskip

The nonlinear convective terms in \eqref{eq:ErrEq2} can be written as follows:
\begin{eqnarray}\label{eq:convErrEq}
&&b(\uv^{n+1},\uv^{n+1},\boldsymbol{\phi}_{r}^{n+1})-b(\uv_{r}^{n+1},\uv_{r}^{n+1},\boldsymbol{\phi}_{r}^{n+1})\nonumber \\
&=&b(\uv_{r}^{n+1},\boldsymbol{\eta}^{n+1}-\boldsymbol{\phi}_{r}^{n+1},\boldsymbol{\phi}_{r}^{n+1})
+b(\boldsymbol{\eta}^{n+1}-\boldsymbol{\phi}_{r}^{n+1},\uv^{n+1},\boldsymbol{\phi}_{r}^{n+1})\nonumber \\
&=&b(\uv_{r}^{n+1},\boldsymbol{\eta}^{n+1},\boldsymbol{\phi}_{r}^{n+1})
+b(\boldsymbol{\eta}^{n+1},\uv^{n+1},\boldsymbol{\phi}_{r}^{n+1})
-b(\boldsymbol{\phi}_{r}^{n+1},\uv^{n+1},\boldsymbol{\phi}_{r}^{n+1}),
\end{eqnarray}
where we have used $b(\uv_{r}^{n+1},\boldsymbol{\phi}_{r}^{n+1},\boldsymbol{\phi}_{r}^{n+1})=0$, which follows from \eqref{eq:formb}. Next, we estimate each term on the r.h.s. of \eqref{eq:convErrEq}. Since $\uv_{r}^{n+1},\boldsymbol{\eta}^{n+1},\boldsymbol{\phi}_{r}^{n+1},\uv^{n+1}\in\Xv$, 
we can apply the standard bound \eqref{eq:formb1} for the trilinear form $b(\cdot,\cdot,\cdot)$ and use Young's inequality to get:
\begin{eqnarray}\label{eq:conv1}
b(\uv_{r}^{n+1},\boldsymbol{\eta}^{n+1},\boldsymbol{\phi}_{r}^{n+1})&\leq&
C\nor{\nabla\uv_{r}^{n+1}}{{\bf L}^{2}}\nor{\nabla\boldsymbol{\eta}^{n+1}}{{\bf L}^{2}}\nor{\nabla\boldsymbol{\phi}_{r}^{n+1}}{{\bf L}^{2}}\nonumber \\ 
&\leq& \frac{\varepsilon_{1}^{-1}C^{2}}{4}\nor{\nabla\uv_{r}^{n+1}}{{\bf L}^{2}}^{2}\nor{\nabla\boldsymbol{\eta}^{n+1}}{{\bf L}^{2}}^{2} +\varepsilon_{1}\nor{\nabla\boldsymbol{\phi}_{r}^{n+1}}{{\bf L}^{2}}^{2};
\end{eqnarray}
\begin{eqnarray}\label{eq:conv2}
b(\boldsymbol{\eta}^{n+1},\uv^{n+1},\boldsymbol{\phi}_{r}^{n+1})&\leq&
C\nor{\nabla\boldsymbol{\eta}^{n+1}}{{\bf L}^{2}}\nor{\nabla\uv^{n+1}}{{\bf L}^{2}}\nor{\nabla\boldsymbol{\phi}_{r}^{n+1}}{{\bf L}^{2}}\nonumber \\ 
&\leq& \frac{\varepsilon_{1}^{-1}C^{2}}{4}\nor{\nabla\uv^{n+1}}{{\bf L}^{2}}^{2}\nor{\nabla\boldsymbol{\eta}^{n+1}}{{\bf L}^{2}}^{2} +\varepsilon_{1}\nor{\nabla\boldsymbol{\phi}_{r}^{n+1}}{{\bf L}^{2}}^{2}.
\end{eqnarray}
For the last nonlinear convective term, applying H\"older's inequality, Sobolev embedding theorem, and Young's inequality yields:
\begin{eqnarray}\label{eq:conv3}
b(\boldsymbol{\phi}_{r}^{n+1},\uv^{n+1},\boldsymbol{\phi}_{r}^{n+1})&\leq&
C\nor{\boldsymbol{\phi}_{r}^{n+1}}{{\bf L}^{2}}(\nor{\nabla\uv^{n+1}}{{\bf L}^{3}}+\nor{\uv^{n+1}}{{\bf L}^{\infty}})\nor{\nabla\boldsymbol{\phi}_{r}^{n+1}}{{\bf L}^{2}}\nonumber \\  
&\leq& \frac{\varepsilon_{1}^{-1}C^{2}}{4}\nor{\nabla\uv^{n+1}}{{\bf H}^{1}}^{2}\nor{\boldsymbol{\phi}_{r}^{n+1}}{{\bf L}^{2}}^{2} 
+\varepsilon_{1}\nor{\nabla\boldsymbol{\phi}_{r}^{n+1}}{{\bf L}^{2}}^{2}.
\end{eqnarray}
By Cauchy--Schwarz and Young's inequalities, we bound the fourth and fifth terms on the r.h.s. of \eqref{eq:ErrEq2}:
\begin{eqnarray}
\nu(\nabla\boldsymbol{\eta}^{n+1},\boldsymbol{\phi}_{r}^{n+1})&\leq& \frac{\varepsilon_{1}^{-1}\nu^2}{4}\nor{\nabla\boldsymbol{\eta}^{n+1}}{{\bf L}^{2}}^{2}+\varepsilon_{1}\nor{\nabla\boldsymbol{\phi}_{r}^{n+1}}{{\bf L}^{2}}^{2},\label{eq:diffErrEq}\\
-(\rho^{n+1},\div\boldsymbol{\phi}_{r}^{n+1})&\leq& \frac{\varepsilon_{1}^{-1}}{4}\nor{\rho^{n+1}}{{\bf L}^{2}}^{2}+\varepsilon_{1}\nor{\nabla\boldsymbol{\phi}_{r}^{n+1}}{{\bf L}^{2}}^{2}.\label{eq:presErrEq}
\end{eqnarray}
For the sixth term on the r.h.s of \eqref{eq:ErrEq2}, applying triangle, Cauchy--Schwarz and Young's inequalities plus Lemma \ref{lm:AngleCond}, we have:
\begin{eqnarray}\label{eq:divErrEq}
(\div\boldsymbol{\eta}^{n+1},s_{r}^{n+1})&\leq& |(\div\uv^{n+1},s_{r}^{n+1})|+|(\div\vv_{r}^{n+1},s_{r}^{n+1})|\nonumber \\
&\leq& \sigma\nor{p^{n+1}}{L^2}\nor{s_{r}^{n+1}}{L^2}+\alpha\nor{\div\vv_{r}^{n+1}}{L^2}\nor{s_{r}^{n+1}}{L^2}\nonumber \\
&\leq& \sigma\left(\nor{\rho^{n+1}}{L^2}+\nor{z_{r}^{n+1}}{L^2}\right)\nor{s_{r}^{n+1}}{L^2}\nonumber \\
&&+\alpha\left(\nor{\div\boldsymbol{\eta}^{n+1}}{L^2}+\nor{\div\uv^{n+1}}{L^2}\right)\nor{s_{r}^{n+1}}{L^2}\nonumber \\
&\leq& \frac{\varepsilon_{2}^{-1}\sigma^{2}}{4}\left(\nor{\rho^{n+1}}{L^2}^{2}+\nor{z_{r}^{n+1}}{L^2}^{2}\right)
+\frac{\varepsilon_{2}^{-1}\alpha^{2}}{4}\left(\nor{\nabla\boldsymbol{\eta}^{n+1}}{{\bf L}^2}^{2}+\nor{\div\uv^{n+1}}{L^2}^{2}\right)\nonumber \\
&&+4\varepsilon_{2}\nor{s_{r}^{n+1}}{L^2}^{2},
\end{eqnarray}
for some small positive constant $\varepsilon_{2}$ (to be determined later). In particular, we have applied Lemma \ref{lm:AngleCond} with $Y=\text{span}\{\div\boldsymbol{\varphi}_{1},\ldots,\div\boldsymbol{\varphi}_{r}\}$ and $Z=Q_r$.

\medskip

Finally, for the last term on the r.h.s. of \eqref{eq:ErrEq2} we have:
\begin{eqnarray}\label{eq:stabErrEq}
-\left(\Pi_{h}^{*}(\nabla p_{r}^{n+1}),\Pi_{h}^{*}(\nabla s_{r}^{n+1})\right)_{\tau}&=&\left(\Pi_{h}^{*}(\nabla \rho_{r}^{n+1}),\Pi_{h}^{*}(\nabla s_{r}^{n+1})\right)_{\tau}-\left(\Pi_{h}^{*}(\nabla s_{r}^{n+1}),\Pi_{h}^{*}(\nabla s_{r}^{n+1})\right)_{\tau}\nonumber \\
&&-\left(\Pi_{h}^{*}(\nabla p^{n+1}),\Pi_{h}^{*}(\nabla s_{r}^{n+1})\right)_{\tau}\nonumber \\
&\leq& \nor{\Pi_{h}^{*}(\nabla \rho_{r}^{n+1})}{\tau}\nor{\Pi_{h}^{*}(\nabla s_{r}^{n+1})}{\tau}-\nor{\Pi_{h}^{*}(\nabla s_{r}^{n+1})}{\tau}^{2}\nonumber \\
&&+\nor{\Pi_{h}^{*}(\nabla p^{n+1})}{\tau}\nor{\Pi_{h}^{*}(\nabla s_{r}^{n+1})}{\tau}\nonumber \\
&\leq& \nor{\Pi_{h}^{*}(\nabla \rho_{r}^{n+1})}{\tau}^{2}-\frac{1}{2}\nor{\Pi_{h}^{*}(\nabla s_{r}^{n+1})}{\tau}^{2}\nonumber \\
&&+\nor{\Pi_{h}^{*}(\nabla p^{n+1})}{\tau}^{2},
\end{eqnarray}
where we have used Cauchy--Schwarz and Young's inequalities. 

\medskip

Substituting inequalities \eqref{eq:timederErrEq}-\eqref{eq:stabErrEq} in \eqref{eq:ErrEq2} and taking $\varepsilon_{1}=\nu/12$, $\varepsilon_{2}=\sigma/8$, we obtain:
\begin{eqnarray}\label{eq:ErrEq3}
&&\nor{\boldsymbol{\phi}_{r}^{n+1}}{{\bf L}^{2}}^{2}-\nor{\boldsymbol{\phi}_{r}^{n}}{{\bf L}^{2}}^{2}
+\nor{\boldsymbol{\phi}_{r}^{n+1}-\boldsymbol{\phi}_{r}^{n}}{{\bf L}^2}^{2}
+\Delta t\nu\nor{\nabla\boldsymbol{\phi}_{r}^{n+1}}{{\bf L}^{2}}^{2}
+\Delta t\sigma\nor{s_{r}^{n+1}}{L^{2}}^{2}\nonumber \\
&&+\Delta t\nor{\Pi_{h}^{*}(\nabla s_{r}^{n+1})}{\tau}^{2}\nonumber \\
&\leq 2\Delta t&\left[\frac{3}{\nu}\nor{\boldsymbol{c}^{n}}{{\bf H}^{-1}}^{2}
+\frac{3 C^{2}}{\nu}\nor{\nabla\uv_{r}^{n+1}}{{\bf L}^{2}}^{2}\nor{\nabla\boldsymbol{\eta}^{n+1}}{{\bf L}^{2}}^{2}
+\frac{3 C^{2}}{\nu}\nor{\nabla\uv^{n+1}}{{\bf L}^{2}}^{2}\nor{\nabla\boldsymbol{\eta}^{n+1}}{{\bf L}^{2}}^{2}\right.
\nonumber \\
&&+\frac{3 C^{2}}{\nu}\nor{\nabla\uv^{n+1}}{{\bf H}^{1}}^{2}\nor{\boldsymbol{\phi}_{r}^{n+1}}{{\bf L}^{2}}^{2}
+3\nu\nor{\nabla\boldsymbol{\eta}^{n+1}}{{\bf L}^{2}}^{2}
+\frac{3}{\nu}\nor{\rho^{n+1}}{L^{2}}^{2}\nonumber \\
&&+2\sigma\left(\nor{\rho^{n+1}}{L^2}^{2}+\nor{z_{r}^{n+1}}{L^2}^{2}\right)
+\frac{2\alpha^2}{\sigma}\left(\nor{\nabla\boldsymbol{\eta}^{n+1}}{{\bf L}^2}^{2}+\nor{\div\uv^{n+1}}{L^2}^{2}\right)\nonumber \\
&&\left.+\nor{\Pi_{h}^{*}(\nabla \rho_{r}^{n+1})}{\tau}^{2}+\nor{\Pi_{h}^{*}(\nabla p^{n+1})}{\tau}^{2}\right]\nonumber \\
&\leq \disp\frac{C\Delta t}{\nu}&\left[\nor{\boldsymbol{c}^{n}}{{\bf H}^{-1}}^{2}
+\left(\nor{\nabla\uv_{r}^{n+1}}{{\bf L}^{2}}^{2}+\nor{\nabla\uv^{n+1}}{{\bf L}^{2}}^{2}+\nu^2+\frac{\alpha^2\nu}{\sigma}\right)\nor{\nabla\boldsymbol{\eta}^{n+1}}{{\bf L}^{2}}^{2}\right.\nonumber \\
&&\left.+\nor{\nabla\uv^{n+1}}{{\bf H}^{1}}^{2}\nor{\boldsymbol{\phi}_{r}^{n+1}}{{\bf L}^{2}}^{2}+\left(1+\sigma\nu\right)\nor{\rho^{n+1}}{L^2}^{2}
+\sigma\nu\nor{z_{r}^{n+1}}{L^2}^{2}+\frac{\alpha^2\nu}{\sigma}\nor{\div\uv^{n+1}}{L^2}^{2}\right]
\nonumber \\
&+2\Delta t&\left(\nor{\Pi_{h}^{*}(\nabla \rho_{r}^{n+1})}{\tau}^{2}+\nor{\Pi_{h}^{*}(\nabla p^{n+1})}{\tau}^{2}\right)\nonumber \\
&\leq \disp\frac{C\Delta t}{\nu}&\left[\nor{\boldsymbol{c}^{n}}{{\bf H}^{-1}}^{2}
+\left(\nor{\nabla\uv_{r}^{n+1}}{{\bf L}^{2}}^{2}+\nor{\nabla\uv^{n+1}}{{\bf L}^{2}}^{2}+\nu^2+\frac{\alpha^2\nu}{\sigma}\right)\nor{\nabla\boldsymbol{\eta}^{n+1}}{{\bf L}^{2}}^{2}\right.\nonumber \\
&&\left.+\nor{\nabla\uv^{n+1}}{{\bf H}^{1}}^{2}\nor{\boldsymbol{\phi}_{r}^{n+1}}{{\bf L}^{2}}^{2}+\left(1+\sigma\nu\right)\nor{\rho^{n+1}}{L^2}^{2}
+(1+\alpha^2)\sigma\nu\nor{p^{n+1}}{L^2}^{2}\right]
\nonumber \\
&+2\Delta t&\left(\nor{\Pi_{h}^{*}(\nabla \rho_{r}^{n+1})}{\tau}^{2}+\nor{\Pi_{h}^{*}(\nabla p^{n+1})}{\tau}^{2}\right),
\end{eqnarray}
where in the last inequality we have used the identity $z_{r}^{n+1}=p^{n+1}-\rho^{n+1}$ from \eqref{eq:splitErrPres}, triangle inequality, and the inequality $\nor{\div\uv^{n+1}}{L^2}^{2}\leq \sigma^2\nor{p^{n+1}}{L^2}^{2}$ from \eqref{eq:fvNS}.

\medskip

Summing \eqref{eq:ErrEq3} from $n=0$ to $k\leq N-1$, we have:
\begin{eqnarray}\label{eq:ErrEqSum}
&&\max_{0\leq k\leq N}\nor{\boldsymbol{\phi}_{r}^{k}}{{\bf L}^{2}}^{2}
+\sum_{n=0}^{N-1}\nor{\boldsymbol{\phi}_{r}^{n+1}-\boldsymbol{\phi}_{r}^{n}}{{\bf L}^2}^{2}
+\nu\Delta t\sum_{n=0}^{N-1}\nor{\nabla\boldsymbol{\phi}_{r}^{n+1}}{{\bf L}^{2}}^{2}
+\sigma\Delta t\sum_{n=0}^{N-1}\nor{s_{r}^{n+1}}{L^{2}}^{2}\nonumber \\
&&+\Delta t\sum_{n=0}^{N-1}\nor{\Pi_{h}^{*}(\nabla s_{r}^{n+1})}{\tau}^{2}\nonumber \\
&\leq&\nor{\boldsymbol{\phi}_{r}^{0}}{{\bf L}^{2}}^{2}+\frac{C\Delta t}{\nu}\sum_{n=0}^{N-1}\left[\nor{\boldsymbol{c}^{n}}{{\bf H}^{-1}}^{2}
+\left(\nor{\nabla\uv_{r}^{n+1}}{{\bf L}^{2}}^{2}+\nor{\nabla\uv^{n+1}}{{\bf L}^{2}}^{2}+\nu^2+\frac{\alpha^2\nu}{\sigma}\right)\nor{\nabla\boldsymbol{\eta}^{n+1}}{{\bf L}^{2}}^{2}\right]\nonumber \\
&&+\frac{C\Delta t}{\nu}\sum_{n=0}^{N-1}\left[\nor{\nabla\uv^{n+1}}{{\bf H}^{1}}^{2}\nor{\boldsymbol{\phi}_{r}^{n+1}}{{\bf L}^{2}}^{2}+\left(1+\sigma\nu\right)\nor{\rho^{n+1}}{L^2}^{2}
+(1+\alpha^2)\sigma\nu\nor{p^{n+1}}{L^2}^{2}\right]
\nonumber \\
&&+2\Delta t\sum_{n=0}^{N-1}\left(\nor{\Pi_{h}^{*}(\nabla \rho_{r}^{n+1})}{\tau}^{2}+\nor{\Pi_{h}^{*}(\nabla p^{n+1})}{\tau}^{2}\right).
\end{eqnarray}
Next, we estimate each term on the r.h.s. of \eqref{eq:ErrEqSum}.

\medskip

The first term on the r.h.s. of \eqref{eq:ErrEqSum} can be estimated as follows:
\BEQ\label{eq:ErrIneq2RHS1}
\nor{\boldsymbol{\phi}_{r}^{0}}{{\bf L}^2}^{2}\leq \nor{\uv^{0}-\uv_{r}^{0}}{{\bf L}^2}^{2}+\nor{\uv^{0}-\vv_{r}^{0}}{{\bf L}^2}^{2}\leq Ch^{2s},
\EEQ
where the last inequality follows from the fact that $\vv_{r}^{0}$ is the $L^2$-orthogonal projection of $\uv^{0}$ on $\Xv_{r}\subset \Xv_{h}$, so that it satisfies optimal approximation properties similar to standard FE interpolations (cf. \cite{Ciarlet02}), and we have supposed $\nor{\uv^{0}-\uv_{r}^{0}}{}=\mathcal{O}(h^{s})$. 

\medskip

By Taylor's theorem, the second term on the r.h.s. of \eqref{eq:ErrEqSum} can be estimated as follows:
\BEQ\label{eq:cn}
\Delta t\sum_{n=0}^{N-1}\nor{\boldsymbol{c}^{n}}{{\bf H}^{-1}}^{2}\leq C\Delta t\sum_{n=0}^{N-1}\nor{\boldsymbol{c}^{n}}{{\bf L}^{2}}^{2}\leq C\Delta t^{2}\nor{\partial_{t}^{2}\uv}{L^{2}({\bf L}^{2})}^{2}\leq C\Delta t^{2},
\EEQ
where we have used the regularity assumption \ref{hp:ContRegularity} on the continuous velocity.

\medskip

To estimate the third term on the r.h.s. of \eqref{eq:ErrEqSum}, we use Theorem \ref{th:stabEst} and the fact that $\vv_{r}^{n+1}$ is the $L^2$-orthogonal projection of $\uv^{n+1}$ on $\Xv_{r}\subset\Xv_{h}$, so that it satisfies optimal approximation properties as standard FE interpolations (cf. \cite{Ciarlet02}): 
\BEQ\label{eq:eta2}
\Delta t\sum_{n=0}^{N-1}\nor{\nabla\uv_{r}^{n+1}}{{\bf L}^{2}}^{2}\nor{\nabla\boldsymbol{\eta}^{n+1}}{{\bf L}^{2}}^{2}
\leq \frac{C}{\nu^2}h^{2s}.
\EEQ

By using the regularity assumption \ref{hp:ContRegularity} on the continuous velocity and $H^{1}$ velocity projection error estimate \eqref{eq:PODlmH01vel}, the fourth term on the r.h.s. of \eqref{eq:ErrEqSum} can be estimated as follows:
\BEQ\label{eq:eta3}
\Delta t\sum_{n=0}^{N-1}\nor{\nabla\uv^{n+1}}{{\bf L}^{2}}^{2}\nor{\nabla\boldsymbol{\eta}^{n+1}}{{\bf L}^{2}}^{2}
\leq C\left[\left(1+\nor{S^v_r}{2}\right)\left(h^{2s}+\Delta t ^{2}\right) + \sum_{i=r+1}^{M_v}\nor{\nabla\boldsymbol{\varphi}_{i}}{{\bf L}^2}^{2}\lambda_{i}\right].
\EEQ

\medskip

Similarly, we have: 
\begin{eqnarray}\label{eq:eta4}
&&\Delta t\sum_{n=0}^{N-1}\left(\nu^2 + \frac{\alpha^2\nu}{\sigma}\right)\nor{\nabla\boldsymbol{\eta}^{n+1}}{{\bf L}^{2}}^{2}\nonumber \\
&\leq& C\frac{\alpha^2\nu}{\sigma}\left[\left(1+\nor{S^v_r}{2}\right)\left(h^{2s}+\Delta t ^{2}\right) + \sum_{i=r+1}^{M_v}\nor{\nabla\boldsymbol{\varphi}_{i}}{{\bf L}^2}^{2}\lambda_{i}\right],
\end{eqnarray}
and, by using $L^2$ pressure projection error estimate \eqref{eq:PODlmL2pres}:
\BEQ\label{eq:rho}
\Delta t \sum_{n=0}^{N-1}(1+\sigma\nu)\nor{\rho^{n+1}}{L^{2}}^{2}\leq C\left(h^{2s}+\Delta t^2+ \sum_{i=r+1}^{M_p}\gamma_{i}\right).
\EEQ

\medskip

The ninth term on the r.h.s. of \eqref{eq:ErrEqSum} makes apparent the convergence order reduction linked to the diffusive nature of the penalty term. Indeed, by using the regularity assumption \ref{hp:ContRegularity} on the continuous pressure, we get: 
\BEQ\label{eq:pres}
\Delta t \sum_{n=0}^{N-1}(1+\alpha^2)\sigma\nu\nor{p^{n+1}}{L^{2}}^{2}\leq C\sigma\nu.
\EEQ

\medskip

For the last two terms on the r.h.s. of \eqref{eq:ErrEqSum}, we obtain:
\begin{eqnarray}\label{eq:presStab1}
\Delta t \sum_{n=0}^{N-1}\nor{\Pi_{h}^{*}(\nabla \rho_{r}^{n+1})}{\tau}^{2}&\leq& C h^2 \Delta t \sum_{n=0}^{N-1}\nor{\nabla \rho_{r}^{n+1}}{{\bf L}^2}^{2}\nonumber \\
&\leq& C\left[\left(1+\nor{S^p_r}{2}\right)\left(h^{2s}+\Delta t ^{2}\right) + \sum_{i=r+1}^{M_p}\nor{\nabla\psi_{i}}{{\bf L}^2}^{2}\gamma_{i}\right],
\end{eqnarray}
and:
\BEQ\label{eq:presStab2}
\Delta t \sum_{n=0}^{N-1}\nor{\Pi_{h}^{*}(\nabla p^{n+1})}{\tau}^{2}\leq C h^2 \Delta t \sum_{n=0}^{N-1}\nor{\Pi_{h}^{*}(\nabla p^{n+1})}{{\bf L}^2}^{2}\leq C h^{2s},
\EEQ
where we have used the assumption on stabilization parameter \ref{hp:stabCoef}, the regularity assumption \ref{hp:ContRegularity} on the continuous pressure, stability and optimal error estimates of $\Pi_{h}$ (cf. \cite{Ciarlet02}),  and $H^{1}$ pressure projection error estimate \eqref{eq:PODlmH1pres}.

\medskip

Collecting \eqref{eq:ErrIneq2RHS1}-\eqref{eq:presStab2}, estimate \eqref{eq:ErrEqSum} becomes:
\begin{eqnarray}\label{eq:ErrEqRiep}
&&\max_{0\leq k\leq N}\nor{\boldsymbol{\phi}_{r}^{k}}{{\bf L}^{2}}^{2}
+\sum_{n=0}^{N-1}\nor{\boldsymbol{\phi}_{r}^{n+1}-\boldsymbol{\phi}_{r}^{n}}{{\bf L}^2}^{2}
+\nu\Delta t\sum_{n=0}^{N-1}\nor{\nabla\boldsymbol{\phi}_{r}^{n+1}}{{\bf L}^{2}}^{2}
+\sigma\Delta t\sum_{n=0}^{N-1}\nor{s_{r}^{n+1}}{L^{2}}^{2}\nonumber \\
&&+\Delta t\sum_{n=0}^{N-1}\nor{\Pi_{h}^{*}(\nabla s_{r}^{n+1})}{\tau}^{2}\nonumber \\
&\leq&\frac{C\Delta t}{\nu}\sum_{n=0}^{N-1}\nor{\nabla\uv^{n+1}}{{\bf H}^{1}}^{2}\nor{\boldsymbol{\phi}_{r}^{n+1}}{{\bf L}^{2}}^{2}
+C\left(\sigma+\frac{h^{2s}}{\nu^3}\right)\nonumber \\
&&+C\left[ \left(\frac{1}{\nu}+\frac{\alpha^2}{\sigma}\right)\left(h^{2s}+\Delta t^2\right)\left(1+\nor{S^v_r}{2}\right)+\frac{1}{\nu}\left(h^{2s}+\Delta t^2\right)\nor{S^p_r}{2}\right]\nonumber \\
&&+C\left[\left(\frac{1}{\nu}+\frac{\alpha^2}{\sigma}\right)\sum_{i=r+1}^{M_v}\nor{\nabla\boldsymbol{\varphi}_{i}}{{\bf L}^2}^{2}\lambda_i + \frac{1}{\nu}\sum_{i=r+1}^{M_p}\left(1+\nor{\nabla\psi_{i}}{{\bf L}^2}^{2}\right)\gamma_i\right].
\end{eqnarray}
The discrete Gr\"onwall's lemma (cf. Lemma 27 in \cite{Layton08}) implies the following inequality:
\begin{eqnarray}\label{eq:ErrEqRiepF}
&&\max_{0\leq k\leq N}\nor{\boldsymbol{\phi}_{r}^{k}}{{\bf L}^{2}}^{2}
+\sum_{n=0}^{N-1}\nor{\boldsymbol{\phi}_{r}^{n+1}-\boldsymbol{\phi}_{r}^{n}}{{\bf L}^2}^{2}
+\nu\Delta t\sum_{n=0}^{N-1}\nor{\nabla\boldsymbol{\phi}_{r}^{n+1}}{{\bf L}^{2}}^{2}
+\sigma\Delta t\sum_{n=0}^{N-1}\nor{s_{r}^{n+1}}{L^{2}}^{2}\nonumber \\
&&+\Delta t\sum_{n=0}^{N-1}\nor{\Pi_{h}^{*}(\nabla s_{r}^{n+1})}{\tau}^{2}\nonumber \\
&\leq&C^{*}\left[\sigma+\frac{h^{2s}}{\nu^3}
+ \left(\frac{1}{\nu}+\frac{\alpha^2}{\sigma}\right)\left(h^{2s}+\Delta t^2\right)\left(1+\nor{S^v_r}{2}\right)+\frac{1}{\nu}\left(h^{2s}+\Delta t^2\right)\nor{S^p_r}{2}\right]\nonumber \\
&&+C^{*}\left[\left(\frac{1}{\nu}+\frac{\alpha^2}{\sigma}\right)\sum_{i=r+1}^{M_v}\nor{\nabla\boldsymbol{\varphi}_{i}}{{\bf L}^2}^{2}\lambda_i + \frac{1}{\nu}\sum_{i=r+1}^{M_p}\left(1+\nor{\nabla\psi_{i}}{{\bf L}^2}^{2}\right)\gamma_i\right].
\end{eqnarray}
where $C^{*}=Ce^{C_{1}\Delta t\sum_{n=0}^{N-1}\nor{\nabla\uv^{n+1}}{{\bf H}^{1}}^{2}}$, and $C_1$ is a positive constant depending on $\nu^{-1}$. Finally, using in \eqref{eq:ErrEqRiepF} the inequality:
$$
\max_{0\leq k\leq N}\nor{\boldsymbol{\phi}_{r}^{k}}{{\bf L}^{2}}^{2}\geq
C\Delta t\sum_{n=0}^{N-1}\nor{\boldsymbol{\phi}_{r}^{n+1}}{{\bf L}^{2}}^{2},
$$
triangle inequality, $L^2-H^1$ projection error estimates for velocity \eqref{eq:PODlmL2vel}-\eqref{eq:PODlmH01vel} and $L^2$ projection error estimate for pressure \eqref{eq:PODlmL2pres}, we get:
\begin{eqnarray}\label{eq:PODEEFend}
&&\Delta t\sum_{n=0}^{N-1}\nor{\uv^{n+1}-\uv_{r}^{n+1}}{{\bf L}^{2}}^{2}+\nu\Delta t\sum_{n=0}^{N-1}\nor{\nabla(\uv^{n+1}-\uv_{r}^{n+1})}{{\bf L}^{2}}^{2} + \sigma\Delta t\sum_{n=0}^{N-1}\nor{p^{n+1}-p_{r}^{n+1}}{L^{2}}^{2}\nonumber \\
&\leq& C^{*}\left[\sigma+\frac{h^{2s}}{\nu^3}
+ \left(\frac{1}{\nu}+\frac{\alpha^2}{\sigma}\right)\left(h^{2s}+\Delta t^2\right)\left(1+\nor{S^v_r}{2}\right)+\frac{1}{\nu}\left(h^{2s}+\Delta t^2\right)\nor{S^p_r}{2}\right]\nonumber \\
&&+C^{*}\left[\left(\frac{1}{\nu}+\frac{\alpha^2}{\sigma}\right)\sum_{i=r+1}^{M_v}\left(1+\nor{\nabla\boldsymbol{\varphi}_{i}}{{\bf L}^2}^{2}\right)\lambda_i + \frac{1}{\nu}\sum_{i=r+1}^{M_p}\left(1+\nor{\nabla\psi_{i}}{{\bf L}^2}^{2}\right)\gamma_i\right].
\end{eqnarray}
This concludes the proof.
\qed

\begin{remark}\label{rm:OptErrEst}
Note that in error estimate \eqref{eq:PODEEF} the convergence order is limited by the penalty constant $\sigma$. To remove this order limitation, one could simply replace bound \eqref{eq:divErrEq} in Theorem \ref{th:PODEE} by: 
\BEQ\label{eq:divErrEqOpt}
(\div\boldsymbol{\eta}^{n+1},s_{r}^{n+1})\leq \nor{\div\boldsymbol{\eta}^{n+1}}{L^2}\nor{s_{r}^{n+1}}{L^2}\leq \frac{\varepsilon_{2}^{-1}}{4}\nor{\nabla\boldsymbol{\eta}^{n+1}}{{\bf L}^2}^{2}
+\varepsilon_{2}\nor{s_{r}^{n+1}}{L^2}^{2},
\EEQ
applying Cauchy--Schwarz and Young's inequalities. Taking $\varepsilon_2 = \sigma/2$, this leads to the error estimate:
\begin{eqnarray}\label{eq:PODEEFopt}
&&\Delta t\sum_{n=0}^{N-1}\nor{\uv^{n}-\uv_{r}^{n}}{{\bf L}^{2}}^{2}+\nu\Delta t\sum_{n=0}^{N-1}\nor{\nabla(\uv^{n+1}-\uv_{r}^{n+1})}{{\bf L}^{2}}^{2} + \sigma\Delta t\sum_{n=0}^{N-1}\nor{p^{n+1}-p_{r}^{n+1}}{L^{2}}^{2}\nonumber \\
&\leq& C^{*}\left[\frac{h^{2s}}{\nu^3}
+ \left(\frac{1}{\nu}+\frac{1}{\sigma}\right)\left(h^{2s}+\Delta t^2\right)\left(1+\nor{S^v_r}{2}\right)+\frac{1}{\nu}\left(h^{2s}+\Delta t^2\right)\nor{S^p_r}{2}\right]\nonumber \\
&&+C^{*}\left[\left(\frac{1}{\nu}+\frac{1}{\sigma}\right)\sum_{i=r+1}^{M_v}\left(1+\nor{\nabla\boldsymbol{\varphi}_{i}}{{\bf L}^2}^{2}\right)\lambda_i + \frac{1}{\nu}\sum_{i=r+1}^{M_p}\left(1+\nor{\nabla\psi_{i}}{{\bf L}^2}^{2}\right)\gamma_i\right]. 
\end{eqnarray}
So, we have removed the order limitation due to the penalty constant $\sigma$, but changed the constant $(1/\nu+\alpha^2/\sigma)$ by $(1/\nu+1/\sigma)$. In any case, note that to prove error estimates, we have to assume to work with sufficiently regular flows, which is a common approach in the derivation of error estimates for POD-ROM (cf. \cite{IliescuWang14,Luo08}). Thus, the convergence orders obtained are generally valid in laminar flow settings or for sufficiently regular flows, but are usually not valid in realistic turbulent flow settings, since the convergence order decreases with the regularity. In this context, the kinematic viscosity is usually of the order $\nu=\mathcal{O}(10^{-3})$. Also, we have observed in the numerical studies performed in Section \ref{sec:NumStud}  that the saturation constant $\alpha$ starts with small values for small $r$ and seems to expe\-rience a flattening effect around $10^{-2}$ when adding more POD modes. Thus, for practical relatively coarse FE mesh size $h$ and time step $\Delta t$, we have that estimate \eqref{eq:PODEEF}, even if contains explicitly the penalty constant $\sigma$, improves estimate \eqref{eq:PODEEFopt}, since the presence of the saturation constant $\alpha$ allows to ease the convergence order reduction due to $\sigma^{-1}$.
\end{remark}

\begin{remark}\label{rm:ErrEstPrimPres}
Following Remark \ref{rm:PresEst}, an alternative unconditional pressure error estimate for LPS-ROM \eqref{eq:discTempROM} could be obtained if we set $\psi=0$ in the error equation \eqref{eq:ErrEq0}. Thus, we obtain for any $\boldsymbol{\varphi}\in\Xv_{r}$:
\begin{eqnarray*}
(s_{r}^{n+1},\div\boldsymbol{\varphi})&=&\frac{1}{\Delta t}(\boldsymbol{\phi}_{r}^{n+1}-\boldsymbol{\phi}_{r}^{n},\boldsymbol{\varphi})+\nu(\nabla\boldsymbol{\phi}_{r}^{n+1},\nabla\boldsymbol{\varphi})+b(\uv_{r}^{n+1},\boldsymbol{\phi}_{r}^{n+1},\boldsymbol{\varphi})\\
&&+b(\boldsymbol{\phi}_{r}^{n+1},\uv^{n+1},\boldsymbol{\varphi})-\langle\boldsymbol{\varepsilon}_{r}^{n+1},\boldsymbol{\varphi}\rangle,
\end{eqnarray*}
with $\boldsymbol{\varepsilon}_{r}^{n+1}$ denoting the consistency error, defined as:
\begin{equation*}
\langle\boldsymbol{\varepsilon}_{r}^{n+1},\boldsymbol{\varphi}\rangle=\langle\boldsymbol{c}^{n},\boldsymbol{\varphi}\rangle+\nu(\nabla\boldsymbol{\eta}^{n+1},\nabla\boldsymbol{\varphi})+b(\uv_{r}^{n+1},\boldsymbol{\eta}^{n+1},\boldsymbol{\varphi})+b(\boldsymbol{\eta}^{n+1},\uv^{n+1},\boldsymbol{\varphi})-(\rho^{n+1},\div\boldsymbol{\varphi}).
\end{equation*}
Let $S_{r}^{n+1}=\disp\sum_{k=0}^{n}\Delta t s_{r}^{k+1}=\disp\sum_{k=0}^{n}\Delta t (p_{r}^{k+1}-z_{r}^{k+1})$, then summation over the discrete times gives:
\begin{eqnarray*}
\left(S_{r}^{n+1},\div\boldsymbol{\varphi}\right)&=&(\boldsymbol{\phi}_{r}^{n+1}-\boldsymbol{\phi}_{r}^{0},\boldsymbol{\varphi})+\sum_{k=0}^{n}\Delta t\left[\nu(\nabla\boldsymbol{\phi}_{r}^{n+1},\nabla\boldsymbol{\varphi})-\langle\boldsymbol{\varepsilon}_{r}^{n+1},\boldsymbol{\varphi}\rangle\right]\\
&&+\sum_{k=0}^{n}\Delta t\left[b(\uv_{r}^{n+1},\boldsymbol{\phi}_{r}^{n+1},\boldsymbol{\varphi})+b(\boldsymbol{\phi}_{r}^{n+1},\uv^{n+1},\boldsymbol{\varphi})\right].
\end{eqnarray*}
Thus, we get:
\begin{eqnarray}\label{eq:ModInfSupEE}
&&\sup_{\boldsymbol{\varphi}\in\Xv_r}\frac{(S_{r}^{n+1},\div\boldsymbol{\varphi})}{\nor{\nabla\boldsymbol{\varphi}}{{\bf L}^2}}+\nor{\Pi_{h}^{*}(\nabla S_{r}^{n+1})}{\tau}\nonumber \\
&\leq& C\left[\nor{\boldsymbol{\phi}_{r}^{n+1}}{{\bf L}^{2}}+\nor{\boldsymbol{\phi}_{r}^{0}}{{\bf L}^{2}}+\sum_{k=0}^{n}\Delta t\left(\nu\nor{\nabla\boldsymbol{\phi}_{r}^{k+1}}{{\bf L}^{2}}+\nor{\boldsymbol{\varepsilon}_{r}^{k+1}}{{\bf H}^{-1}}\right)\right]\nonumber \\
&&+C\sum_{k=0}^{n}\Delta t\left(\nor{\nabla\uv_{r}^{k+1}}{{\bf L}^2}\nor{\nabla\boldsymbol{\phi}_{r}^{k+1}}{{\bf L}^2}+\nor{\nabla\boldsymbol{\phi}_{r}^{k+1}}{{\bf L}^2}\nor{\nabla\uv^{k+1}}{{\bf L}^2}\right)\nonumber \\
&&+\nor{\Pi_{h}^{*}(\nabla S_{r}^{n+1})}{\tau}.
\end{eqnarray}
where we have applied triangle inequality, Cauchy--Schwarz inequality, the definition of the dual norm, and the standard estimate \eqref{eq:formb1} for the convective term. Then, using the norm defined in \eqref{eq:NewNorm}, by Cauchy--Schwarz inequality, the stability result \eqref{eq:stabEst} for the reduced order velocity and the regularity assumption \ref{hp:ContRegularity} on the continuous velocity, from \eqref{eq:ModInfSupEE} we have:
\begin{eqnarray}\label{eq:ModInfSup1EE}
|||S_{r}^{n+1}|||&\leq\disp\frac{C}{\nu}&\left[\max_{k=0,\ldots,N}\nor{\boldsymbol{\phi}_{r}^{k}}{{\bf L}^{2}}+\left(\sum_{k=0}^{N-1}\Delta t\nor{\nabla\boldsymbol{\phi}_{r}^{k+1}}{{\bf L}^{2}}^{2}\right)^{1/2}+\left(\sum_{k=0}^{N-1}\Delta t\nor{\boldsymbol{\varepsilon}_{r}^{k+1}}{{\bf H}^{-1}}^{2}\right)^{1/2}\right. \nonumber \\
&&+\left.\left(\sum_{k=0}^{N-1}\Delta t\nor{\Pi_{h}^{*}(\nabla s_{r}^{k+1})}{\tau}^{2}\right)^{1/2}\right].
\end{eqnarray}
Using estimates \eqref{eq:cn}-\eqref{eq:eta4} to bound the third term on \eqref{eq:ModInfSup1EE}, and estimate \eqref{eq:ErrEqRiepF} to bound the rest of terms in \eqref{eq:ModInfSup1EE}, we obtain:
\BEQ\label{eq:ModInfSup2EE}
|||S_{r}^{n+1}|||\leq\frac{\sqrt{C^{*}E}}{\nu^{3/2}},
\EEQ
where $E$ denotes the sum of all terms within brackets in \eqref{eq:PODEEF}. From triangle inequality and the $L^2$ projection error estimates for pressure \eqref{eq:PODlmL2pres}, it follows:
\BEQ\label{eq:ModInfSup3EE}
\nor{P-P_{r}}{\ell^{\infty}(\widetilde{Q}_{r})}:=\max_{n=1,\ldots,N}|||P^{n}-P_{r}^{n}|||\leq \frac{\sqrt{C^{*}E}}{\nu^{3/2}},
\EEQ
and we have denoted $P(\cdot,t)=\disp\int_{0}^{t}\widetilde{p}(\cdot,s)\,ds$, being $\widetilde{p}$ the piecewise constant in time function that takes the value $p^{n+1}$ on $(t_{n},t_{n+1})$. Thus, we have derived a new error estimate on the time-average of the reduced order pressure with respect to the one obtained in Theorem \ref{th:PODEE}:
\BEQ\label{eq:stabPresEstEE}
\nor{p-p_{r}}{\ell^{2}(L^2)}:=\left[\sum_{n=1}^{N}\Delta t\nor{p^{n}-p_{r}^{n}}{L^2}^{2}\right]^{1/2}\leq \frac{\sqrt{C^{*}E}}{\sigma^{1/2}}.
\EEQ
\end{remark}

\section{Numerical studies}\label{sec:NumStud}

In this section, we present numerical results for the LPS-ROM introduced and analyzed in the previous sections, for which the standard discrete inf-sup condition is circumvented and neither strongly nor weakly divergence-free POD modes are required. The numerical experiments are performed on the benchmark problem of the 2D laminar unsteady flow around a cylinder with circular cross-section \cite{SchaferTurek96}. The open-source FE software FreeFEM \cite{Hecht12} has been used to run the numerical experiments.

\subsection{Setup for numerical simulations}

\begin{figure}[htb]
\begin{center}
\centerline{\includegraphics[scale=1.]{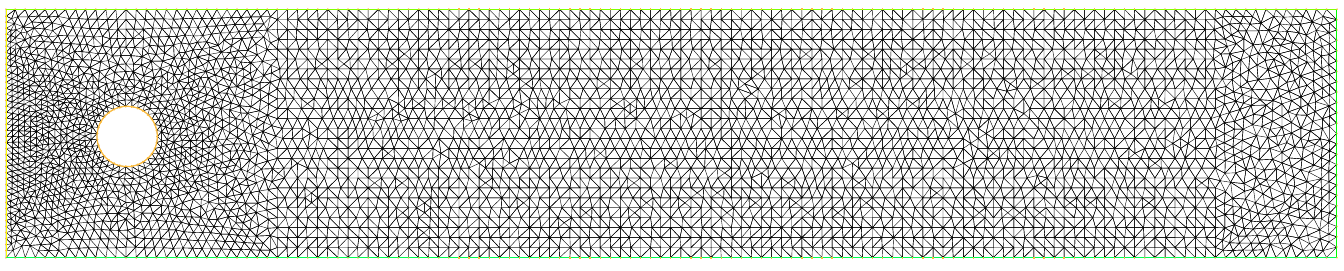}}
\caption{Computational grid.}\label{fig:Mesh}
\end{center}
\end{figure}

Following \cite{SchaferTurek96}, the computational domain is given by a rectangular channel with a circular hole (see Figure \ref{fig:Mesh} for the computational grid used):
$$
\Om=\{(0,0.2)\times(0,0.41)\}\backslash \{\xv : (\xv-(0.2,0.2))^2 \leq 0.05^2\}.
$$

No slip boundary conditions are prescribed on the horizontal walls and on the cylinder, and a parabolic inflow profile is provided at the inlet:
$$
\uv(0,y,t)=(4U_{m}y(H-y)/H^2, 0)^{T},
$$
with $U_{m}=\uv(0,H/2,t)=1.5\,\rm{m/s}$, and $H=0.41\,\rm{m}$ the channel height. At the outlet, we perform a comparison using on one side outflow (do nothing) boundary conditions $(\nu\nabla\uv - p\,Id)\nv={\bf 0}$, with $\nv$ the outward normal to the domain, for which we can remove the penalty term with factor $\sigma$ to the variational formulation, since the constant the pressure is determined up through the formulation is already fixed by the prescribed boundary conditions. On the other side, we impose the same parabolic inflow profile for the outflow velocity, for which we use the penalty term with factor $\sigma$ to fix the constant the pressure is determined up through the formulation. Although imposing Dirichlet boundary conditions at the outlet is unphysical, we have seen that it simplifies the theoretical analysis performed. Also, we have observed in the numerical studies that this do not influence too much quantities of interest such as the upstream drag and lift coefficients, and it is often used \cite{John04paper,RebholzXiao17}, also in the ROM framework \cite{RebholzIliescu17,RebholzIliescu18}.

\medskip

\begin{figure}[htb]
\begin{center}
\includegraphics[width=4.8in]{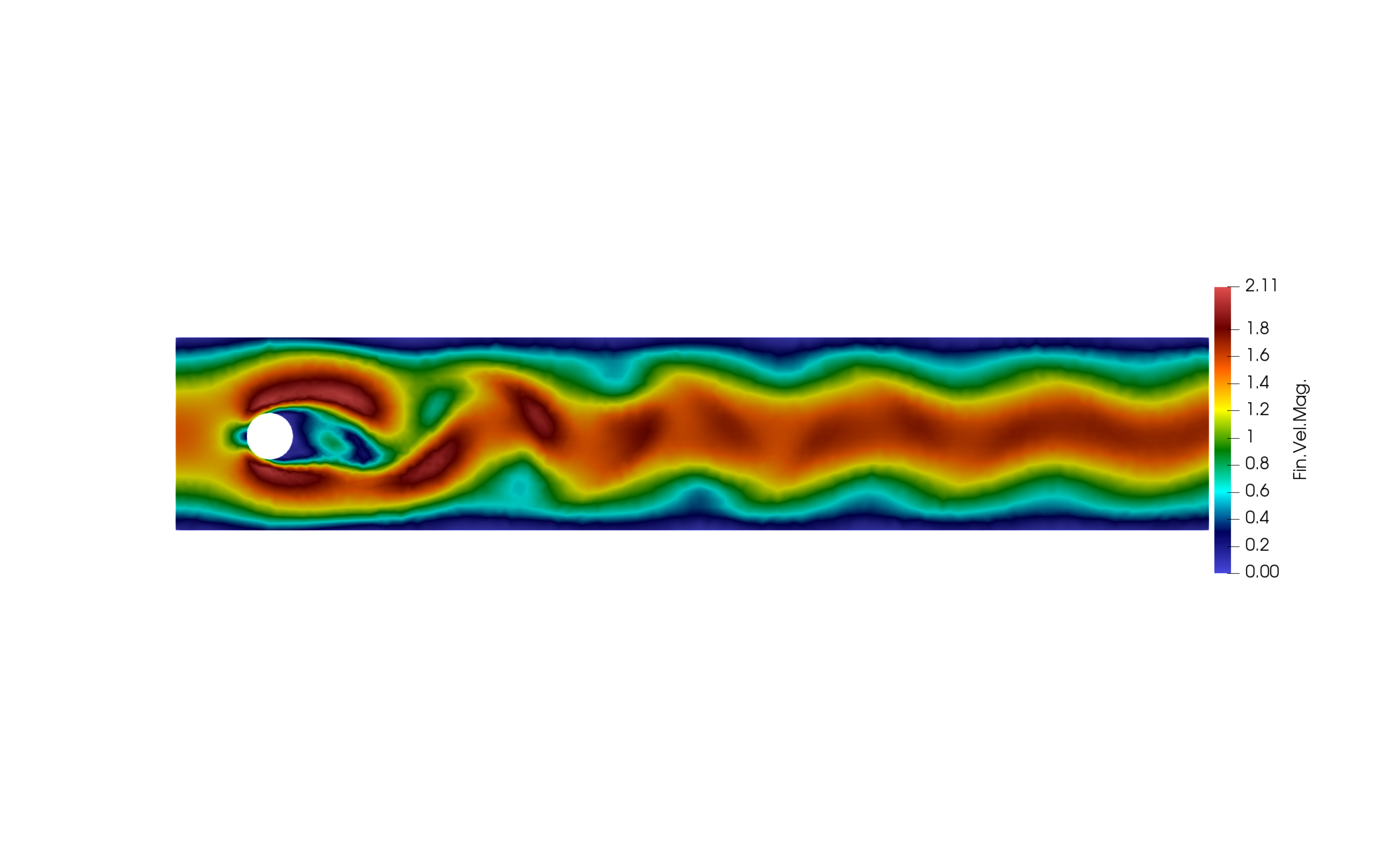}
\includegraphics[width=4.8in]{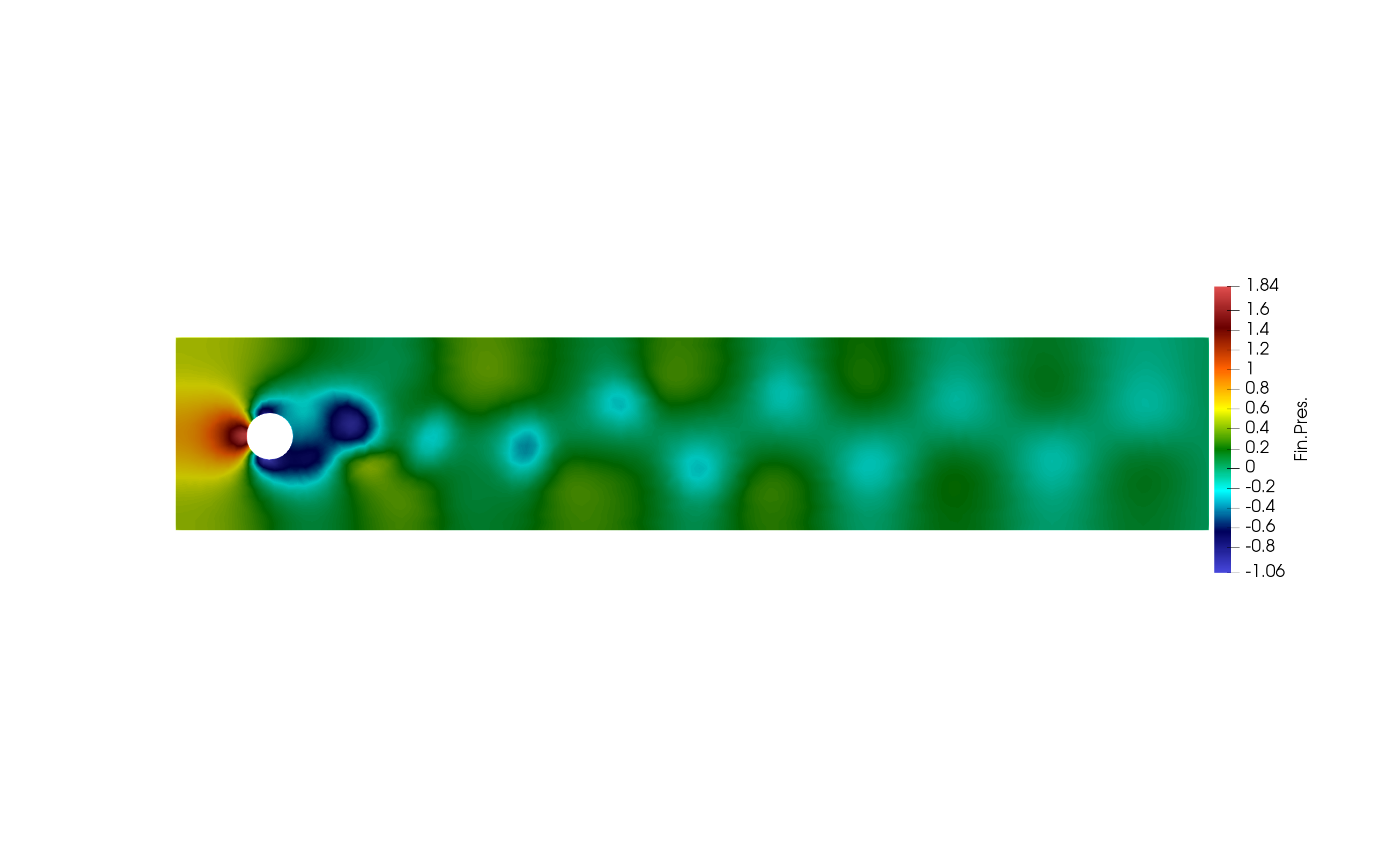}
\includegraphics[width=4.8in]{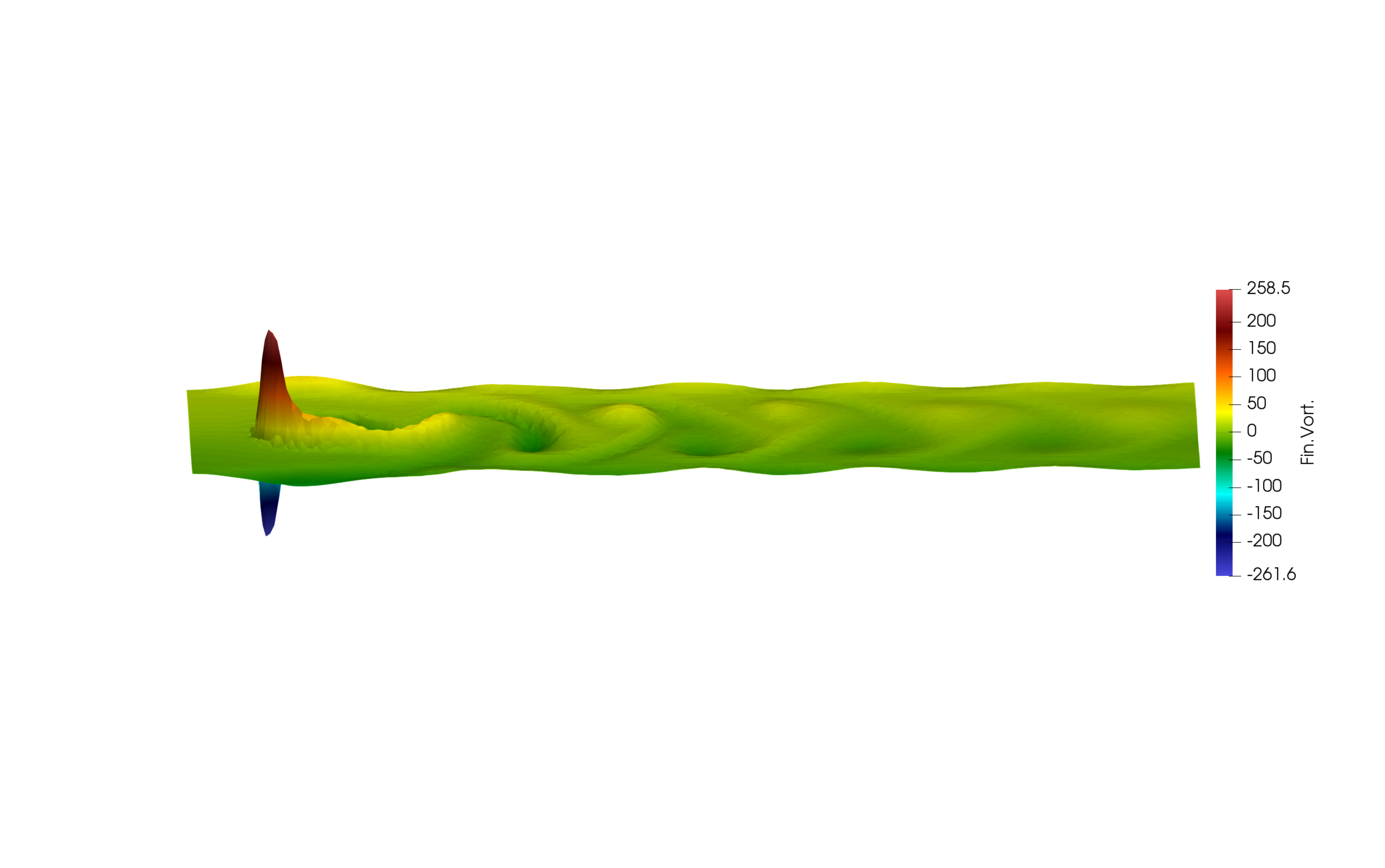}
\caption{Final FOM solution: velocity magnitude, pressure and vorticity (3D plot) from top to bottom with do nothing BC at the outlet.}\label{fig:FinFOMSol}
\end{center}
\end{figure}

\begin{figure}[htb]
\begin{center}
\includegraphics[width=4.8in]{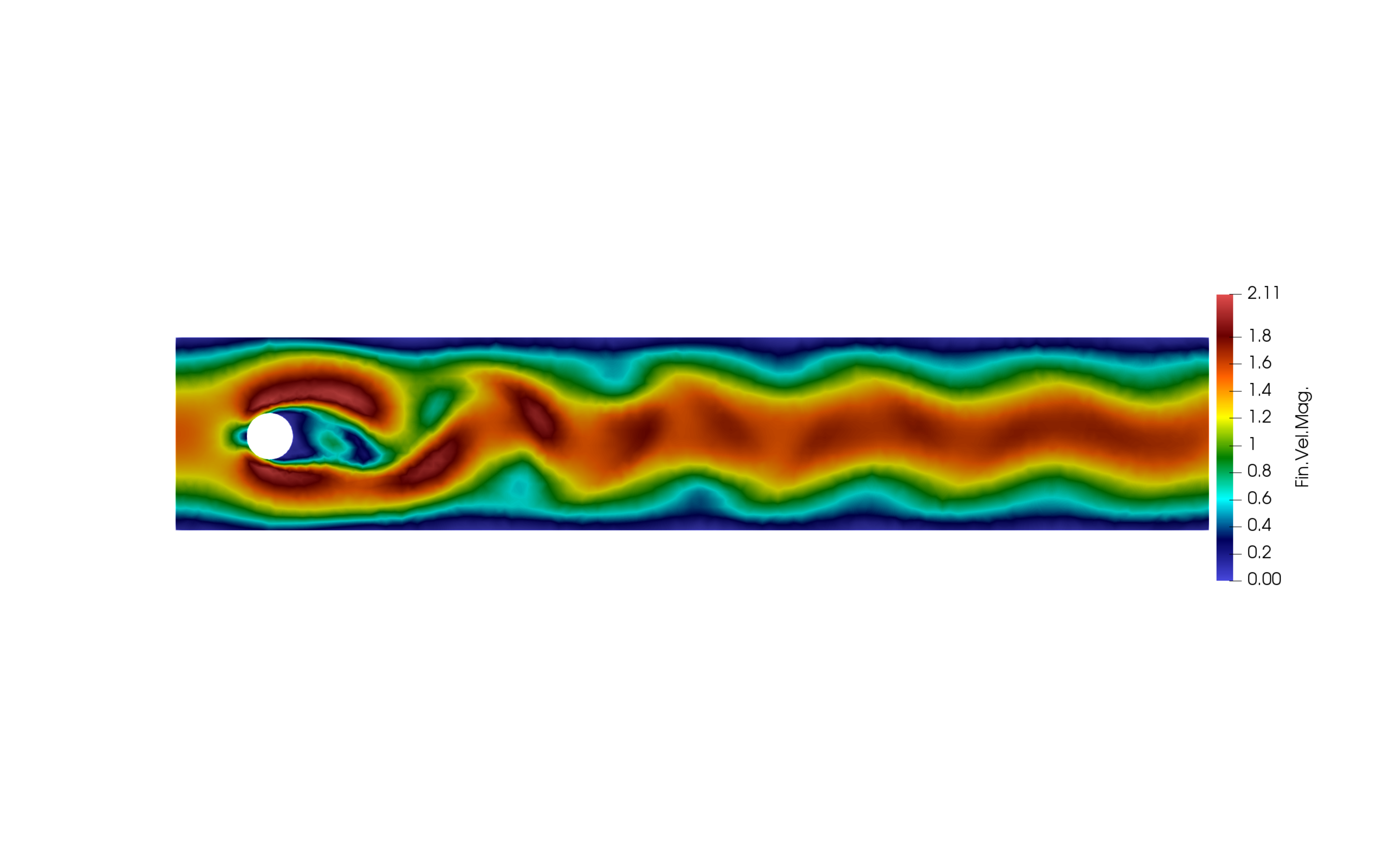}
\includegraphics[width=4.8in]{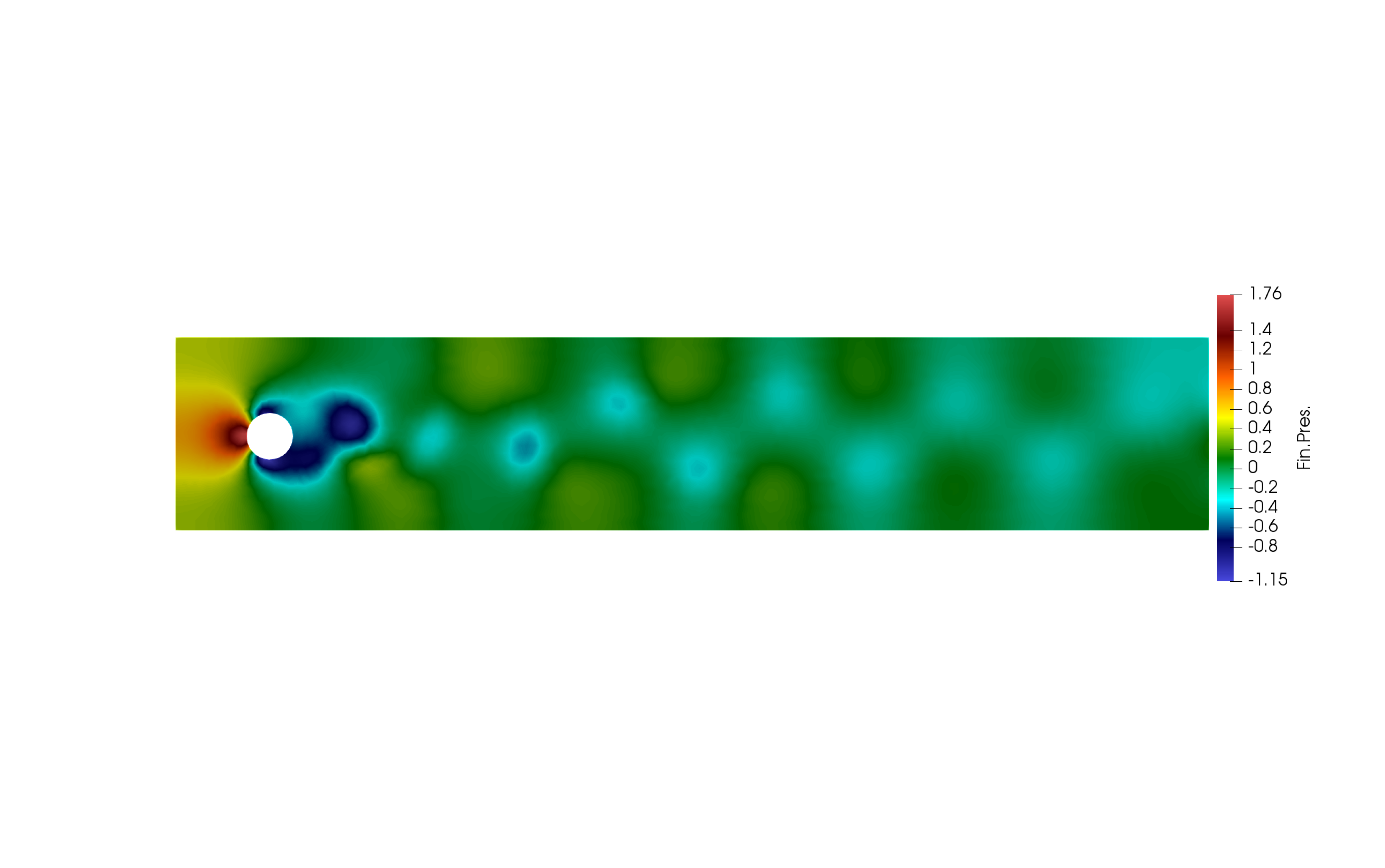}
\includegraphics[width=4.8in]{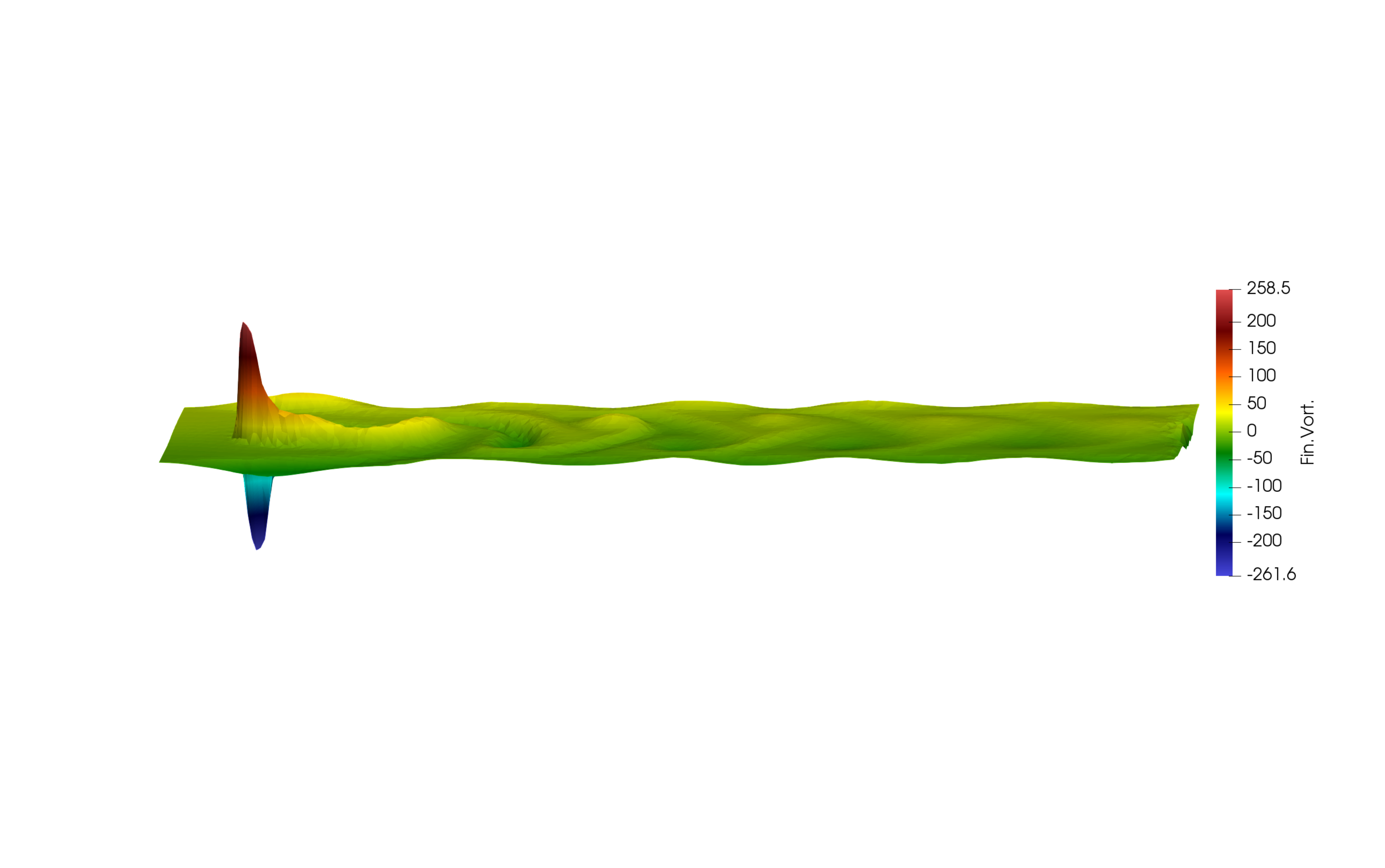}
\caption{Final FOM solution: velocity magnitude, pressure and vorticity (3D plot) from top to bottom with Dirichlet BC at the outlet.}\label{fig:FinFOMSolDir}
\end{center}
\end{figure}

The kinematic viscosity of the fluid is $\nu=10^{-3}\,\rm{m^{2}/s}$, and there is no external (gravity) forcing, i.e. $\fv={\bf 0}\,\rm{m/s^2}$. Based on the mean inflow velocity $\overline{U}=2U_{m}/3=1\,\rm{m/s}$, the cylinder diameter $D=0.1\,\rm{m}$ and the kinematic viscosity of the fluid $\nu=10^{-3}\,\rm{m^{2}/s}$, the Reynolds number considered is $Re=\overline{U}D/\nu=100$. In the fully developed periodic regime, a vortex shedding can be observed behind the obstacle, resulting in the well-known von K\'arm\'an vortex street (see Figures \ref{fig:FinFOMSol}-\ref{fig:FinFOMSolDir}).

\medskip

For the evaluation of computational results, we are interested in studying the temporal evolution of the following
quantities of interest. The kinetic energy of the flow is the most frequently monitored quantity, given by:
$$
E_{Kin}=\frac{1}{2}\nor{\uv}{{\bf L}^2}^2.
$$
Other relevant quantities of interest are the drag and lift coefficients. In order to reduce the boundary approximation influences, in the present work these quantities are computed as volume integrals:
$$
c_{D}=-\frac{2}{D\overline{U}^2}\left[(\partial_t \uv,\vv_D)+b(\uv,\uv,\vv_D)+\nu(\nabla\uv,\nabla\vv_D)-(p,\div\vv_D)\right],
$$
$$
c_{L}=-\frac{2}{D\overline{U}^2}\left[(\partial_t \uv,\vv_L)+b(\uv,\uv,\vv_L)+\nu(\nabla\uv,\nabla\vv_L)-(p,\div\vv_L)\right],
$$
for arbitrary test functions $\vv_{D},\vv_{L}\in {\bf H}^{1}$ such that $\vv_{D}=(1,0)^{T}$ on the boundary of the cylinder and vanishes on the other boundaries, $\vv_{L}=(0,1)^{T}$ on the boundary of the cylinder and vanishes on the other boundaries. In the actual computations, we have used the approach described in \cite{JohnMatthies01} to fix the test functions $\vv_D,\vv_L$ and evaluate the drag and lift coefficients $c_D,c_L$. Reference intervals for these coefficients were given in \cite{SchaferTurek96} (see last row of Table \ref{tab:DragLiftCoef}), together with the Strouhal number $St=Df/\overline{U}$, where $f$ is the frequency of the vortex shedding.

\begin{table}[htb]
$$\hspace{-0.1cm}
\begin{tabular}{|c|c|c|c|}
\hline
 & $c_{D}^{max}$ & $c_{L}^{max}$ & $St$\\
\hline
FOM (Outlet do nothing BC) & $3.22$ & $0.96$ & $0.303$\\
\hline
FOM (Outlet Dirichlet BC) & $3.16$ & $0.92$ & $0.303$\\
\hline
Reference results from \cite{SchaferTurek96} & $[3.22,3.24]$ & $[0.99,1.01]$ & $[0.295,0.305]$\\
\hline
\end{tabular}$$\caption{Maximum drag coefficient $c_{D}^{max}$, maximum lift coefficient $c_{L}^{max}$, and Strouhal number for the FOM solution with do nothing BC (first row) and Dirichlet BC (second row) at the outlet, compared with reference intervals from \cite{SchaferTurek96} (third row).}\label{tab:DragLiftCoef}
\end{table}

\subsection{LPS-FOM for snapshots generation}

\begin{figure}[htb]
\begin{center}
\includegraphics[width=4.8in]{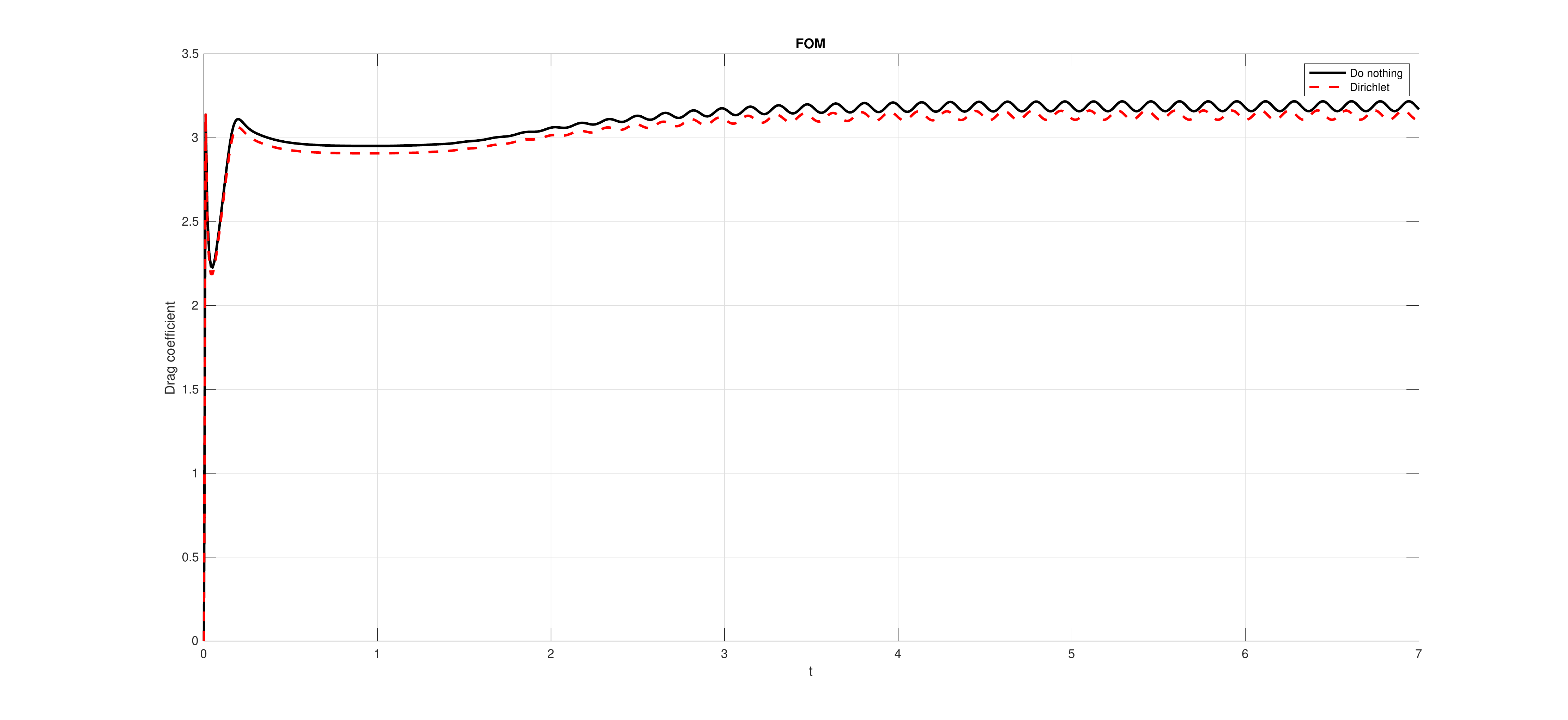}
\includegraphics[width=4.8in]{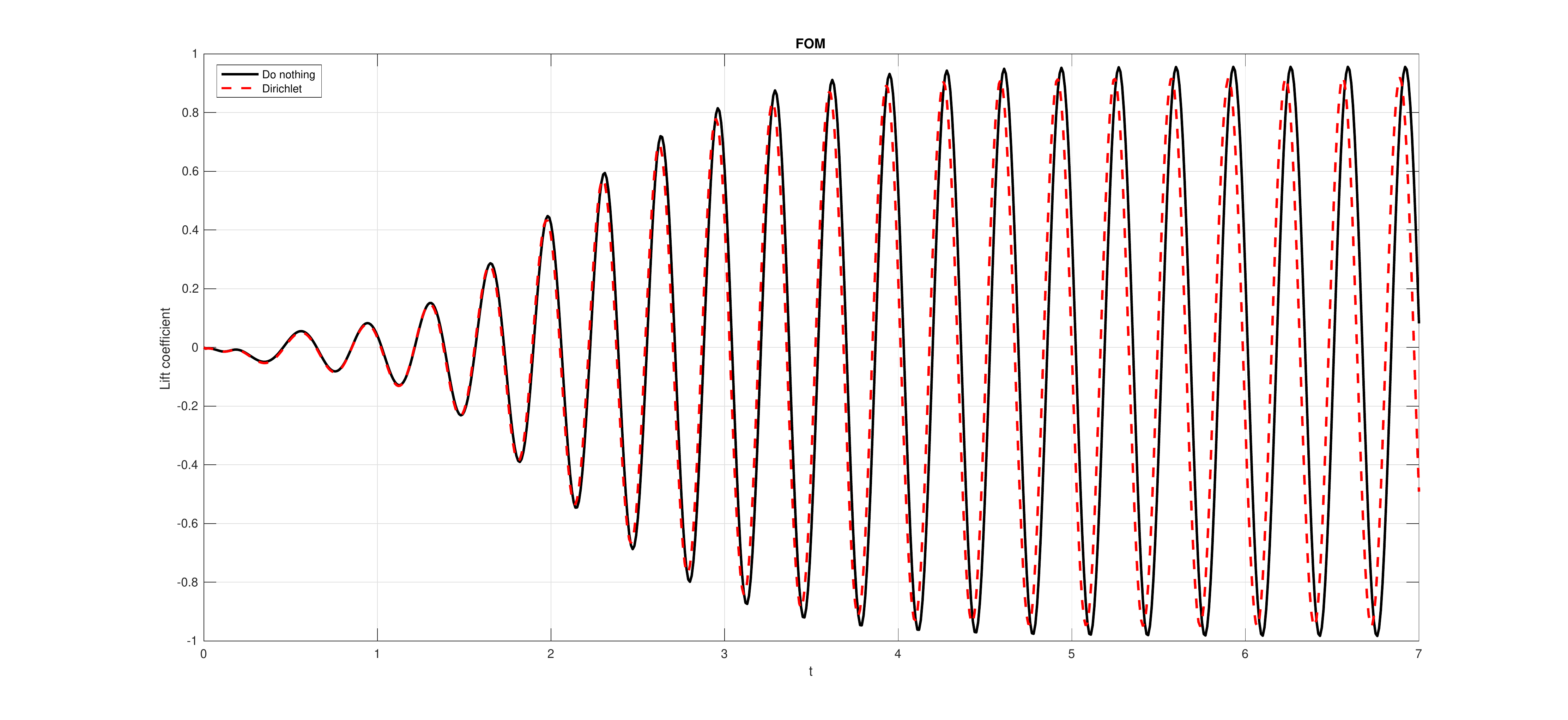}
\includegraphics[width=4.8in]{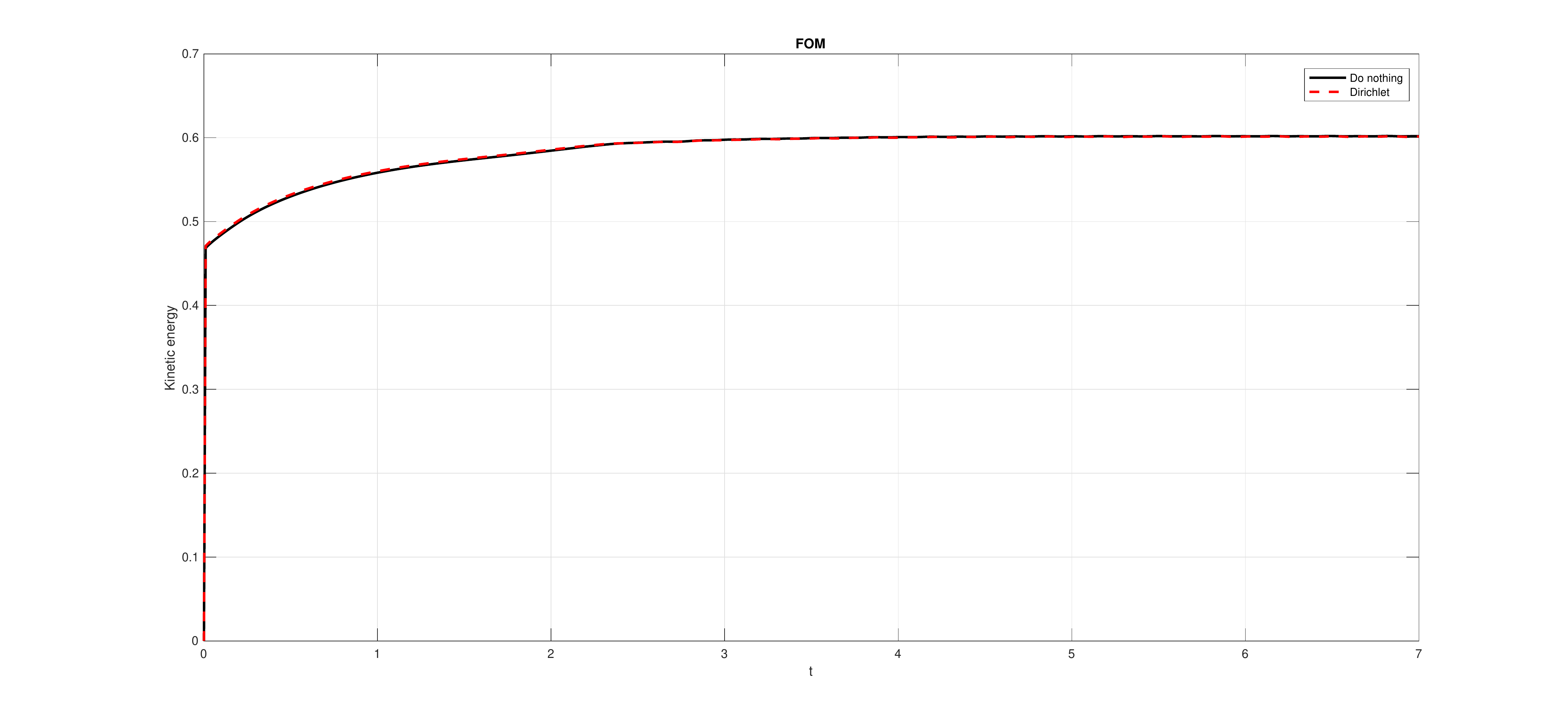}
\includegraphics[width=4.8in]{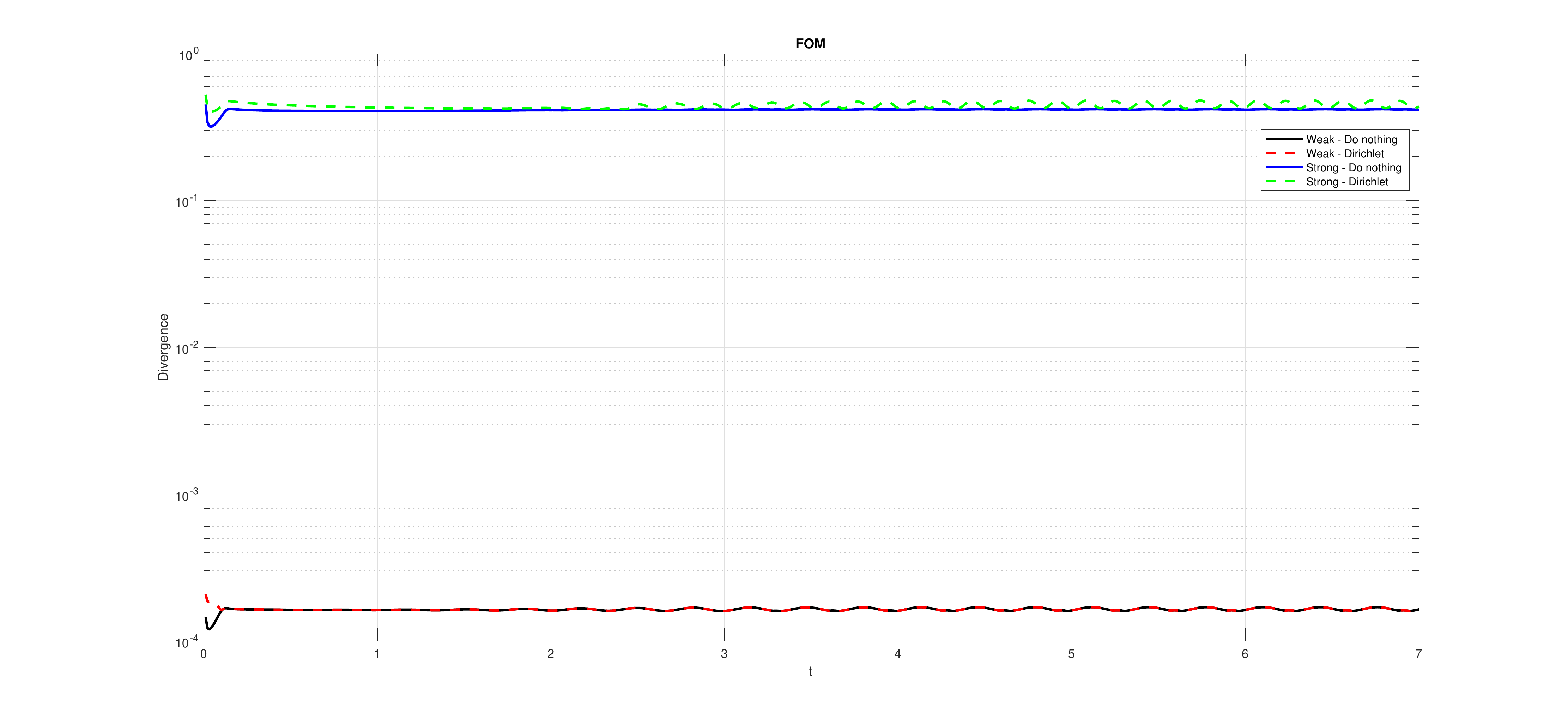}
\caption{Temporal evolution of drag coefficient, lift coefficient, kinetic energy and ``weak-strong'' divergence (from top to bottom) for the FOM solution with do nothing BC and Dirichlet BC at the outlet.}\label{fig:DragLiftCoef}
\end{center}
\end{figure}

The numerical method used to compute the snapshots is the LPS-FOM described in Section \ref{sec:FOM}, with a spatial discretization using EO ${\bf P}^{2}-\mathbb{P}^2$ FE for the pair velocity-pressure on a relatively coarse computational grid (see Figure \ref{fig:Mesh}, $h=2.76\times 10^{-2}\,\rm{m}$), resulting in $32\,488$ d.o.f. for velocities and $16\,244$ d.o.f. for pressure. We perform a comparison using do nothing Boundary Conditions (BC) at the outlet, for which $\sigma=0$, and Dirichlet BC at the outlet, for which $\sigma=10^{-6}$. The following expression of the stabilization parameters is used in the computations:
$$
\tau_{K}=\left(\frac{4}{\Delta t^2} + 32\frac{\nu^2}{(h_{K}/2)^4} + 4\frac{\overline{U}}{(h_{K}/2)^2}\right)^{-1/2},
$$
by adapting the form proposed in \cite{Codina02,CodinaBlasco02}, designed by a specific Fourier analysis applied in the framework of
stabilized methods. For the time discretization, a semi-implicit Backward Differentiation Formula of order 2 (BDF2) has been applied, which guarantees a good balance between numerical accuracy and computational complexity (cf. \cite{AhmedRubino19}). In particular, the discrete time derivative has been approximated by the operator $D_{t}^{2}$ defined as: 
$$
D_{t}^{2}\uhv^{n+1}=\frac{3\uhv^{n+1}-4\uhv^{n}+\uhv^{n-1}}{2\Delta t},\quad n\geq 1,
$$
and we have considered the following extrapolation for the convection velocity by means of Newton–Gregory backward polynomials (cf. \cite{Cellier91}): $\widehat{\uv}_{h}^{n}=2\uhv^{n}-\uhv^{n-1}$, $n\geq 1$, in order to achieve a second-order accuracy in time. For the initialization $(n=0)$, we have considered $\uhv^{-1}=\uhv^{0}=\uv_{0h}$, being $\uv_{0h}$ the initial condition, so that the time scheme reduces to the semi-implicit Euler method for the first time step $(\Delta t)^{0}=(2/3)\Delta t$. In the FOM simulations, an impulsive start is performed, i.e. the initial condition is a zero velocity field, and the time step is $\Delta t = 2\times 10^{-3}\,\rm{s}$. Time integration is performed till a final time $T=7\,\rm{s}$. In the time period $[0,5]\,\rm s$, after an initial spin-up, the flow is expected to develop to full extent, including a subsequent relaxation time. Afterwards, it reaches a periodic-in-time (statistically- or quasi-steady) state (see Figure \ref{fig:DragLiftCoef}, where we plot in particular drag and lift coefficients, and kinetic energy temporal evolution for the FOM solutions). 

\medskip

Observe from Figures \ref{fig:FinFOMSol}-\ref{fig:FinFOMSolDir} that results are qualitatively very close using do-nothing and Dirichlet BC at the outlet. The only very little difference can be noticed in the vorticity plot, where some numerical oscillations of the vorticity field appear when using Dirichlet BC at the outlet. Nevertheless, they have low magnitude and just slightly influence drag and lift coefficients reported in Table \ref{tab:DragLiftCoef}. Indeed, results with do nothing BC for all quantities of interest agree quite well with reference results from \cite{SchaferTurek96} and other numerical studies \cite{IliescuJohn14,RebholzIliescu17,RebholzIliescu18}, and are only slightly less accurate when using Dirichlet BC at the outlet. This behavior is reflected also in Figure \ref{fig:DragLiftCoef}, where in addition we have also plotted the ``weak-strong'' divergence temporal evolution for the FOM solution. In parrticular, the curve of the ``strong'' divergence has been obtained by plotting $\nor{\div\uhv^{n}}{L^2}$, while the curve of the ``weak'' divergence has been obtained by plotting $\max_{q_h \in Q_h}|(\div\uhv^n,q_h)|$. From this figure, it is evident that the computed snapshots are neither strongly nor weakly divergence-free.

\medskip

\begin{figure}[htb]
\begin{center}
\includegraphics[width=3.6in]{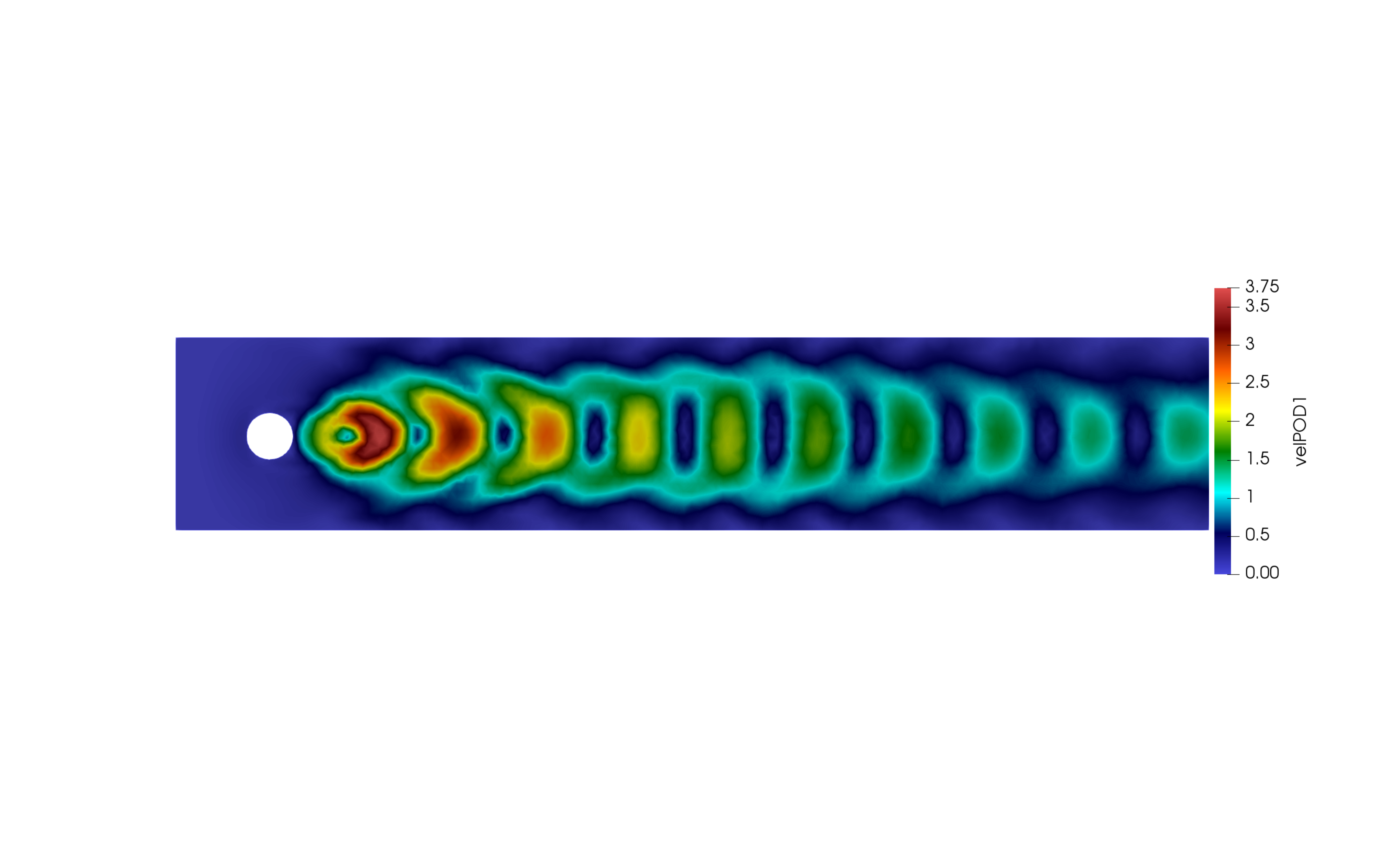}
\includegraphics[width=3.6in]{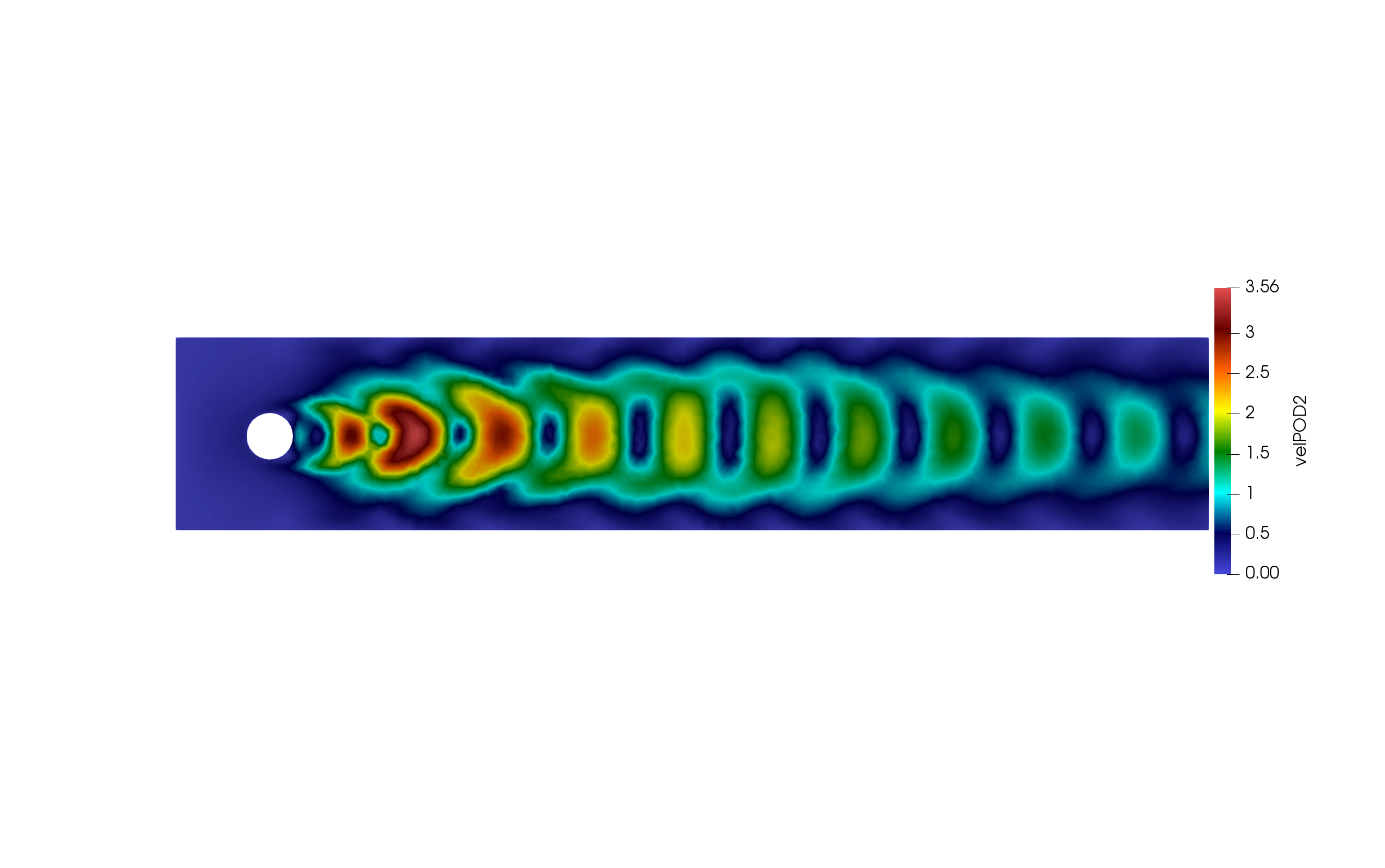}
\includegraphics[width=3.6in]{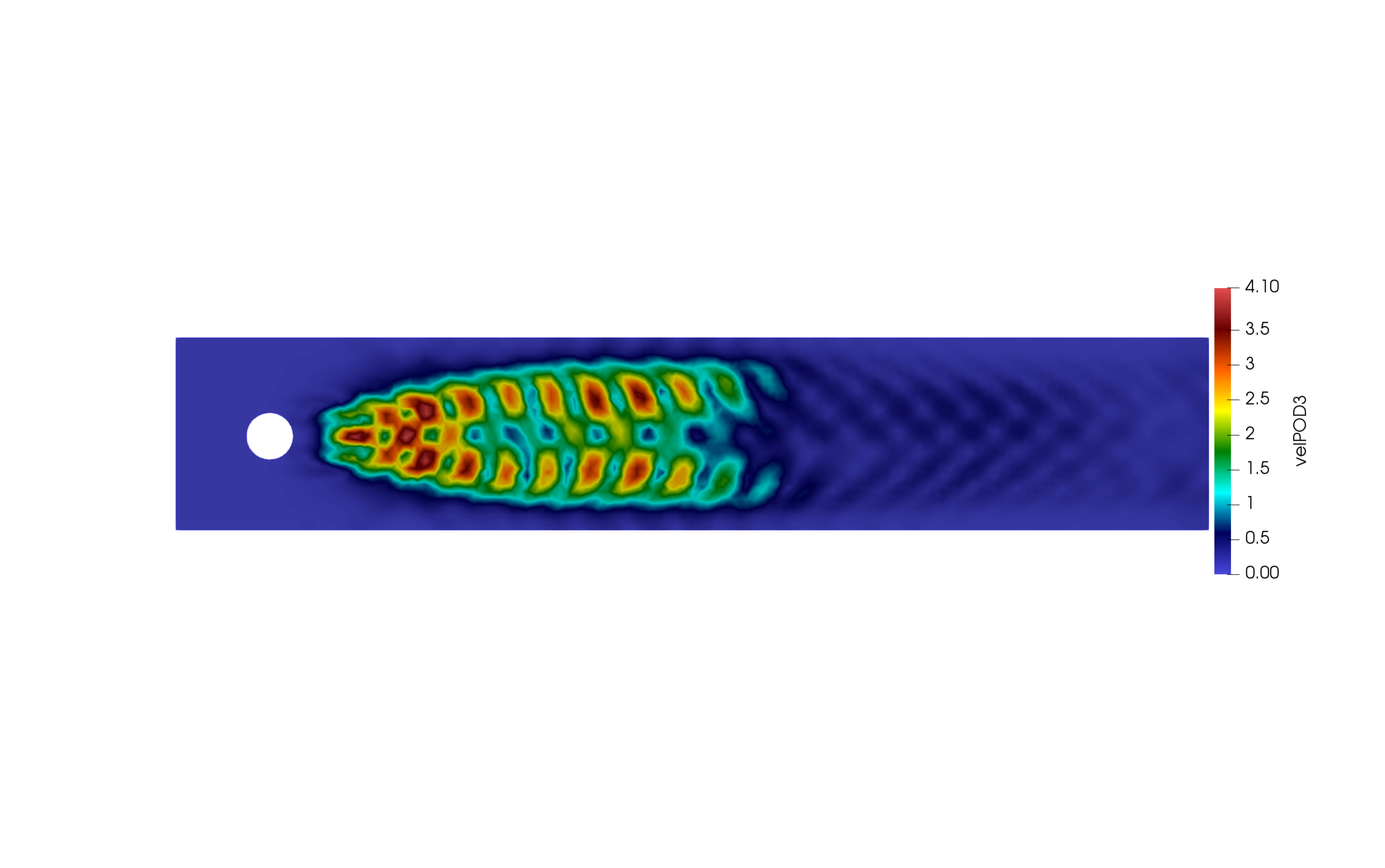}
\includegraphics[width=3.6in]{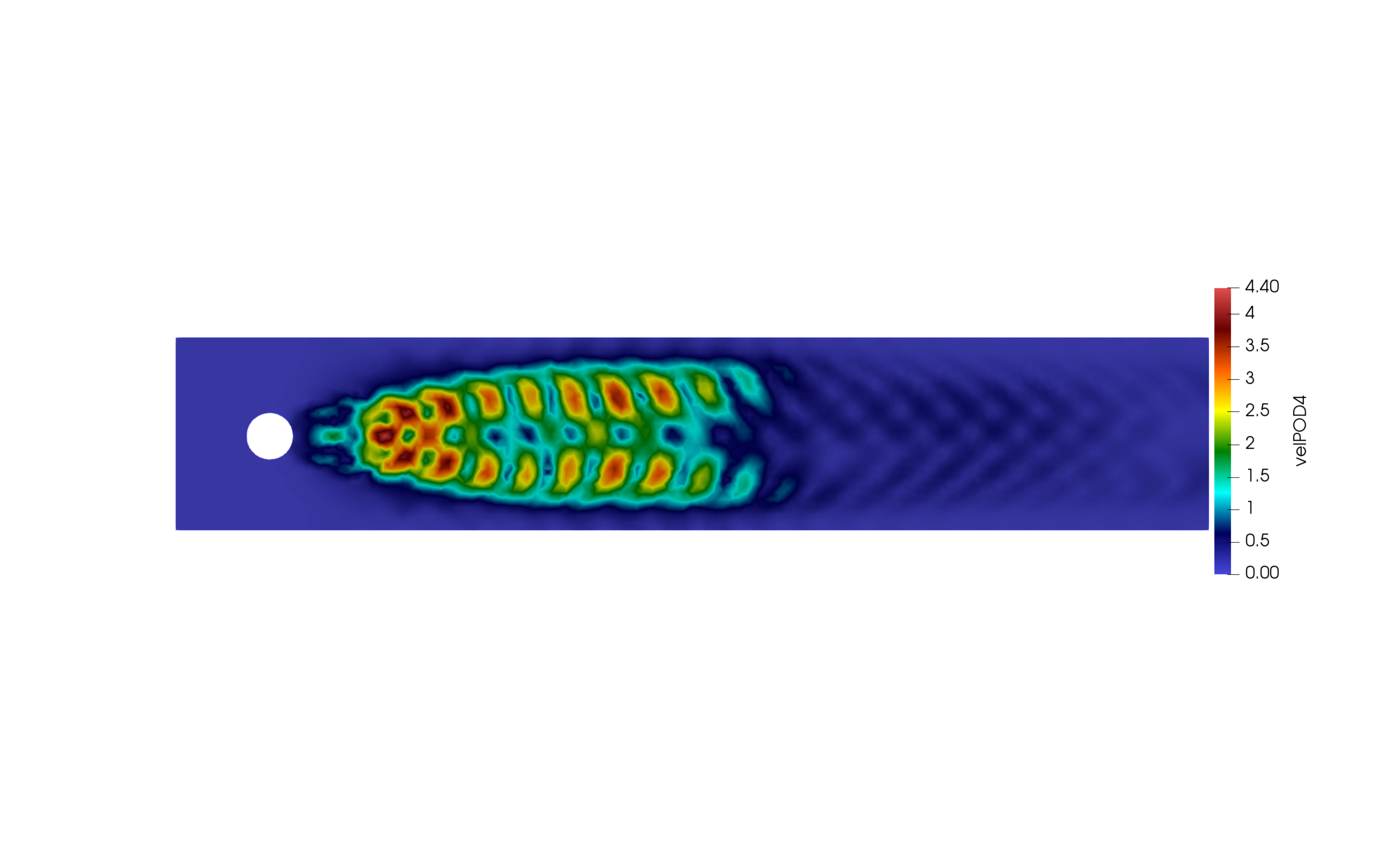}
\caption{First POD velocity modes (Euclidean norm): do nothing BC at the outlet.}\label{fig:velPODmodes}
\end{center}
\end{figure}

\begin{figure}[htb]
\begin{center}
\includegraphics[width=3.6in]{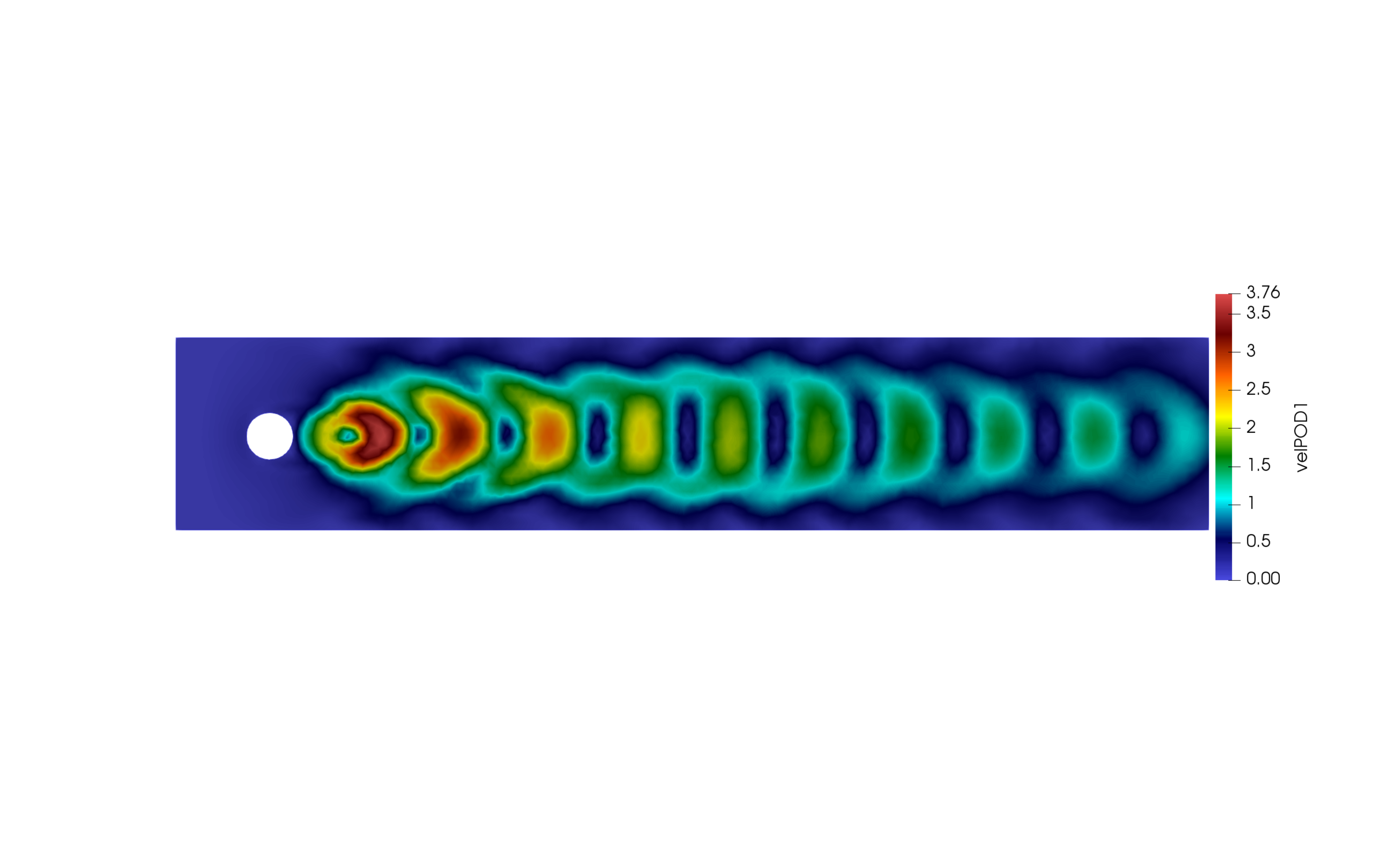}
\includegraphics[width=3.6in]{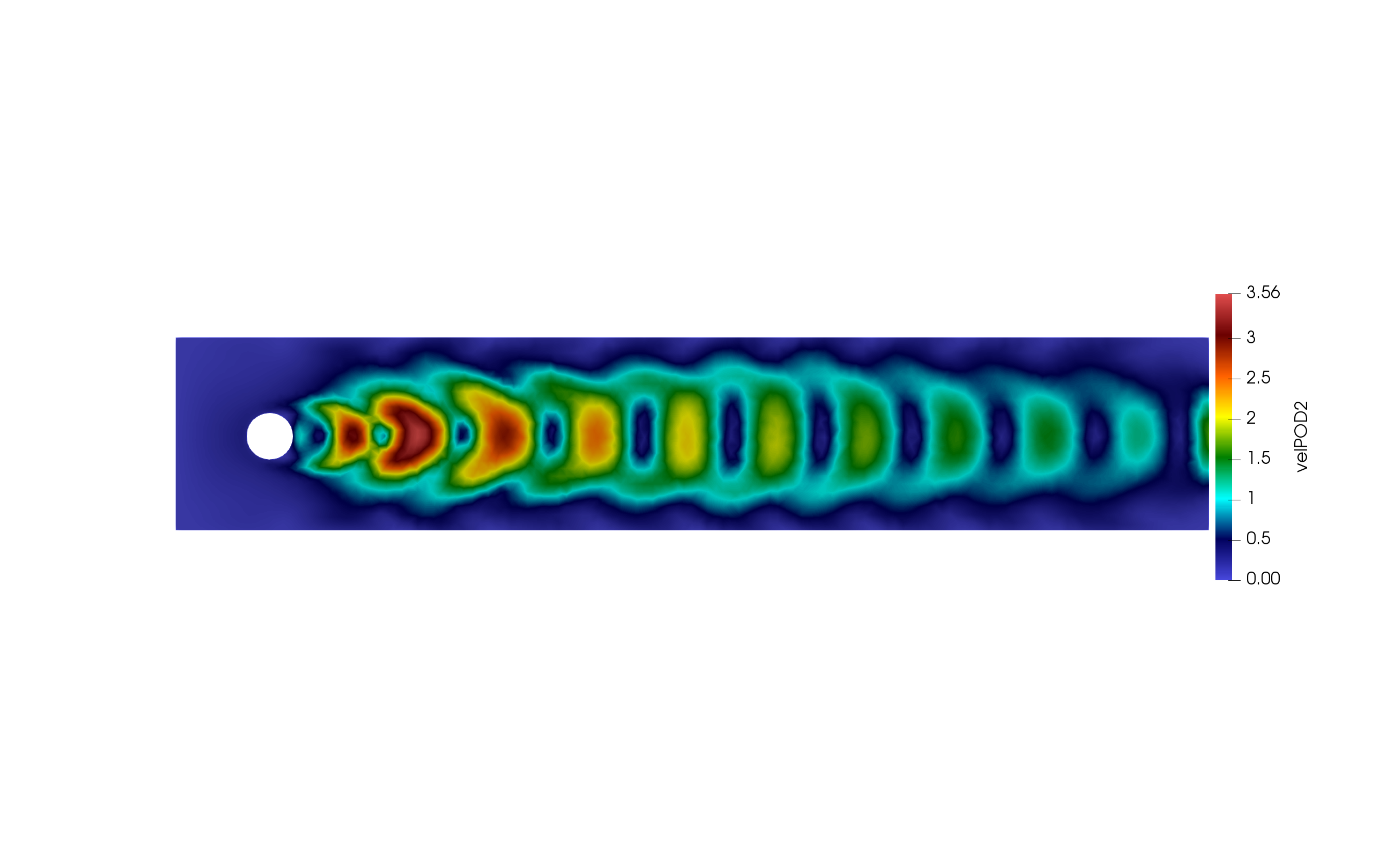}
\includegraphics[width=3.6in]{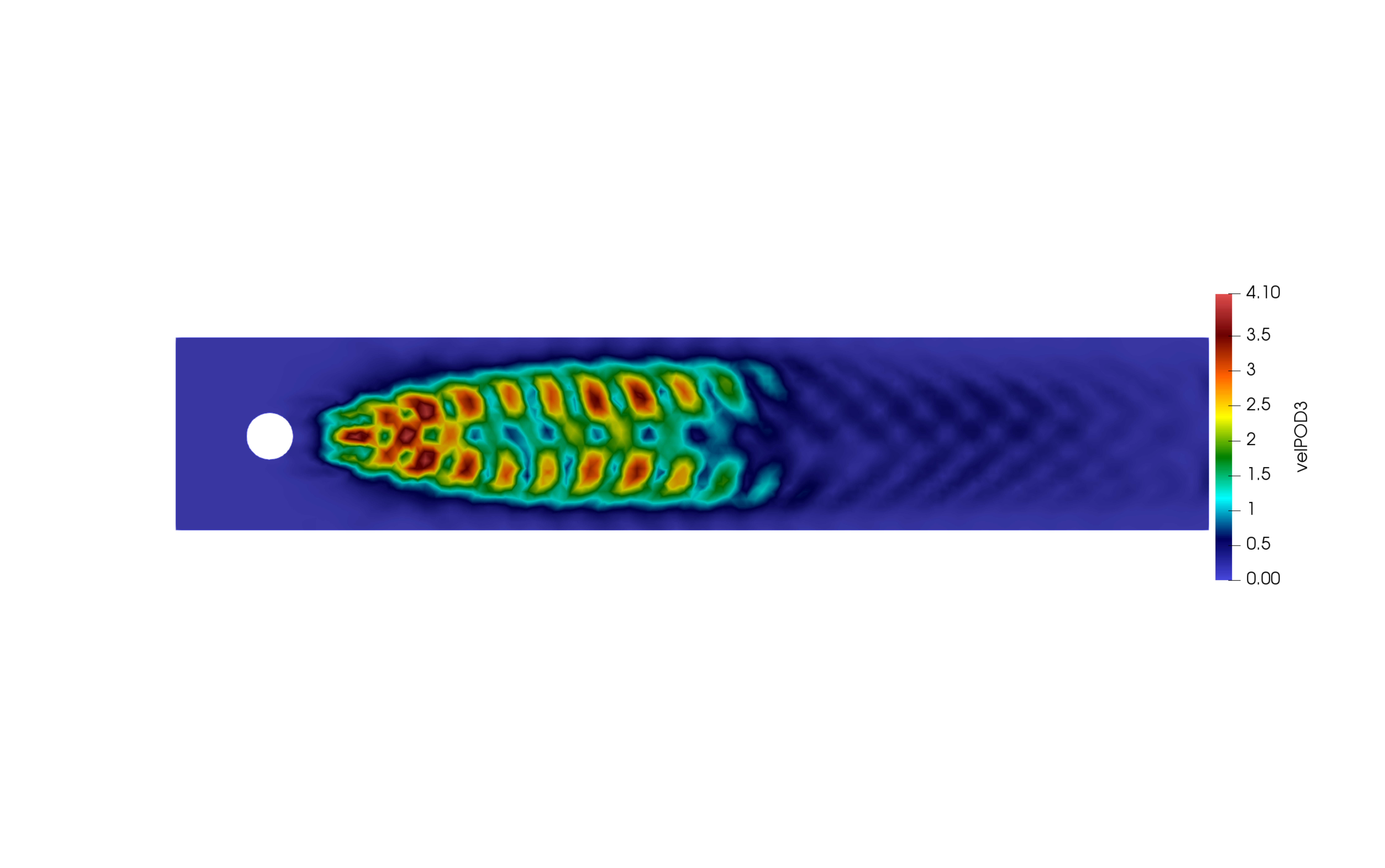}
\includegraphics[width=3.6in]{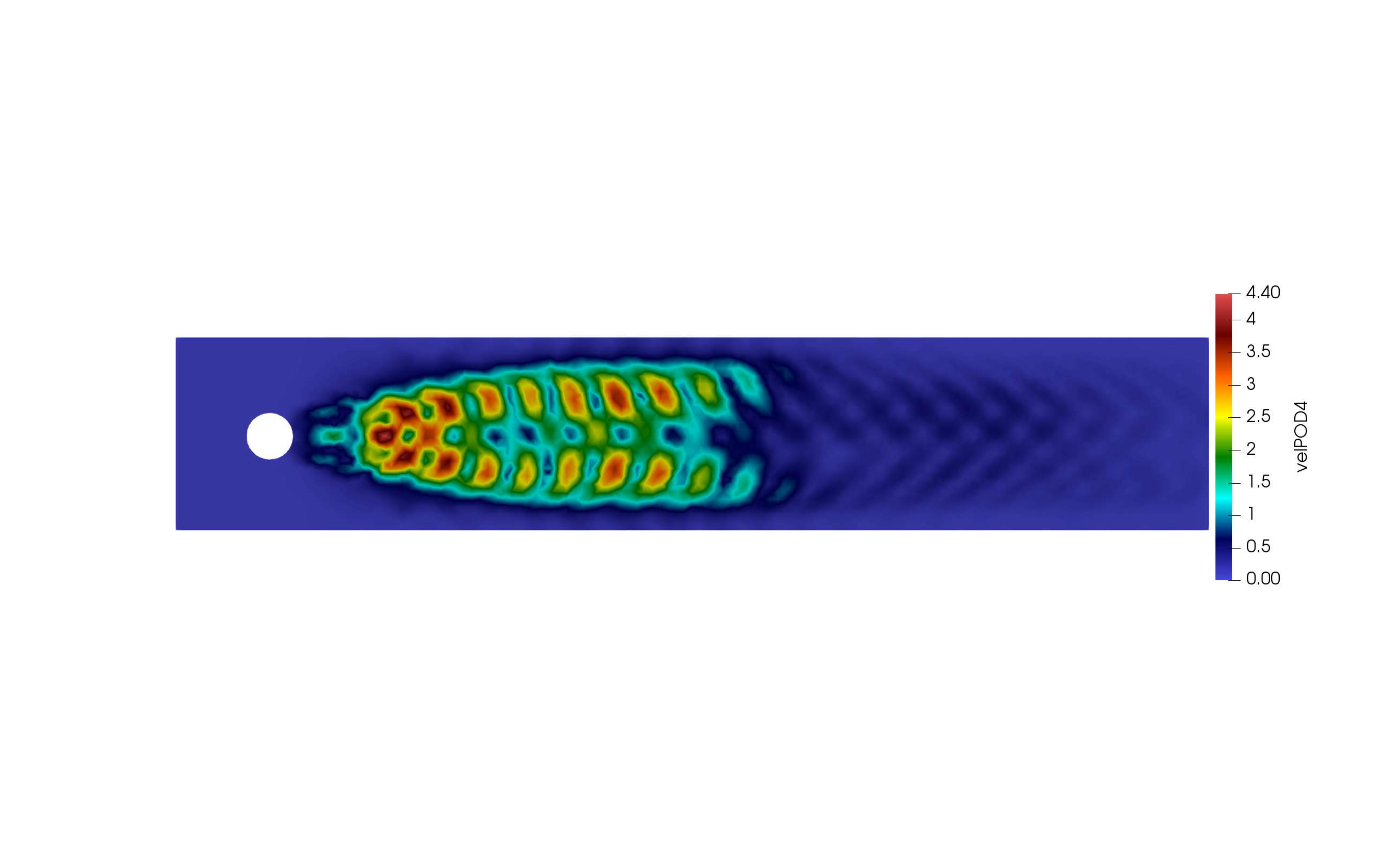}
\caption{First POD velocity modes (Euclidean norm): Dirichlet BC at the outlet.}\label{fig:velPODmodesDir}
\end{center}
\end{figure}

\begin{figure}[htb]
\begin{center}
\includegraphics[width=3.6in]{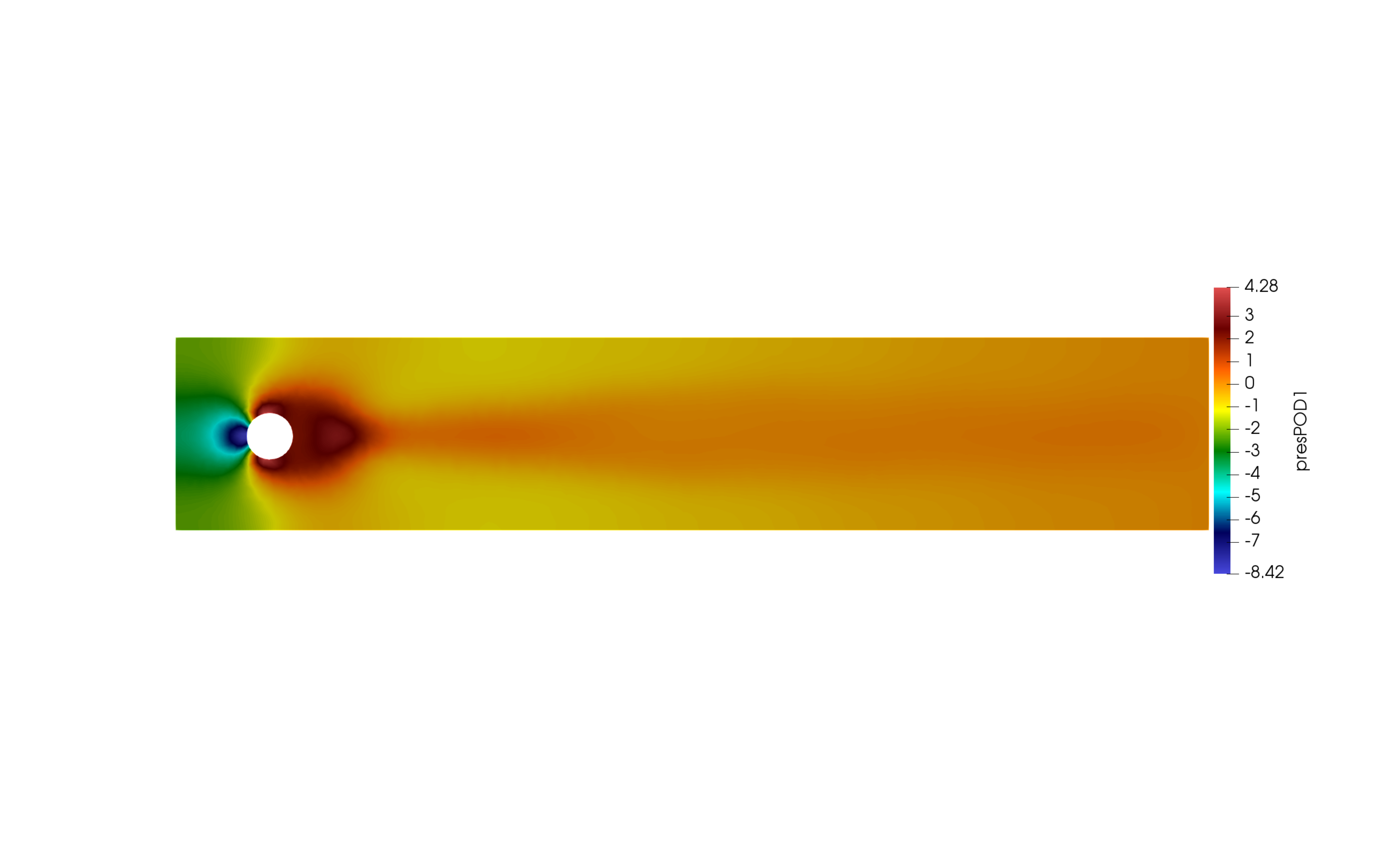}
\includegraphics[width=3.6in]{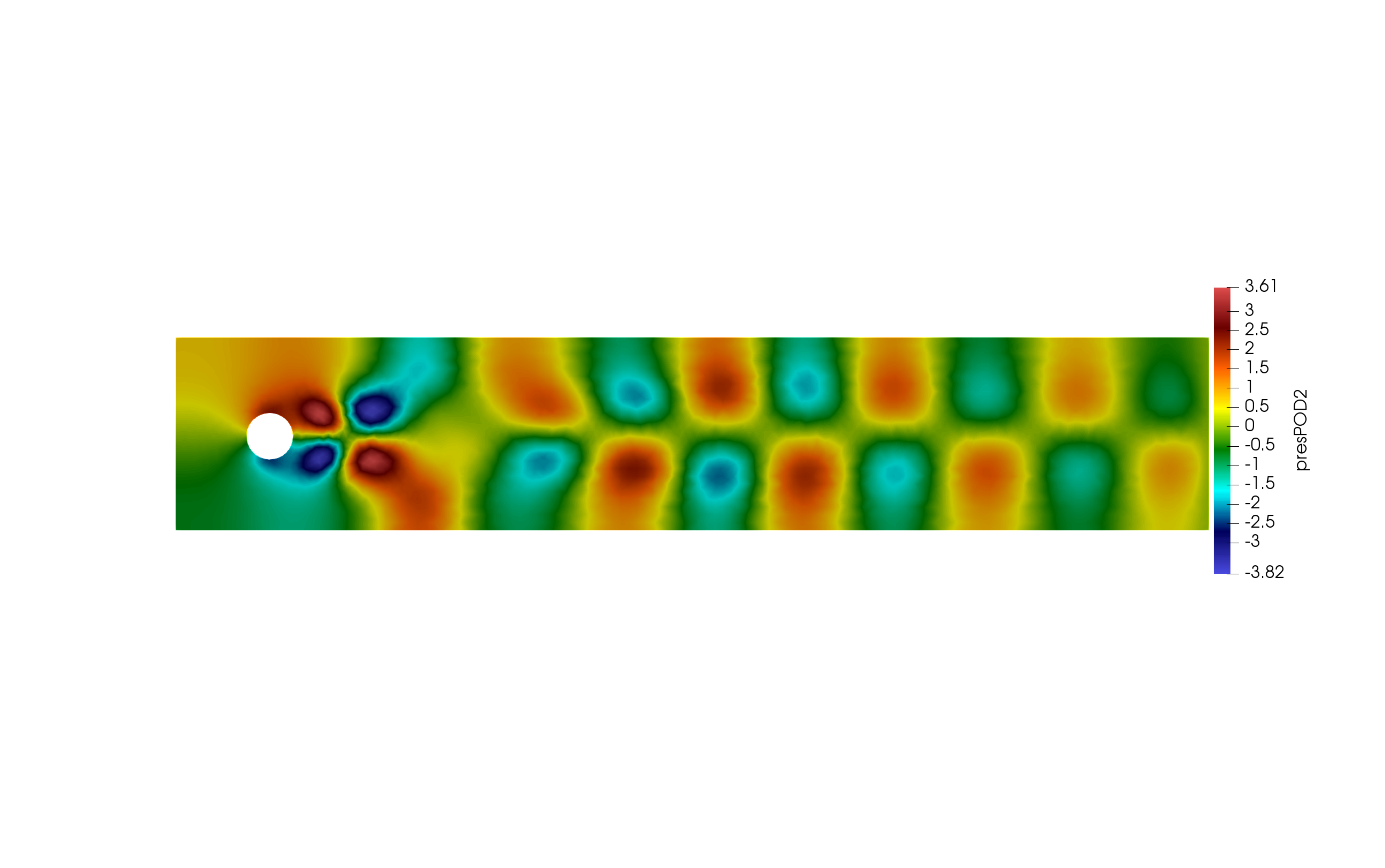}
\includegraphics[width=3.6in]{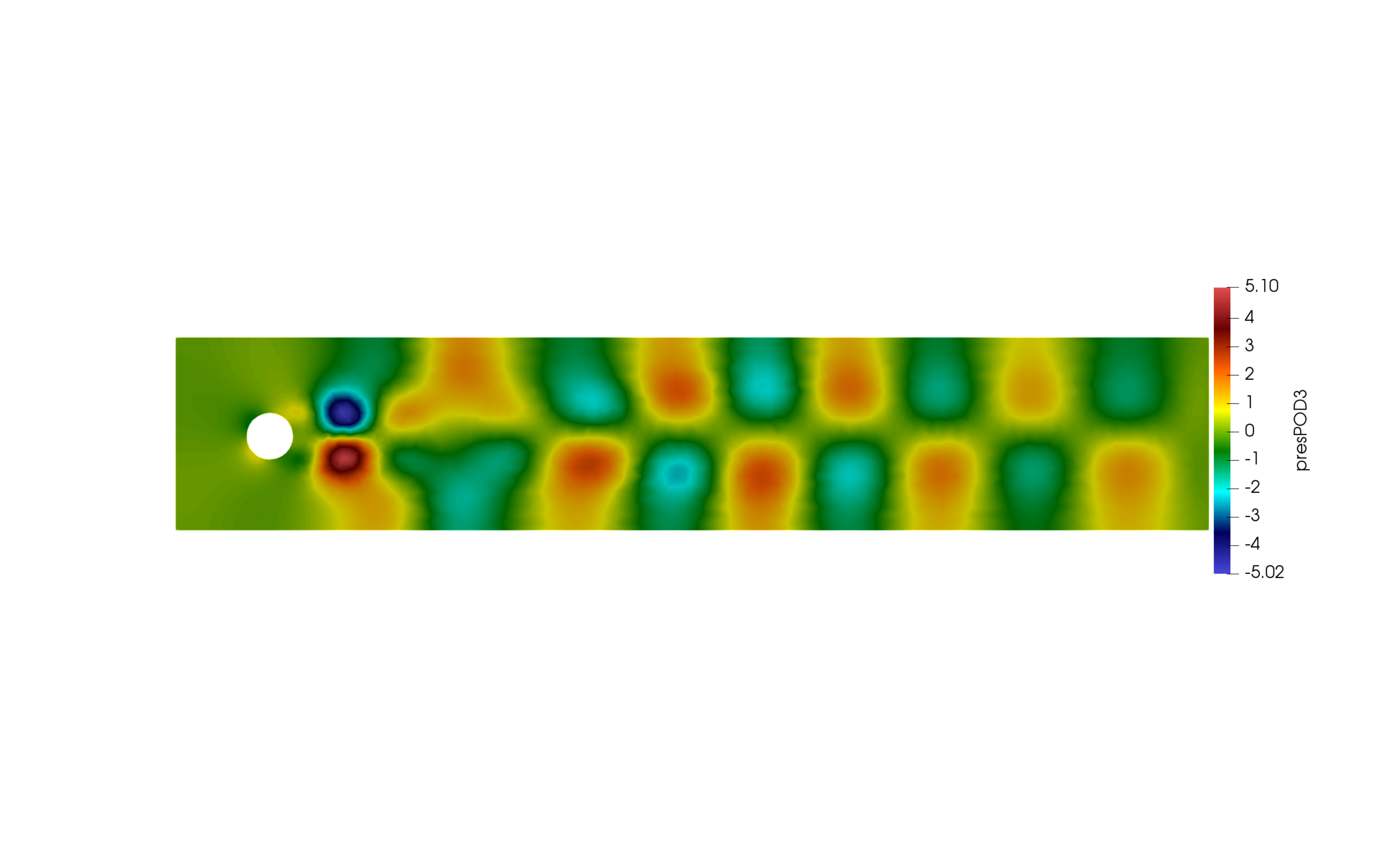}
\includegraphics[width=3.6in]{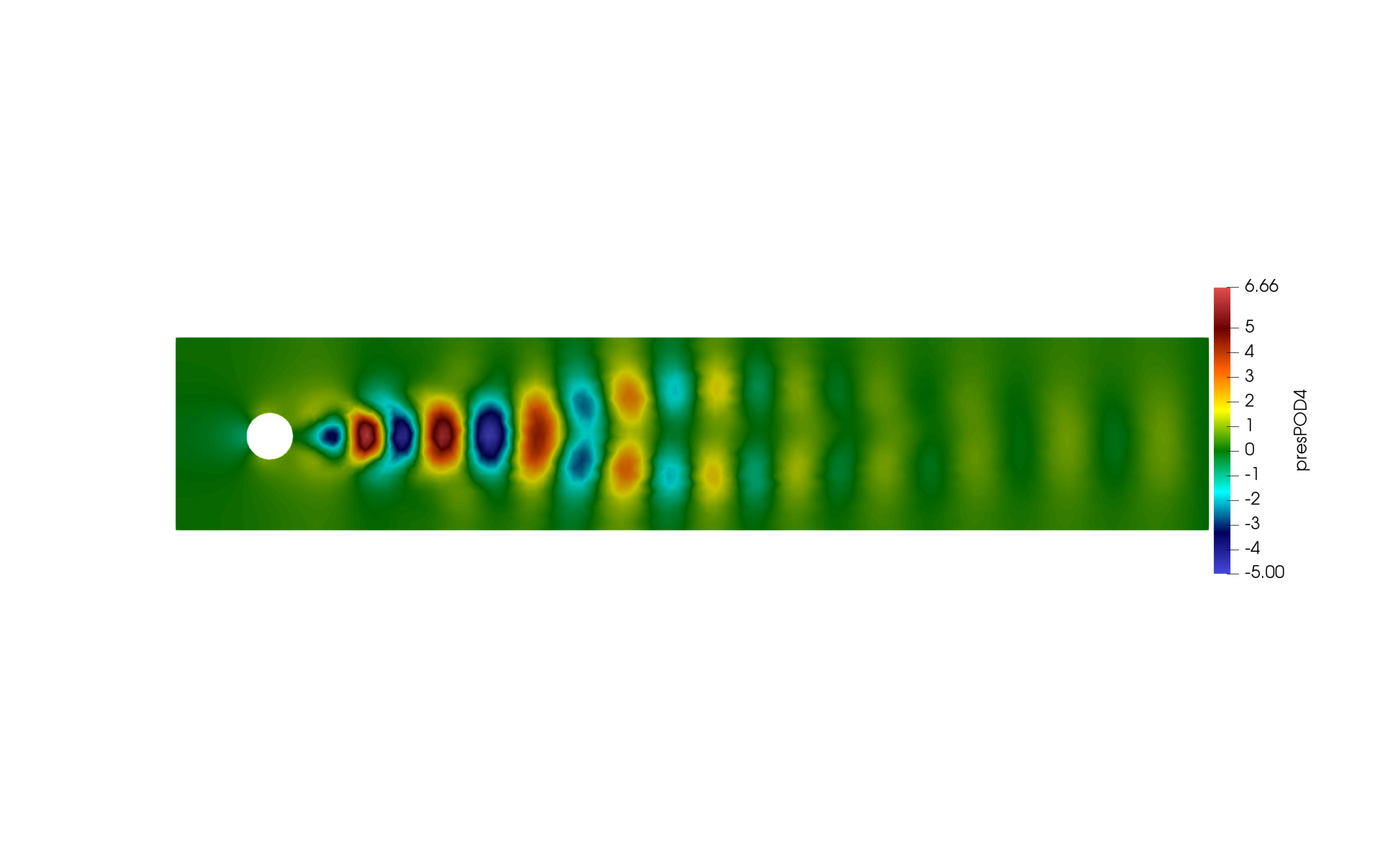}
\caption{First POD pressure modes: do nothing BC at the outlet.}\label{fig:presPODmodes}
\end{center}
\end{figure}

\begin{figure}[htb]
\begin{center}
\includegraphics[width=3.6in]{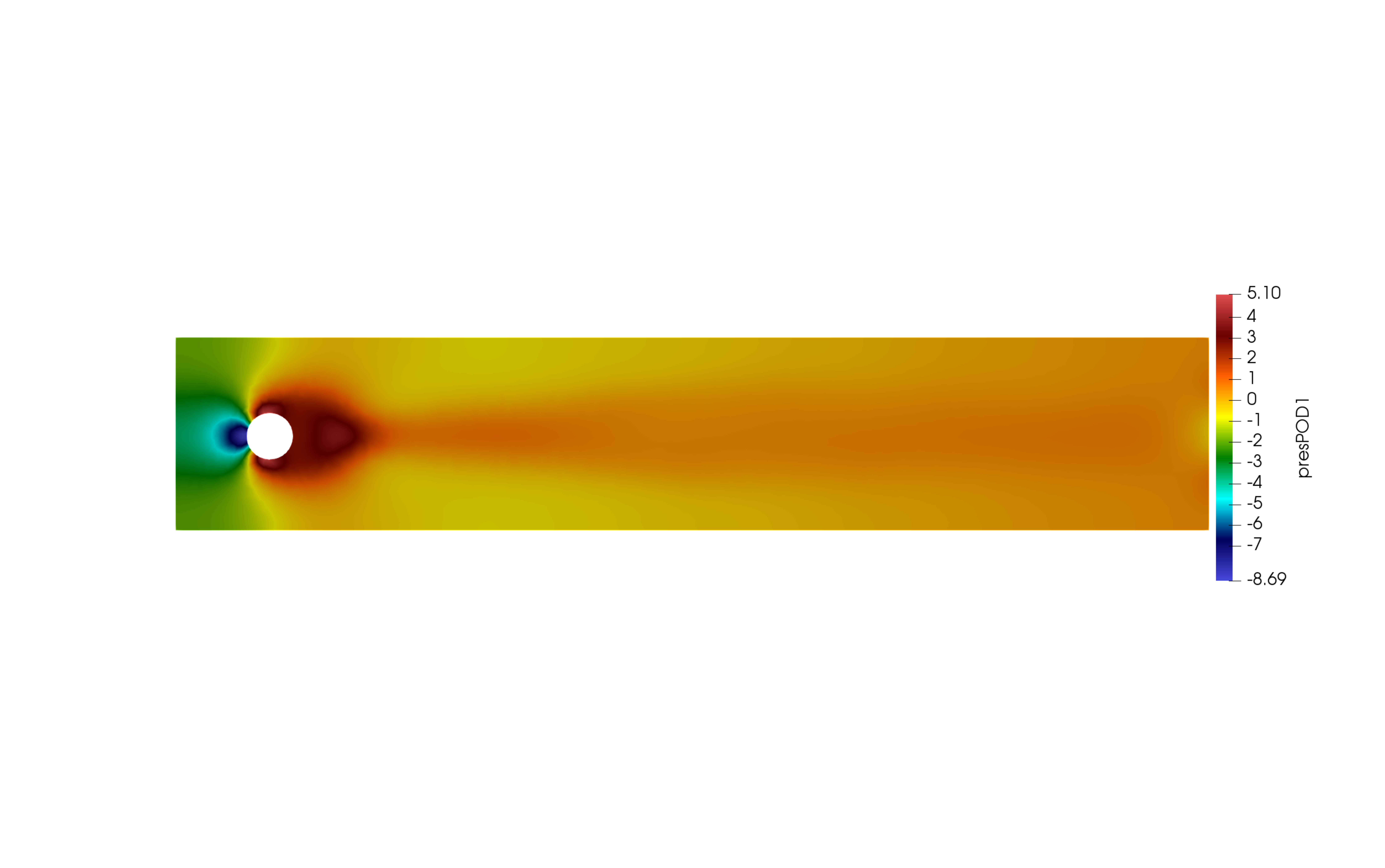}
\includegraphics[width=3.6in]{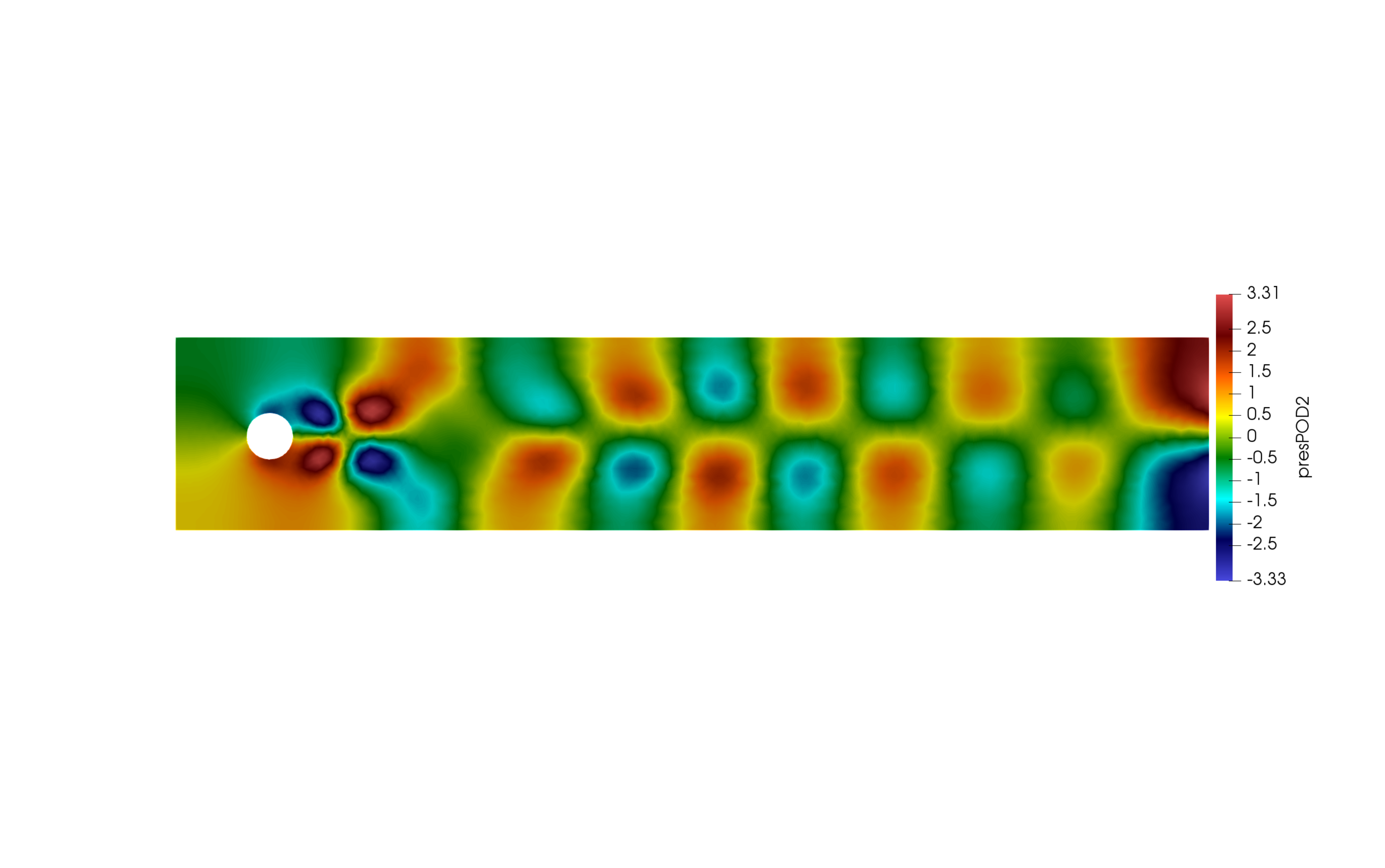}
\includegraphics[width=3.6in]{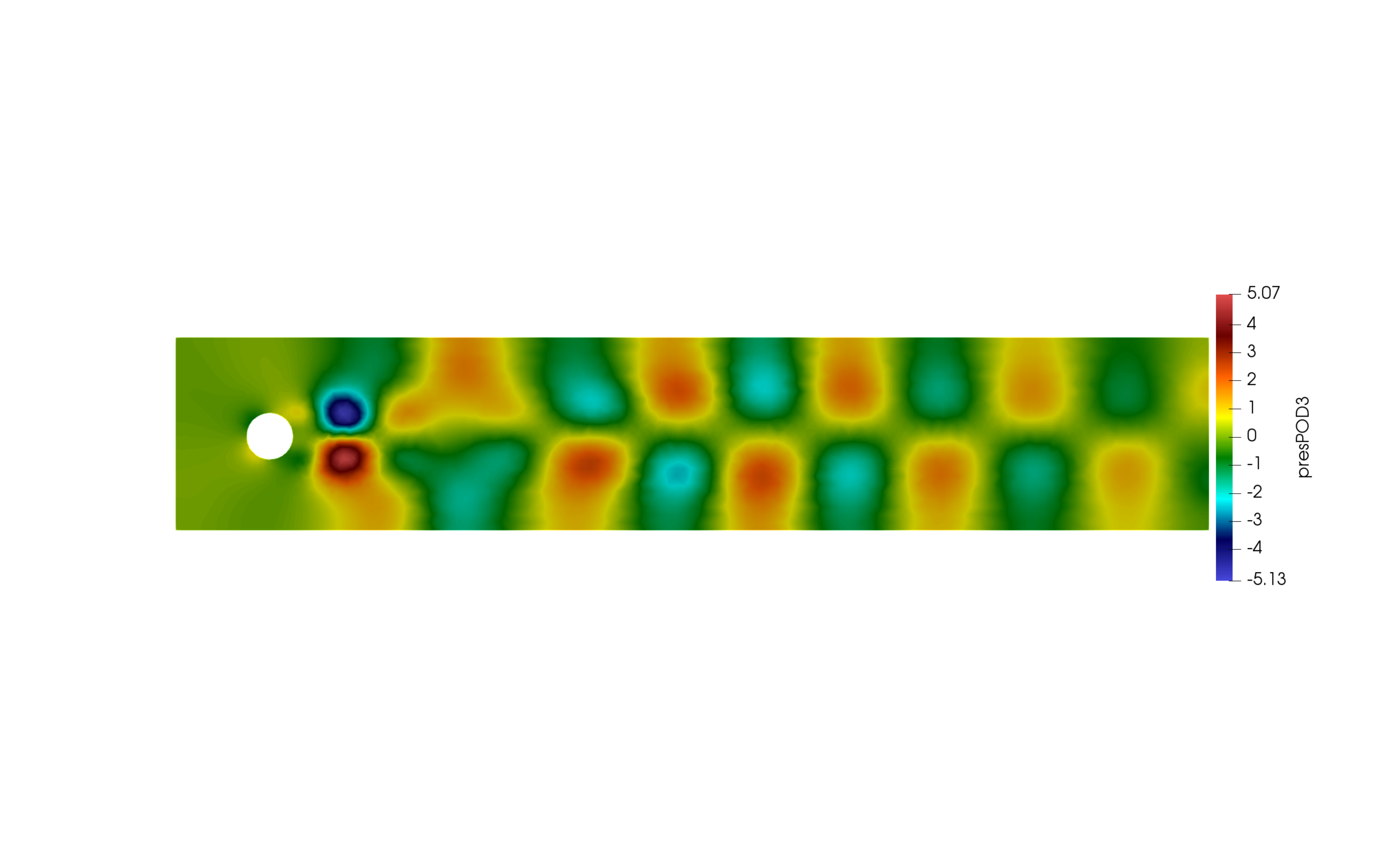}
\includegraphics[width=3.6in]{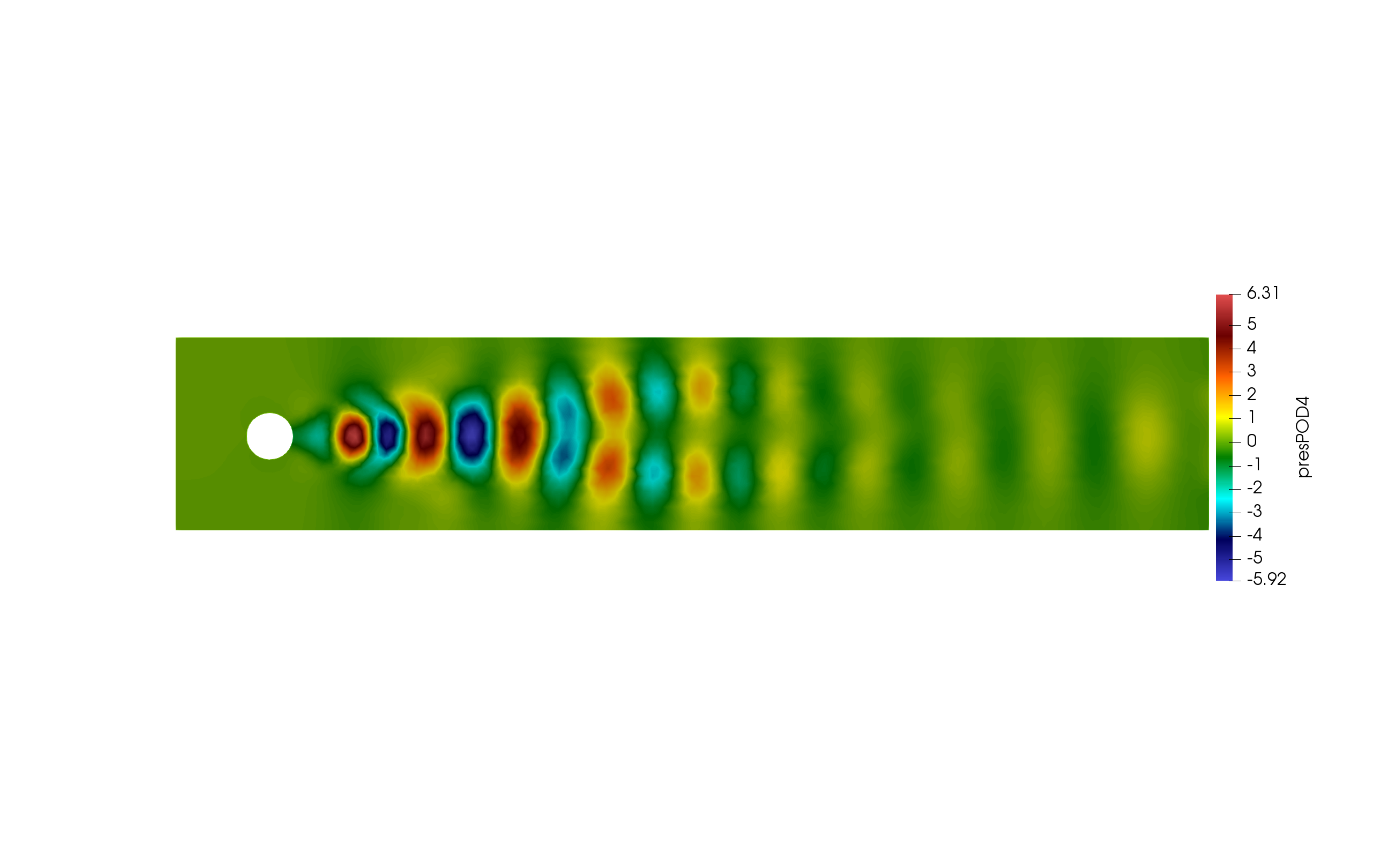}
\caption{First POD pressure modes: Dirichlet BC at the outlet.}\label{fig:presPODmodesDir}
\end{center}
\end{figure}

\begin{figure}[htb]
\begin{center}
\includegraphics[width=6.5in]{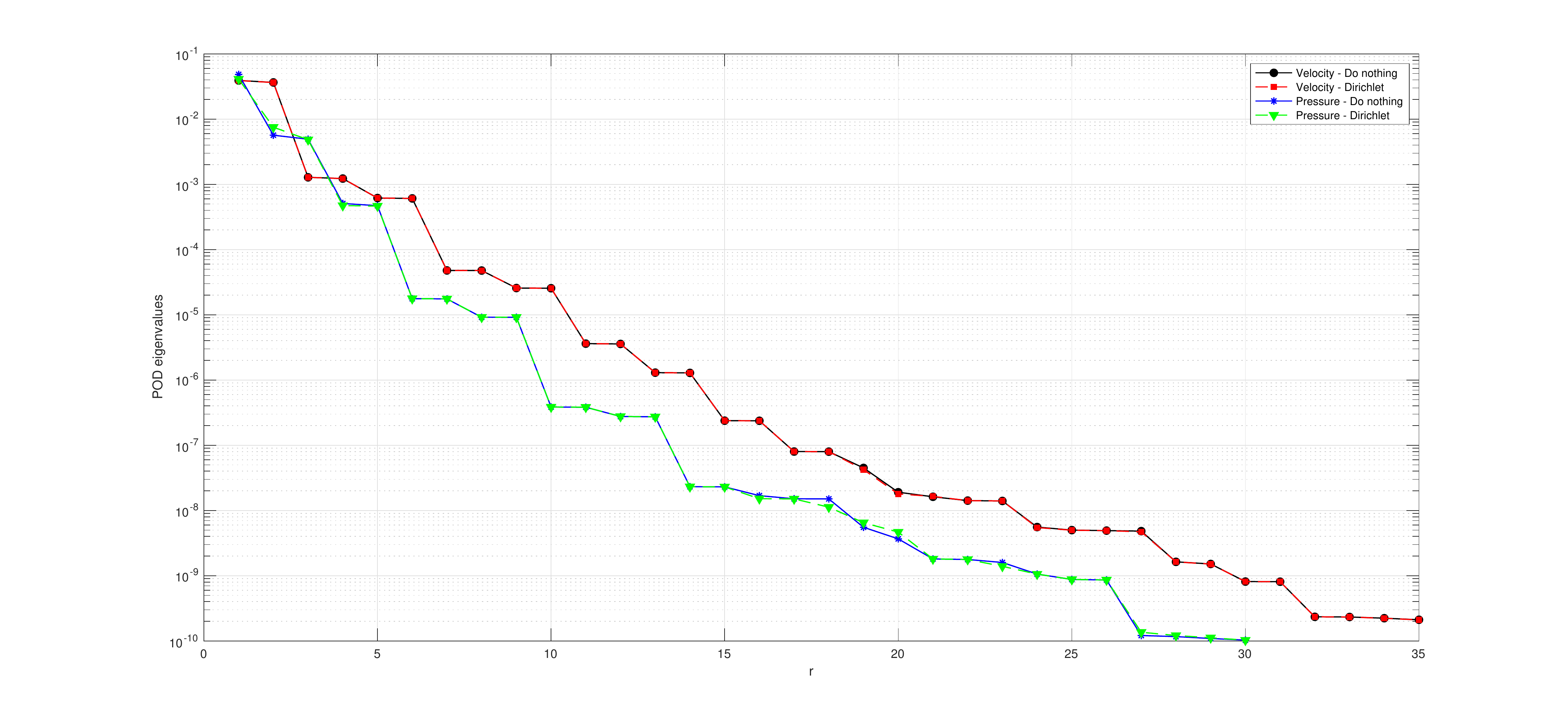}
\includegraphics[width=6.5in]{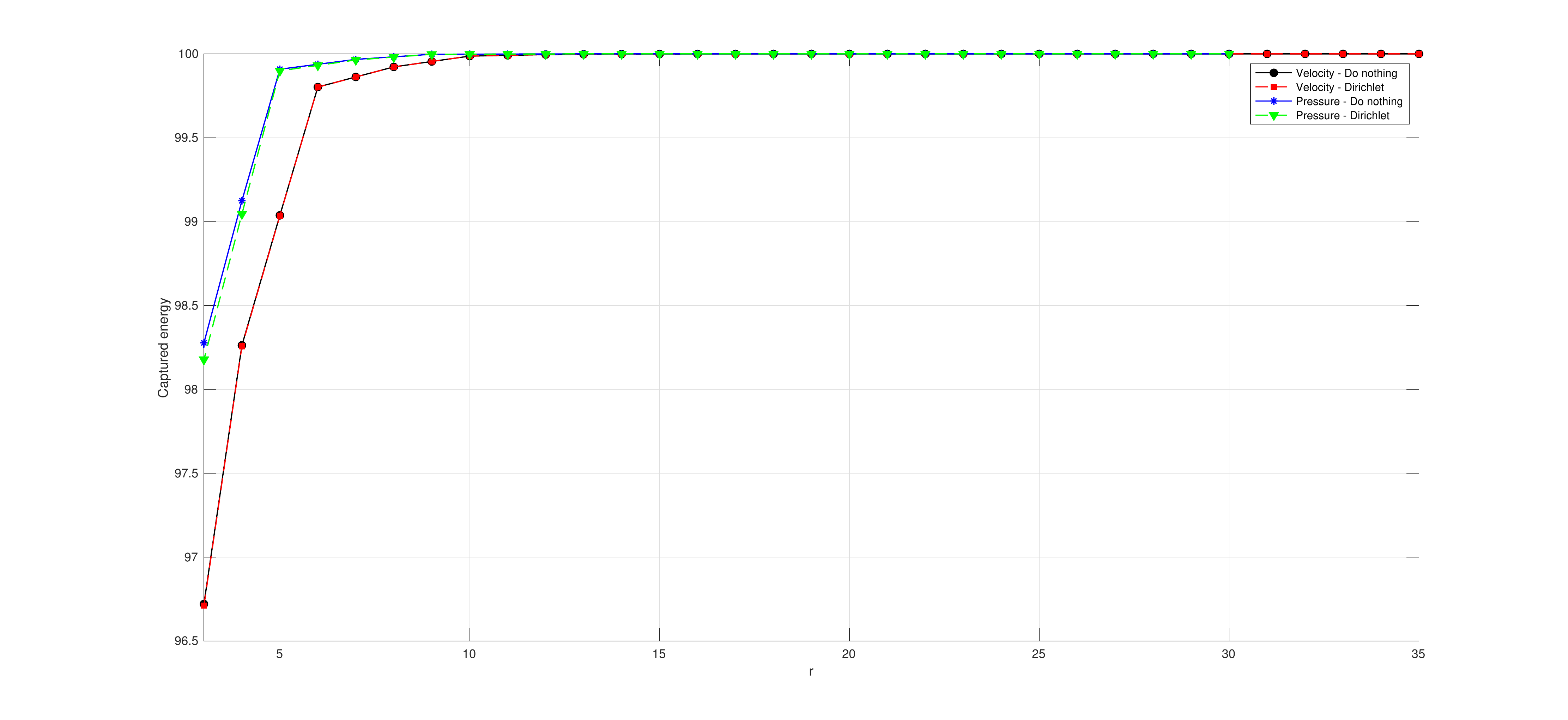}
\caption{POD velocity-pressure eigenvalues (top) and captured system's velocity-pressure energy (bottom) with do nothing BC and Dirichlet BC at the outlet.}\label{fig:EV}
\end{center}
\end{figure}

The POD modes are generated in $L^2$ by the method of snapshots with velocity centered-trajectories \cite{IliescuJohn15} by storing every fifth FOM solution (with do nothing and Dirichlet BC at the outlet) in the stable response time interval $[5,7]\,\rm{s}$, so that $200$ snapshots were used both for velocity and pressure. Figures \ref{fig:velPODmodes}-\ref{fig:velPODmodesDir} display the Euclidean norm of the first POD velocity modes, obtained using do nothing BC and Dirichlet BC at the outlet, respectively. Similarly, Figures \ref{fig:presPODmodes}-\ref{fig:presPODmodesDir} display the first POD pressure modes, obtained using do nothing BC and Dirichlet BC at the outlet, respectively.
Only slight differences can be noticed between the POD modes using do nothing BC and Dirichlet BC at the outlet, being the most noticeable ones given by the numerical noise of the second POD mode using Dirichlet BC at the outlet. In Figure \ref{fig:EV}, we show the decay of POD velocity ($\lambda_i$) and pressure ($\gamma_i$) eigenvalues (top), together with the corresponding captured system's energy (bottom), computed respectively as $100\sum_{i=1}^{r}\lambda_{i}/\sum_{i=1}^{M_v}\lambda_{i}$ and $100\sum_{i=1}^{r}\gamma_{i}/\sum_{i=1}^{M_p}\gamma_{i}$.
Note that the first $r=5$ POD modes already capture more than $99\%$ of the system's velocity-pressure energy. 

\medskip

We recall that the new LPS-ROM uses the same number $r$ of velocity and pressure POD modes. Thus, we expect that also the POD velocity-pressure spaces do not satisfy the standard discrete inf-sup condition and the LPS-ROM try to circumvent it trough nume\-rical stabilization. Effectively, we have checked numerically that, varying $r$, the discrete inf-sup constant $\beta$ for the POD velocity-pressure spaces remains very close to zero (below $10^{-7}$), as we can observe from Figure \ref{fig:IS} (top). In Figure \ref{fig:IS} (bottom), we also display the saturation constant $\alpha$ (see Lemma \ref{lm:AngleCond}) between the spaces $Y=\text{span}\{\div\boldsymbol{\varphi}_i,\ldots,\div\boldsymbol{\varphi}_r\}$ and $Z=Q_r$ in terms of $r$:
\begin{figure}[htb]
\begin{center}
\includegraphics[width=6.5in]{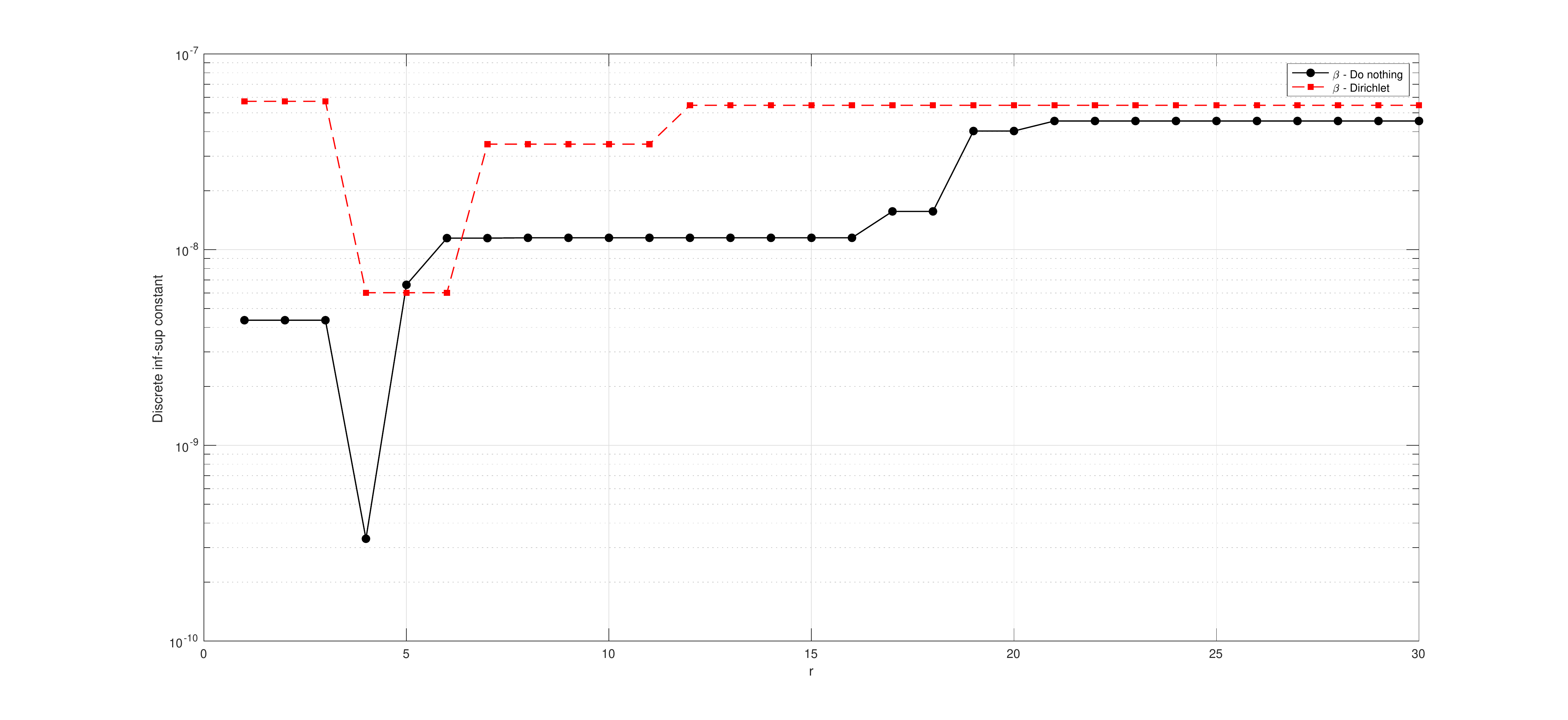}
\includegraphics[width=6.5in]{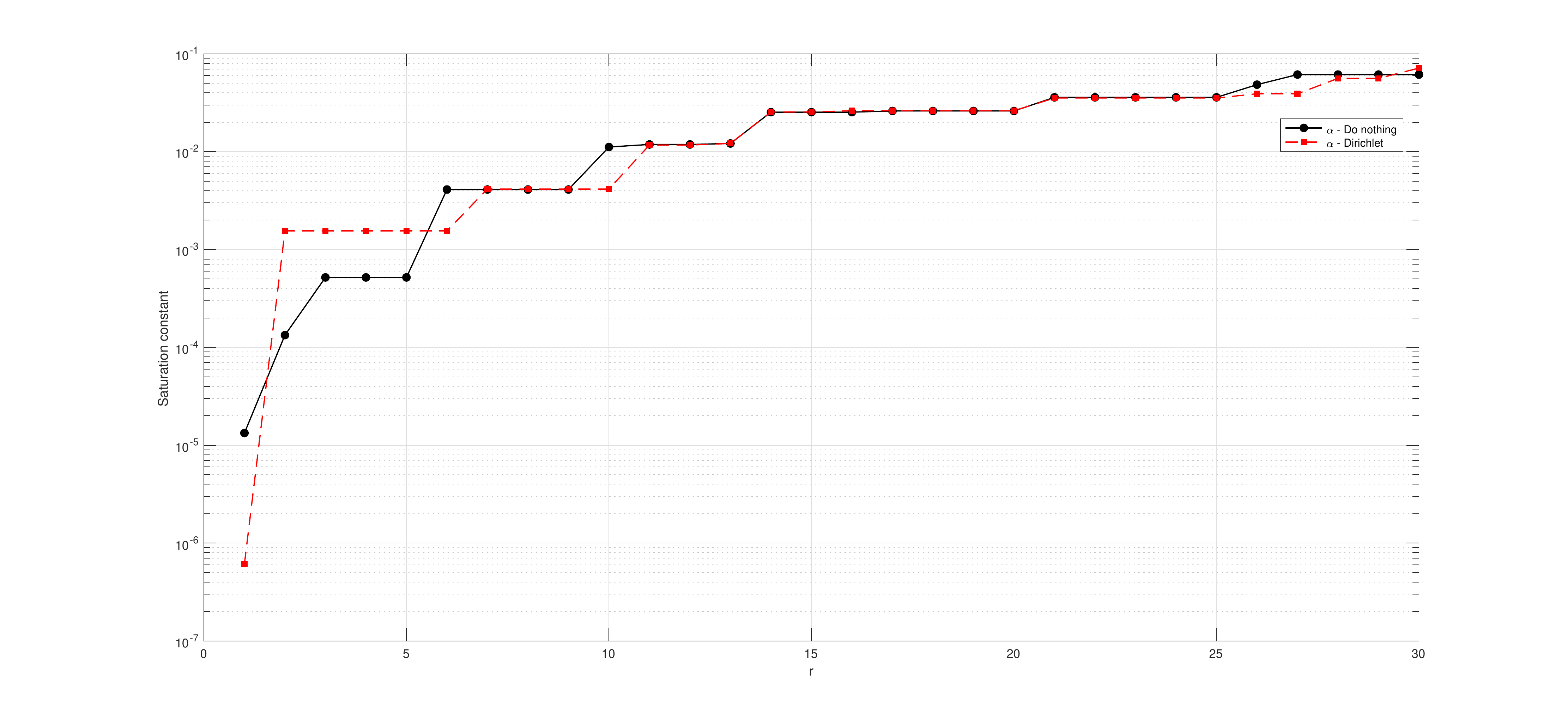}
\caption{Discrete inf-sup constant (top) and saturation constant (bottom) for POD velocity-pressure spaces.}\label{fig:IS}
\end{center}
\end{figure}
We can observe that the saturation constant $\alpha$, for the chosen numerical setup, starts with small values for small $r$ (around $10^{-5}$, $10^{-6}$) and experiences a flattening effect with values within $(10^{-2},10^{-1})$ when adding more POD modes ($r>10$). 
This could explain why we do not experience in practice so much improvement when adding more POD modes and we compare the ROM and the FOM solutions. In terms of the theoretical analysis performed, we have that $\alpha^2<\sigma$ for small $r$, which allows to significantly ease the convergence order reduction due to $\sigma^{-1}$, but then $\alpha^2$ becomes predominant and experiences a flattening effect, which implies no much improvement in the convergence order when adding already more than $7$ POD modes. This could give a sort of criterion (or at least an idea) on how many POD modes are needed to reach, using the new LPS-ROM, a reasonable accuracy with respect to the FOM solution, and thus well catch physical quantities of interest at a very reduced computational cost.

\subsection{Numerical results for LPS-ROM}

With POD velocity-pressure modes generated, the fully discrete LPS-ROM is constructed as discussed in Section \ref{sec:POD-ROM} using the semi-implicit BDF2 time scheme as for the LPS-FOM, and run in the stable response time interval $[5,7]\,\rm{s}$ with $\Delta t=2\times 10^{-3}\,\rm{s}$.

\medskip

\begin{figure}[htb]
\begin{center}
\includegraphics[width=6.5in]{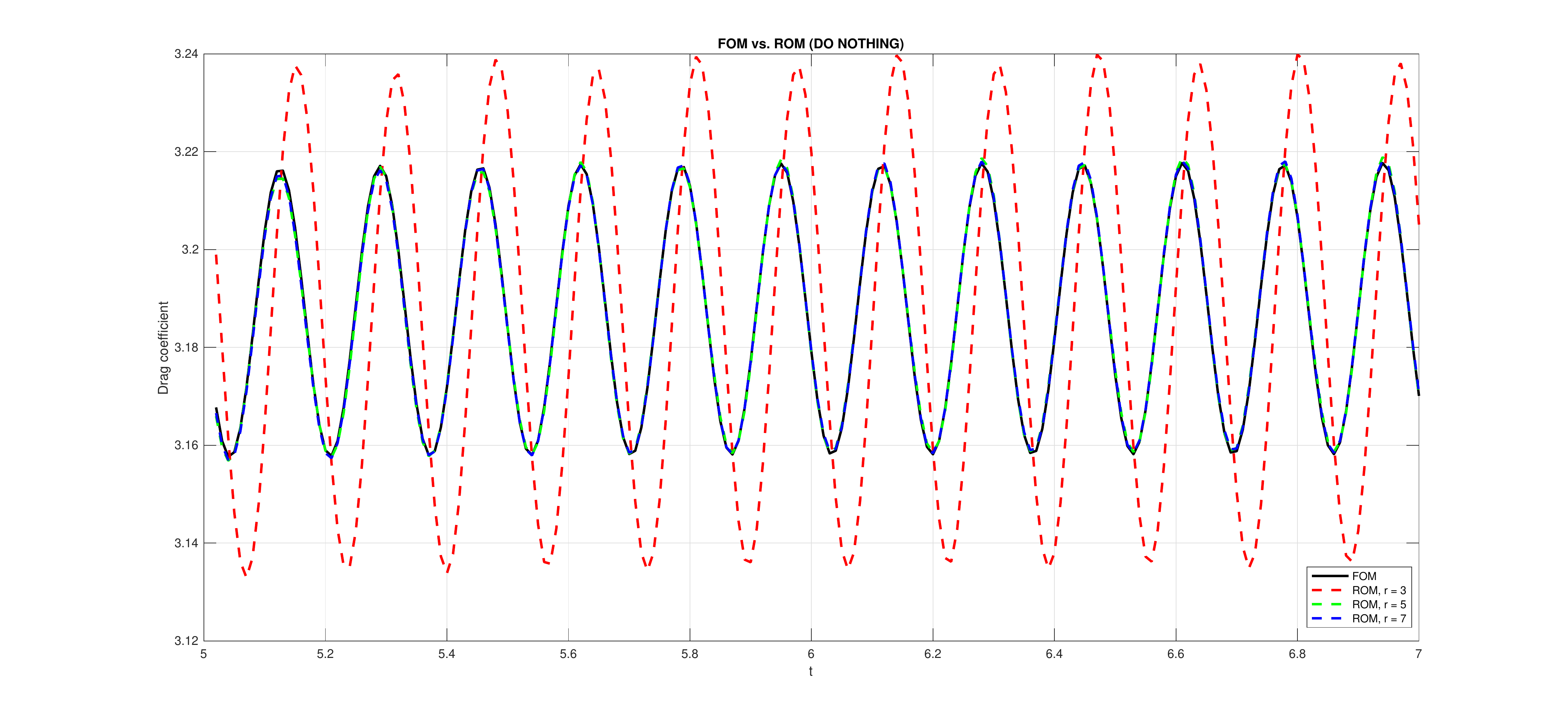}
\includegraphics[width=6.5in]{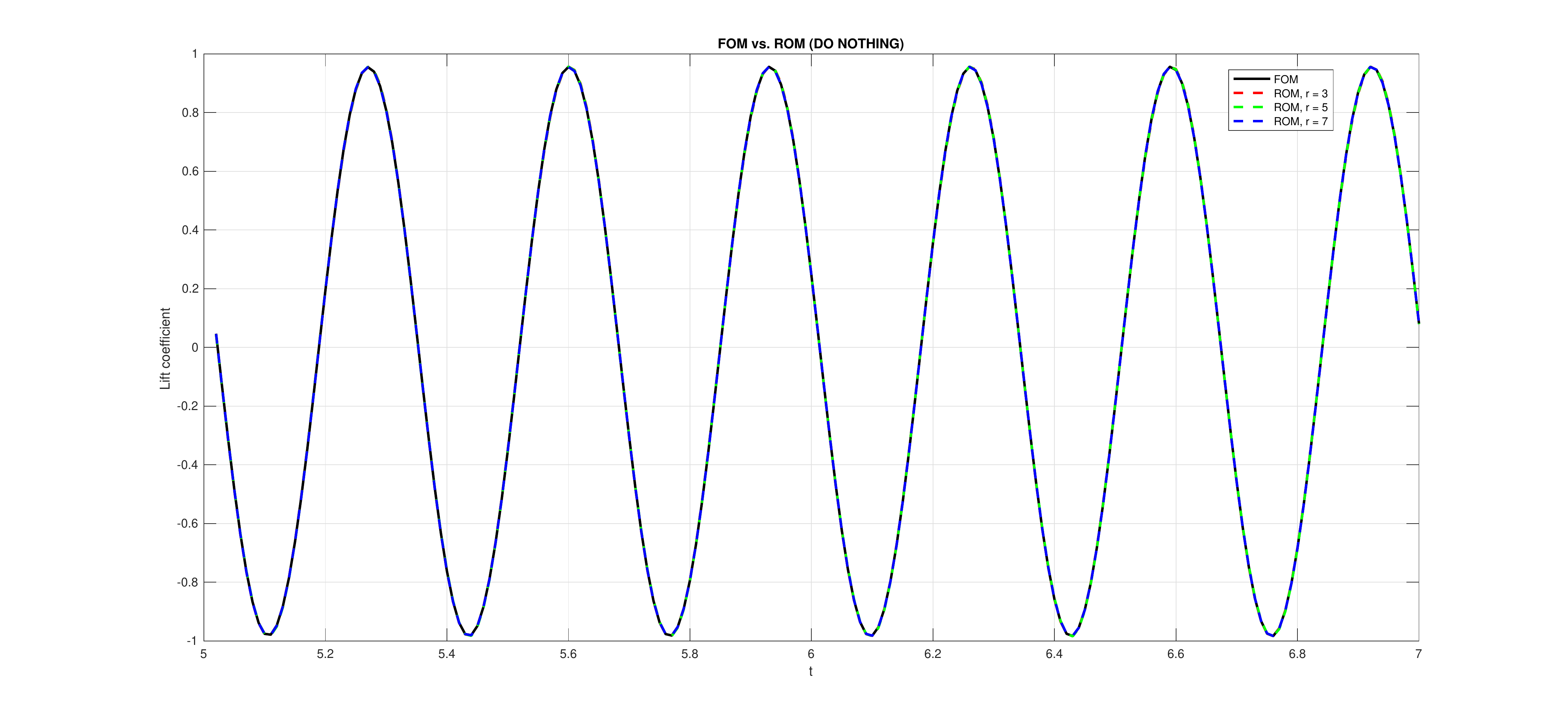}
\includegraphics[width=6.5in]{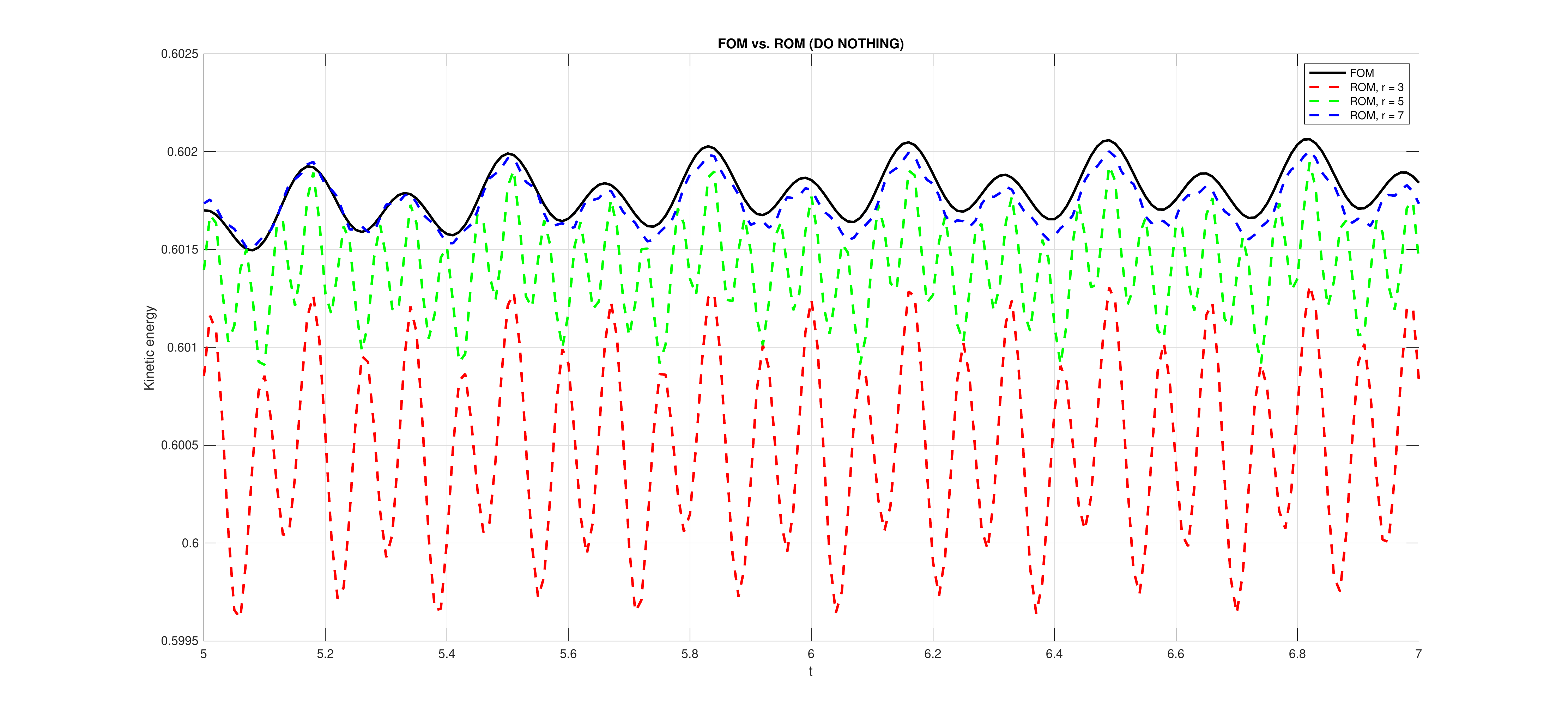}
\caption{Temporal evolution of quantities of interest computed with LPS-ROM ($r=3,5,7$) and compared with LPS-FOM with do nothing BC at the outlet.}\label{fig:DragLiftComp}
\end{center}
\end{figure}

\begin{figure}[htb]
\begin{center}
\includegraphics[width=6.5in]{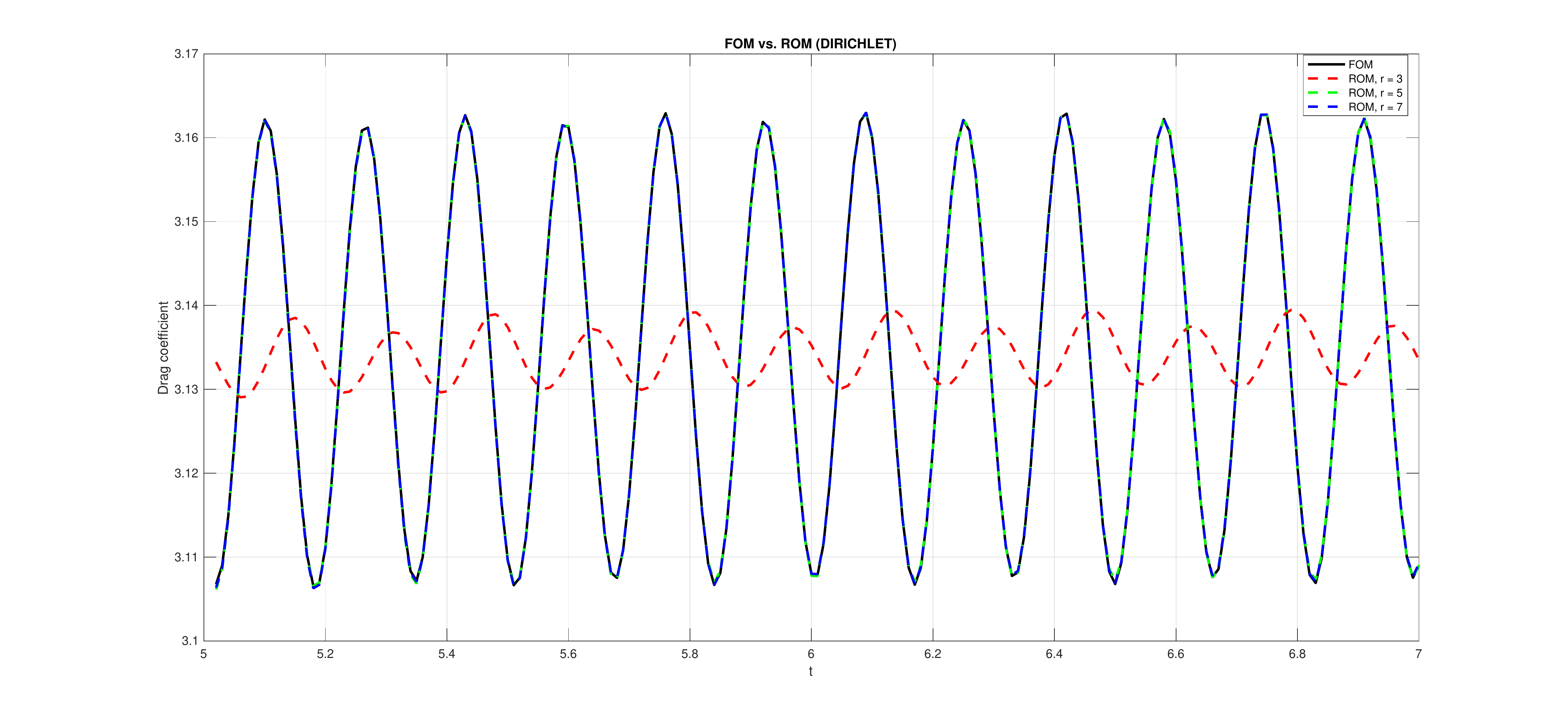}
\includegraphics[width=6.5in]{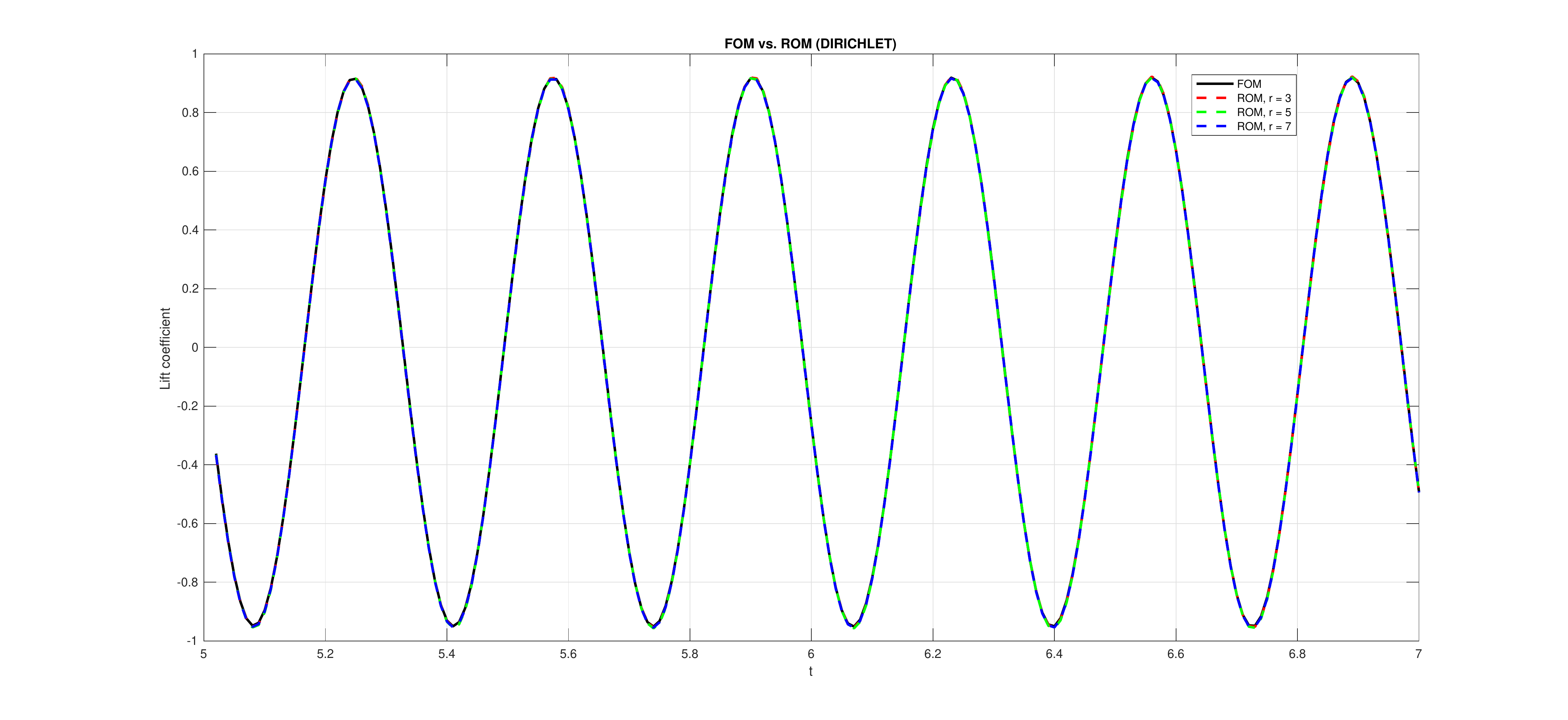}
\includegraphics[width=6.5in]{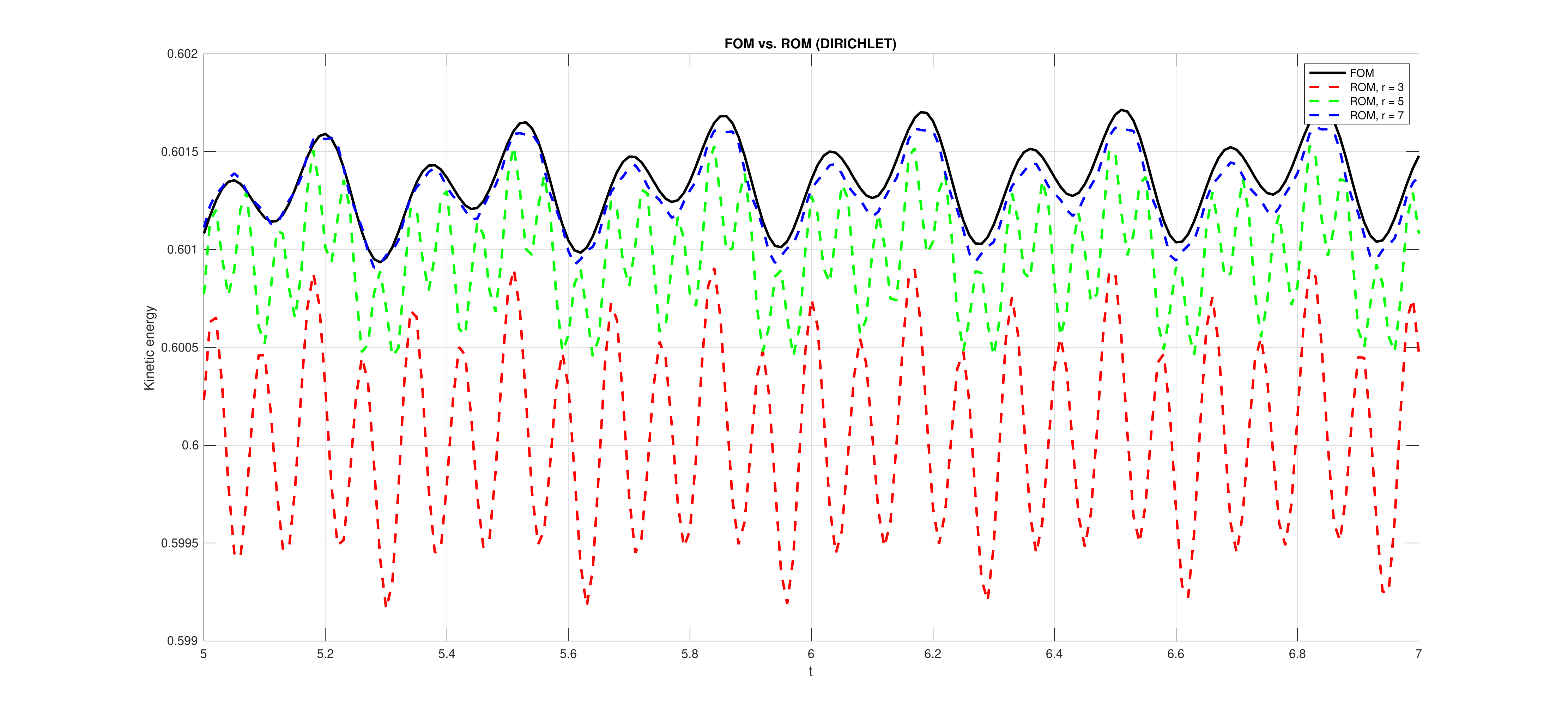}
\caption{Temporal evolution of quantities of interest computed with LPS-ROM ($r=3,5,7$) and compared with LPS-FOM with Dirichlet BC at the outlet.}\label{fig:DragLiftCompDir}
\end{center}
\end{figure}

To assess the behavior of the new LPS-ROM, the temporal evolution of the drag and lift coefficients, and kinetic energy are monitored and compared to the FOM solutions in the stable response time interval $[5,7]\,\rm{s}$, corresponding to six periods for the lift coefficient. The corresponding numerical results with do nothing and Dirichlet BC at the outlet are reported in Figures \ref{fig:DragLiftComp}-\ref{fig:DragLiftCompDir}. The LPS-ROM has been computed with $r=3,5,7$ POD modes for velocity and pressure.
We observe that the new LPS-ROM is able to replicate the temporal evolution of the lift coefficient computed with the LPS-FOM already with $r=3$. As for the drag coefficient, $r=5$ POD modes are already sufficient to capture its corresponding FOM values. More sensitive seems to be the kinetic energy, for which however $r=7$ POD modes are already able to catch up the FOM solution with a reasonable accuracy. Apart from the temporal evolution of the drag coefficient with $r=3$, results with do nothing BC and Dirichlet BC at the outlet are almost similar for the new LPS-ROM. Qualitatively, it is interesting to observe that, as the LPS-FOM, the new LPS-ROM shows a fully periodic time evolution for all monitored quantities of interest at all levels ($r=3,5,7$).

\medskip

\begin{figure}[htb]
\begin{center}
\includegraphics[width=6.5in]{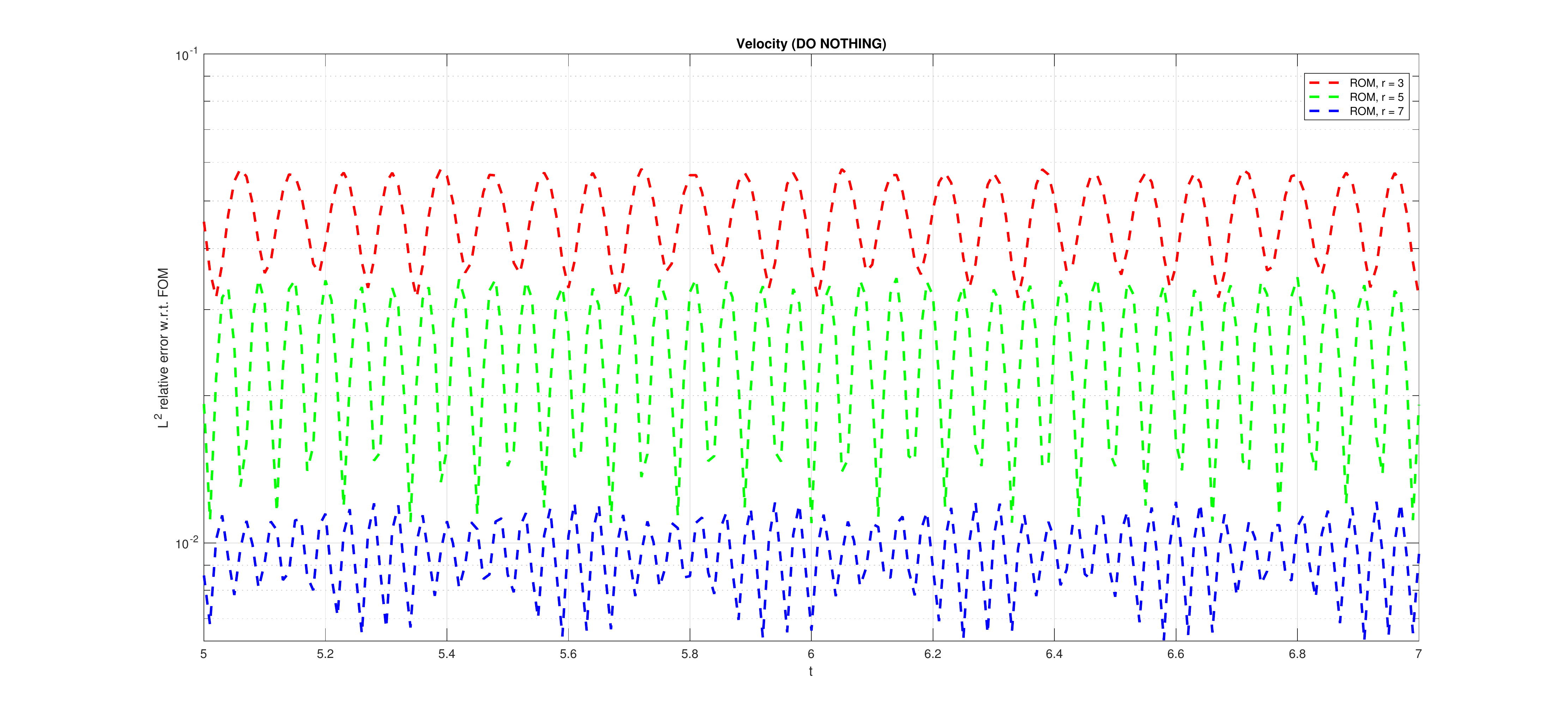}
\includegraphics[width=6.5in]{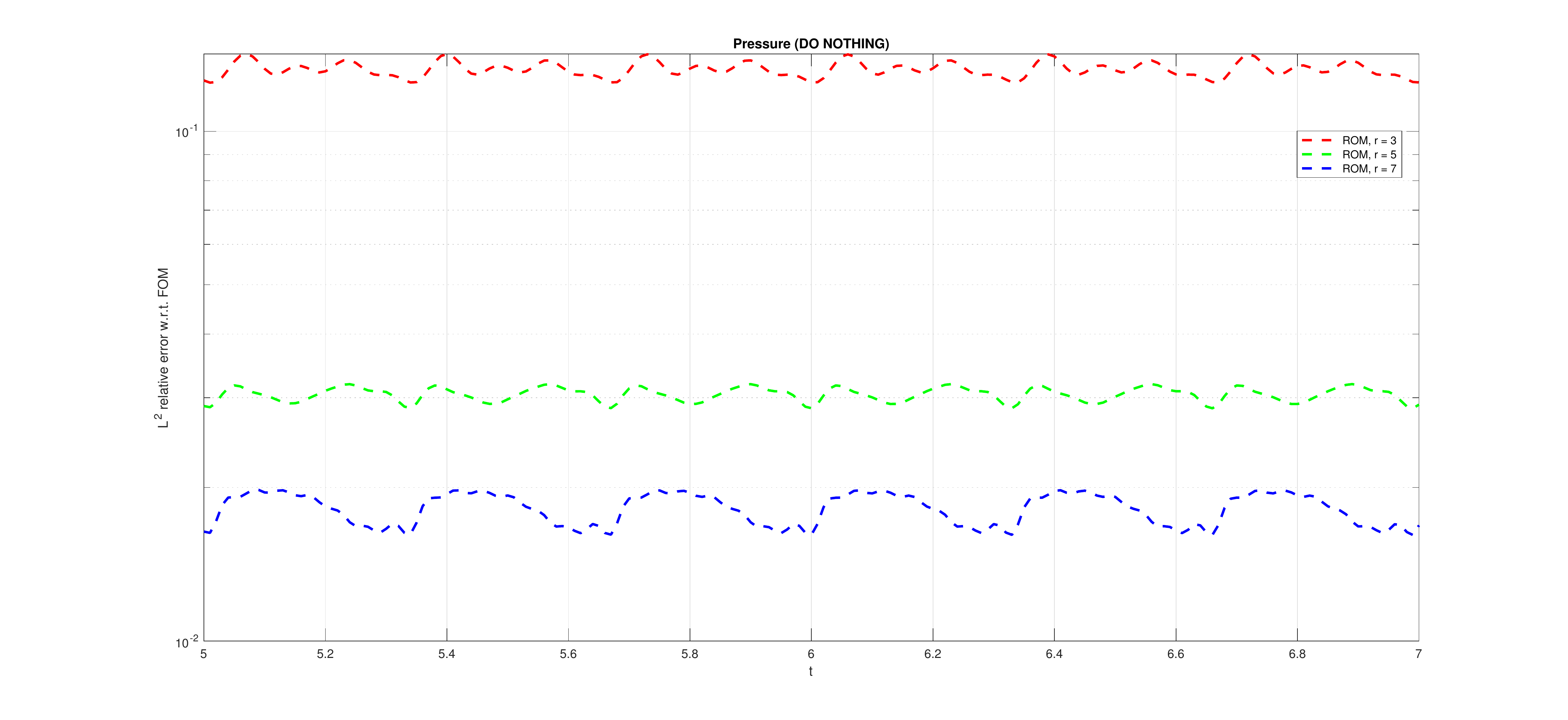}
\caption{Temporal evolution of discrete $L^2$ relative error of LPS-ROM velocity and pressure ($r=3,5,7$) with respect to LPS-FOM ones with do nothing BC at the outlet.}\label{fig:L2RelErr}
\end{center}
\end{figure}

\begin{figure}[htb]
\begin{center}
\includegraphics[width=6.5in]{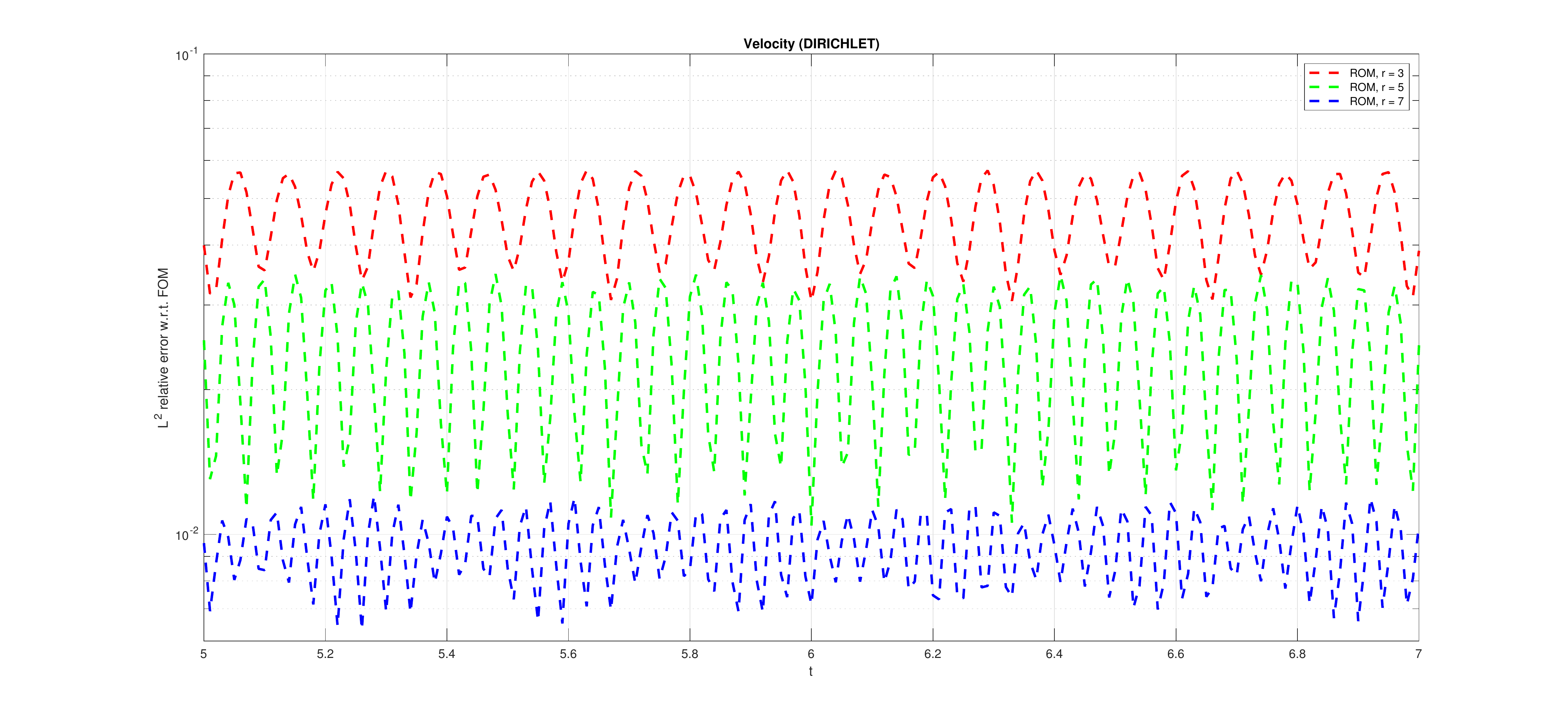}
\includegraphics[width=6.5in]{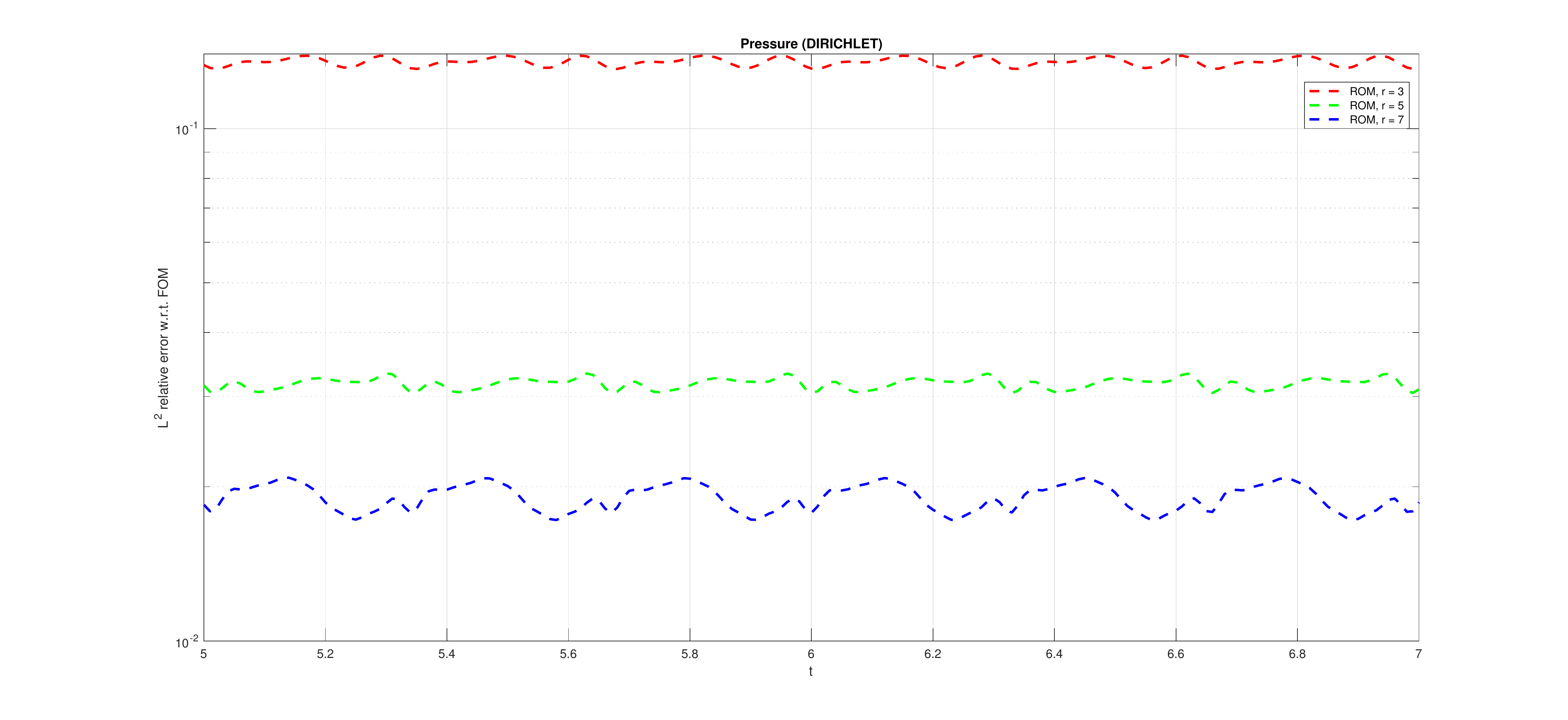}
\caption{Temporal evolution of discrete $L^2$ relative error of LPS-ROM velocity and pressure ($r=3,5,7$) with respect to LPS-FOM ones with Dirichlet BC at the outlet.}\label{fig:L2RelErrDir}
\end{center}
\end{figure}

\begin{figure}[htb]
\begin{center}
\includegraphics[width=6.5in]{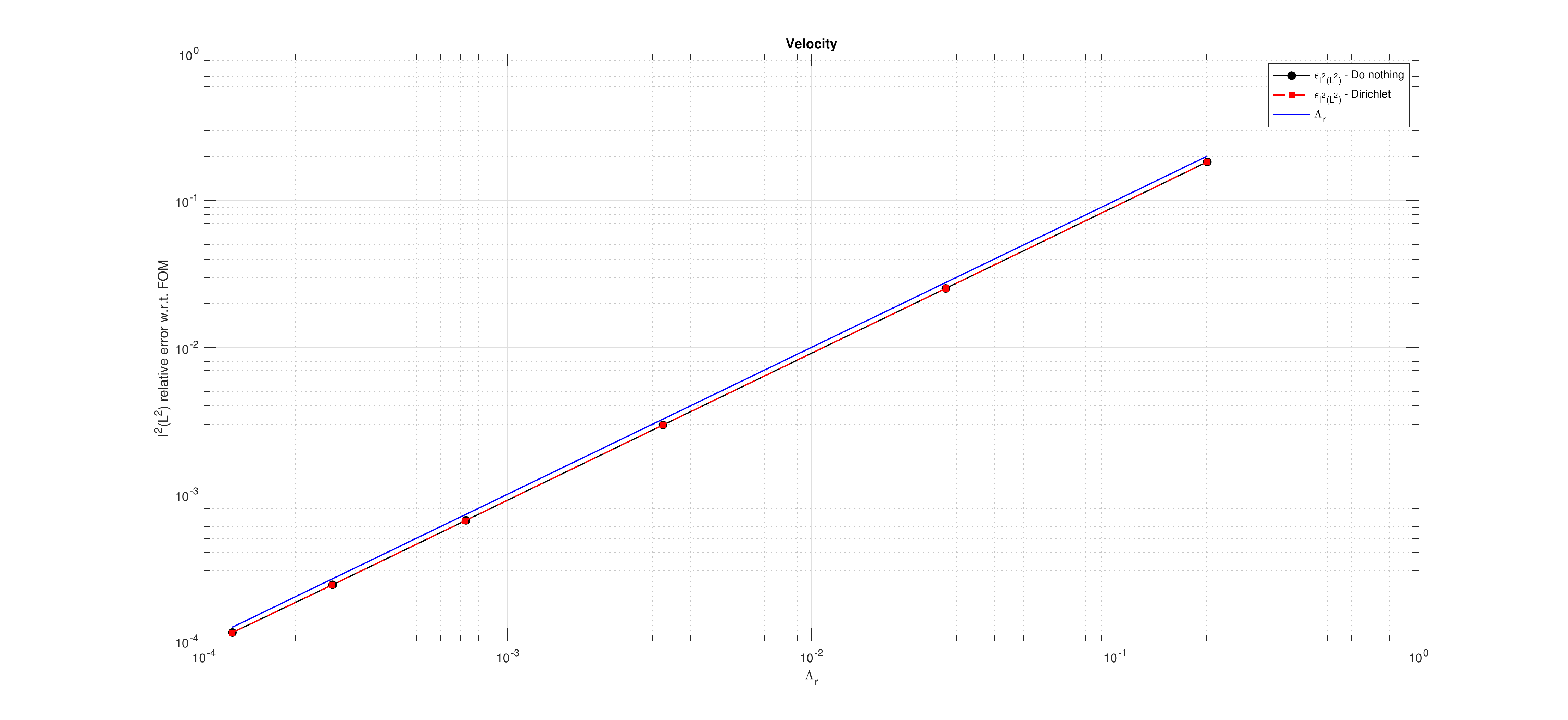}
\includegraphics[width=6.5in]{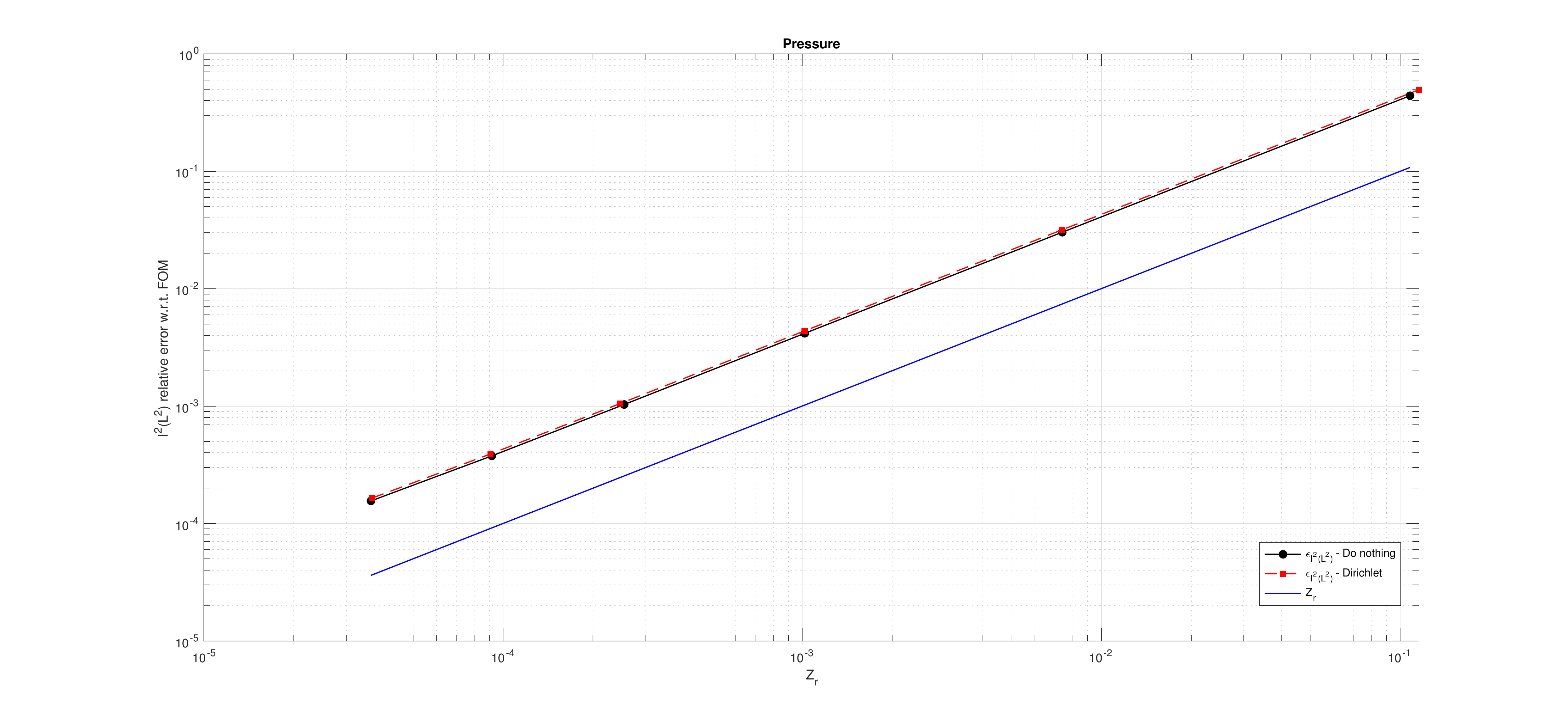}
\caption{Discrete $\ell^2(L^2)$ relative error of LPS-ROM velocity and pressure with respect to $L^2$-projected LPS-FOM ones (on the POD spaces) with do nothing and Dirichlet BC at the outlet in terms of $\Lambda_r = \sum_{i=r+1}^{M_v}\lambda_i$ and $Z_r = \sum_{i=r+1}^{M_p}\gamma_i$, respectively.}\label{fig:l2L2RelErr}
\end{center}
\end{figure}

To better assess on the one hand the behavior of the new LPS-ROM and partially illustrate on the other hand the theoretical convergence order predicted by the numerical analysis and stated in Theorem \ref{th:PODEE}, we plot relative errors with respect to the LPS-FOM solution. In particular, in Figure \ref{fig:L2RelErr}, we first plot the temporal evolution of the discrete $L^2$ relative error (in semilogarithmic scale) of the reduced order velocity and pressure with respect to the full order ones: $\nor{\uhv-\uv_r}{{\bf L}^2}/\nor{\uhv}{{\bf L}^2}$ and $\nor{p_h-p_r}{L^2}/\nor{p_h}{L^2}$, for $r=3,5,7$, obtained with do nothing BC. Similarly, Figure \ref{fig:L2RelErrDir} shows the temporal evolution of the discrete $L^2$ relative error of the reduced order velocity and pressure with respect to the full order ones for $r=3,5,7$, obtained with Dirichlet BC at the outlet.
As expected, the errors decrease when increasing $r$, but then they almost stabilize for $r\geq 7$. It can be observed that the use of do nothing and Dirichlet BC at the outlet provides almost similar results. However, in order to limit the influence of the spatial error due to the FE discretization and the temporal error due to the time-stepping scheme, we have decided to test the $\ell^2(L^2)$ relative errors of the reduced velocity and pressure with respect to the $L^2$ projection of the full order ones on the respective POD spaces, and plot them in terms of $\Lambda_r = \sum_{i=r+1}^{M_v}\lambda_i$ and $Z_r = \sum_{i=r+1}^{M_p}\gamma_i$, respectively. In this way, the analysis is focused on the ROM error due to the POD truncation and one would expect a linear regression behavior (in logarithmic scale), as suggested by the performed numerical analysis. Effectively, this behavior is numerically recovered in Figure \ref{fig:l2L2RelErr} for both reduced velocity and pressure.

\section{Summary and conclusions}\label{sec:Concl}

In this paper, we have proposed a new stabilized projection-based reduced order method (LPS-ROM) for the numerical simulation of incompressible flows.
In particular, the new LPS-ROM is a velocity-pressure ROM that uses pressure modes as well to compute the reduced order
pressure, needed for instance in the computation of relevant quantities, such as drag and
lift forces on bodies in the flow. 
With respect to other approaches existing in the current ROM literature that provides velocity-pressure approximations, the new LPS-ROM circumvents the standard discrete inf-sup condition for the POD velocity-pressure spaces, whose fulfillment can be rather expensive and inefficient in realistic applications in CFD, see for instance \cite{Rozza15,RozzaStabile18}, where an offline strategy based on the supremizer enrichment of the reduced velocity space has been proposed and applied in the POD context, adapted from the RB method framework. Also, the velocity modes for the new LPS-ROM does not have to be
neither strongly nor weakly divergence-free, which allows to use snapshots generated
for instance with penalty or projection-based stabilized methods. This is not the case, for instance, of ROM based on a pressure Poisson equation approach (see, for instance, the first two methods investigated in \cite{IliescuJohn14} and also \cite{RozzaStabile18}), for which the velocity snapshots, and hence the POD velocity modes must be at least weakly divergence-free. This requirement also holds for the last method proposed and investigated in \cite{IliescuJohn14}, which uses a residual-based stabilization mechanism in order to overcome a possible violation of the discrete inf-sup condition in the ROM framework by considering a decoupled approach for the reduced velocity-pressure pair.

\medskip

The main contribution of the present paper has been to perform a stability and convergence analysis of the arising fully discrete LPS-ROM applied to the unsteady incompres\-sible NSE, by mainly deriving the proof of a rigorous error estimate that considers all contributions: the spatial discretization error (due to the FE discretization), the temporal discretization error (due to the backward Euler method), and the POD truncation error. In particular, the numerical analysis corroborated with numerical studies makes apparent an interesting link between the number of POD velocity-pressure modes used in the ROM and the angle $\theta$ between the space spanned by the divergence of the POD velocity modes and the POD pressure space. Indeed, for the numerical setting proposed, where the same numerical stabilization technique used for the ROM is initially applied also to the FOM (LPS-FOM) to generate the snapshots, so that these latter are not weakly divergence-free, we have found that for small values of $r$, which is common in practice, the saturation constant $\alpha=cos(\theta)$ is rather small, and this allows to ease the convergence order reduction due to the violation of the discrete inf-sup stability condition.

\medskip

Numerical studies performed on a two-dimensional laminar unsteady flow past a circular obstacle have also been used to assess the accuracy and efficiency of the new LPS-ROM. Despite the fact that the discrete inf-sup condition is not fulfilled by the new LPS-ROM, using a small equal number of POD velocity-pressure modes already provides accurate approximations, close to the LPS-FOM results, and theoretical scalings suggested by the numerical analysis are recovered in practice.

\medskip

We plan to extend the theoretical work of this paper to the numerical analysis of the last method proposed and investigated in \cite{IliescuJohn14} but for not exactly weakly divergence-free POD velocity modes, by also performing a numerical investigation that both supports the analytical results and illustrates the comparison of that method with the new LPS-ROM here proposed as stabilization-motivated ROM. This theoretical and computational study is today in progress and shall appear in a forthcoming paper.

\medskip

{\sl Acknowledgments:} This work has been partially supported by the Spanish Government-EU Feder grant MTM2015-64577-C2-1-R. The research of the author has been also funded by the Spanish State Research Agency through the national programme Juan de la Cierva-Incorporaci\'on 2017.
The author acknowledges Prof. M. Aza\"iez (Laboratoire I2M CNRS UMR5295, Institut Polytechnique de Bordeaux) for some advices and fruitful suggestions, and Prof. T. Chac\'on (Departamento EDAN \& IMUS, Universidad de Sevilla) for some helpful discussions, especially on Remark \ref{rm:PresEst}.

\bibliographystyle{abbrv}
\bibliography{Biblio_STAB-POD-ROM_NSE}


\end{document}